\documentclass{elsarticle}
\usepackage{a4wide}
\usepackage[english]{babel}
\usepackage{amsmath}
\usepackage{amssymb}
\usepackage{amsfonts}
\usepackage[ruled,linesnumbered]{algorithm2e}
\usepackage{tikz}
\usetikzlibrary{math}
\usetikzlibrary{spy}
\usetikzlibrary{calc}
\usetikzlibrary{patterns}
\usepackage{subcaption}
\usepackage{xcolor}
\biboptions{sort&compress}
\usepackage[perpage]{footmisc}

\newcommand{\appropto}{\mathrel{\vcenter{
  \offinterlineskip\halign{\hfil$##$\cr
    \propto\cr\noalign{\kern2pt}\sim\cr\noalign{\kern-2pt}}}}}

\newtheorem{remark}{Remark}[section]

\journal{Computational Mechanics}

\begin{document}

\begin{frontmatter}

\title{Scalable multigrid methods for immersed finite element methods and immersed isogeometric analysis}

\author[TUE]{F. de Prenter\corref{mycorrespondingauthor}}
\cortext[mycorrespondingauthor]{Corresponding author: fritsdeprenter@gmail.com}

\author[TUE]{C.V.\ Verhoosel}

\author[TUE]{E.H.\ van Brummelen}

\author[UCB]{J.A.\ Evans}

\author[UCB,DLR]{C.\ Messe}

\author[UCB,WD]{J.\ Benzaken}

\author[UCB]{K.\ Maute}

\address[TUE]{Department of Mechanical Engineering, Eindhoven University of Technology, The Netherlands}
\address[UCB]{Ann \& H.J.\ Smead Department of Aerospace Engineering Sciences, University of Colorado Boulder, USA}
\address[DLR]{German Aerospace Center (DLR), Germany}
\address[WD]{Walt Disney Animation Studios, USA}

\begin{abstract}
Ill-conditioning of the system matrix is a well-known complication in immersed finite element methods and trimmed isogeometric analysis. Elements with small intersections with the physical domain yield problematic eigenvalues in the system matrix, which generally degrades efficiency and robustness of iterative solvers. In this contribution we investigate the spectral properties of immersed finite element systems treated by Schwarz-type methods, to establish the suitability of these as smoothers in a multigrid method. Based on this investigation we develop a geometric multigrid preconditioner for immersed finite element methods, which provides mesh-independent and cut-element-independent convergence rates. This preconditioning technique is applicable to higher-order discretizations, and enables solving large-scale immersed systems in parallel, at a computational cost that scales linearly with the number of degrees of freedom. The performance of the preconditioner is demonstrated for conventional Lagrange basis functions and for isogeometric discretizations with both uniform B-splines and locally refined approximations based on truncated hierarchical B-splines.

\end{abstract}


\end{frontmatter}

\section{Introduction}\label{sec:introduction}

Immersed methods are useful tools to avoid laborious and computationally expensive procedures for the generation of body-fitted finite element discretizations or analysis-suitable NURBS geometries in isogeometric analysis, specifically for problems on complex, moving, or implicitly defined geometries. Immersed finite element techniques, such as the finite cell method \cite{Parvizian2007,Duester2008,Schillinger2015}, CutFEM \cite{BurmanHansbo2012,Burman2015} and immersogeometric analysis \cite{Kamensky2015,Hsu2015}, have been successfully applied to a broad range of problems. Noteworthy applications include isogeometric analysis on trimmed CAD objects, e.g., \cite{Schmidt2012,Rank2012,Schillinger2012,Ruess2013,Ruess2014,Marussig2017review}, fluid-structure interaction with large displacements, e.g., \cite{Bazilevs2012,Hsu2014,Burman2014,Massing2015,Kadapa2016,Kadapa2017,Wu2017,Kadapa2018}, scan based analysis \cite{Ruess2012,Verhoosel2015,Elhaddad2017,Duczek2015pore,Wurkner2018,Hoang2019} and topology optimization, e.g., \cite{Parvizian2012,Dijk2013,Nadal2013,Bandara2016,Groen2017,Villanueva2017,Burman2018}.

An essential aspect of finite element methods and isogeometric analysis is the computation of the solution to a system of equations. This is specifically challenging for systems derived from immersed methods, since such methods generally yield severely ill-conditioned system matrices \cite{SIPIC}. For this reason, many researchers resort to direct solvers, e.g., \cite{Duester2008,Schillinger2012,Rank2012,Ruess2013,Schillinger2015,Ruess2014,Elhaddad2017}, the efficiency of which is not affected by the conditioning of the system matrix. Nevertheless, the computational cost of direct solvers, both in terms of memory and floating point operations, scales poorly with the size of the system. With iterative solvers, the scaling between the computational cost and the size of the system is generally better, making these more suitable for large systems of equations \cite{Barrett}. However, the efficiency and reliability of iterative solution methods depends on the conditioning of the system. Without dedicated treatments, the severe ill-conditioning of linear systems derived from immersed finite element methods generally forestalls convergence of iterative solution procedures. Multiple resolutions for these conditioning problems have been proposed, the most prominent of which are the ghost penalty, e.g., \cite{Burman2010,BurmanHansbo2012,Burman2015}, constraining, extending, or aggregation of basis functions, e.g., \cite{Hoellig2001,Hoellig2005,Rueberg2012,Rueberg2014,Rueberg2016,Marussig2017,Marussig2017review,Badia2018,Badia2018b,Marussig2018}, and preconditioning, which is discussed in detail below.

Several dedicated preconditioners have been developed for immersed finite element methods. It is demonstrated in \cite{Lang2014} that a diagonal preconditioner in combination with the constraining of very small basis functions results in an effective treatment for systems with linear bases. With certain restrictions to the cut-element-geometry, \cite{Lehrenfeld2017} derives that a scalable preconditioner for linear bases is obtained by combining diagonal scaling of basis functions on cut elements with standard multigrid techniques for the remaining basis functions. %It is derived in \cite{Lehrenfeld2017} that the diagonal scaling of trimmed basis functions is effective for linear bases and certain restrictions to how elements can be cut.
In \cite{Badia2017} the scaling of a BDDC (Balancing Domain Decomposition by Constraints) is tailored to cut elements, and this is demonstrated to be effective with linear basis functions. An algebraic preconditioning technique is presented in \cite{SIPIC}, which results in an effective treatment for smooth function spaces. References \cite{Prenter2019} and \cite{John} establish that additive Schwarz preconditioners can effectively resolve the conditioning problems of immersed finite element methods with higher-order discretizations, for both isogeometric and $hp$-finite element function spaces and for both symmetric positive definite (SPD) and non-SPD problems. Furthermore, the numerical investigation in \cite{Prenter2019} conveys that the conditioning of immersed systems treated by an additive Schwarz preconditioner is very similar to that of mesh-fitting systems. In particular, the condition number of an additive Schwarz preconditioned immersed system exhibits the same mesh-size dependence as mesh-fitting approaches \cite{Saad,Johnson}, which opens the doors to the application of established concepts of multigrid preconditioning. It should be mentioned that similar conditioning problems as in immersed methods occur in XFEM and GFEM. Dedicated preconditioners have been developed for these problems as well, a survey of which can be found in \cite{Prenter2019}.

Multigrid methods effectively resolve the mesh-size dependence of the conditioning of linear systems and its effect on the convergence of iterative solution methods. In particular, the use of overlapping Schwarz smoothers can lead to multigrid methods with provably mesh-independent convergence rates \cite{Arnold2000,Schoeberl2003,Smith1996}. There exists a rich literature on multigrid techniques, and interested readers are directed to the reference works \cite{Hackbusch1985,Wesseling1992,Briggs2000,Brandt2011}. Multigrid methods have not been studied extensively in the context of immersed finite element methods, but detailed studies regarding closely related aspects are available. In isogeometric analysis, multigrid is an established concept, e.g., \cite{Gahalaut2013,Buffa2013,Donatelli2015,Donatelli2017,Hofreither2017,Hofreither2017correction,Takacs2018,Sogn2019}. In regard of the present manuscript \cite{BeiraoDaVeiga2012,BeiraoDaVeiga2013Schwarz,Coley2018,delaRiva2018} are particularly noteworthy, as these all employ smoothers that are based on Schwarz-type techniques. Another interesting contribution is \cite{Hofreither2016}, which presents a multigrid technique for locally refined function spaces with truncated hierarchical B-splines, that are also employed in this manuscript. Further noteworthy references are multigrid preconditioners for XFEM \cite{Berger-Vergiat2012,Hiriyur2012}, discontinuous Galerkin (DG) \cite{Engwer2016}, unfitted interface problems \cite{Ludescher2018}, and an algebraic multigrid (AMG) preconditioner that is applied to immersed systems which have been treated by an aggregation procedure \cite{Verdugo2019}.

The main objective of this contribution is to develop a geometric multigrid preconditioning technique that is applicable to higher-order immersed finite element methods with conventional, isogeometric, and locally refined basis functions. This preconditioner enables iterative solution methods with a convergence rate that is unaffected by either the cut elements or the grid size. The intrinsic dependence on mesh regularity in geometric multigrid approaches is non-restrictive for immersed finite element methods, as these methods generally employ structured grids. Based on the observations regarding the mesh-size dependence of the additive Schwarz preconditioner developed in \cite{Prenter2019}, this contribution further investigates the spectral properties of immersed systems treated with Schwarz-type preconditioners, in order to establish the suitability of these as smoothers in a multigrid method. This results in a preconditioning technique with the desired properties, the performance of which is demonstrated on a range of test cases with multi-million degrees of freedom. This numerical investigation includes the application to locally refined bases with truncated hierarchical B-splines, and involves a detailed study of aspects that are specific for immersed finite elements.

This contribution only considers SPD problems. Multigrid methods are an established concept in fluid mechanics as well \cite{Brandt2011}, however, and similar Schwarz-type methods have successfully been applied to flow problems with both immersed finite element methods \cite{Prenter2019} and mesh-fitting multigrid solvers \cite{Coley2018}, i.e., Vanka-smoothers \cite{Vanka1986}. Therefore, it is anticipated that the presented preconditioning technique extends matatis mutandis to non-SPD and mixed formulations.

In Section~\ref{sec:methodMulti} of this contribution the employed immersed finite element method is presented, including the quadrature on cut elements and the applied discretization spaces. Section~\ref{sec:multigrid} investigates spectral properties of immersed finite element methods and presents the developed geometric multigrid preconditioner. In Section~\ref{sec:resultsMulti} this preconditioning technique is assessed on a range of test cases and conclusions are drawn in Section~\ref{sec:conclusionMulti}.

\section{Immersed finite element formulation}\label{sec:methodMulti}

\begin{figure}[pt]
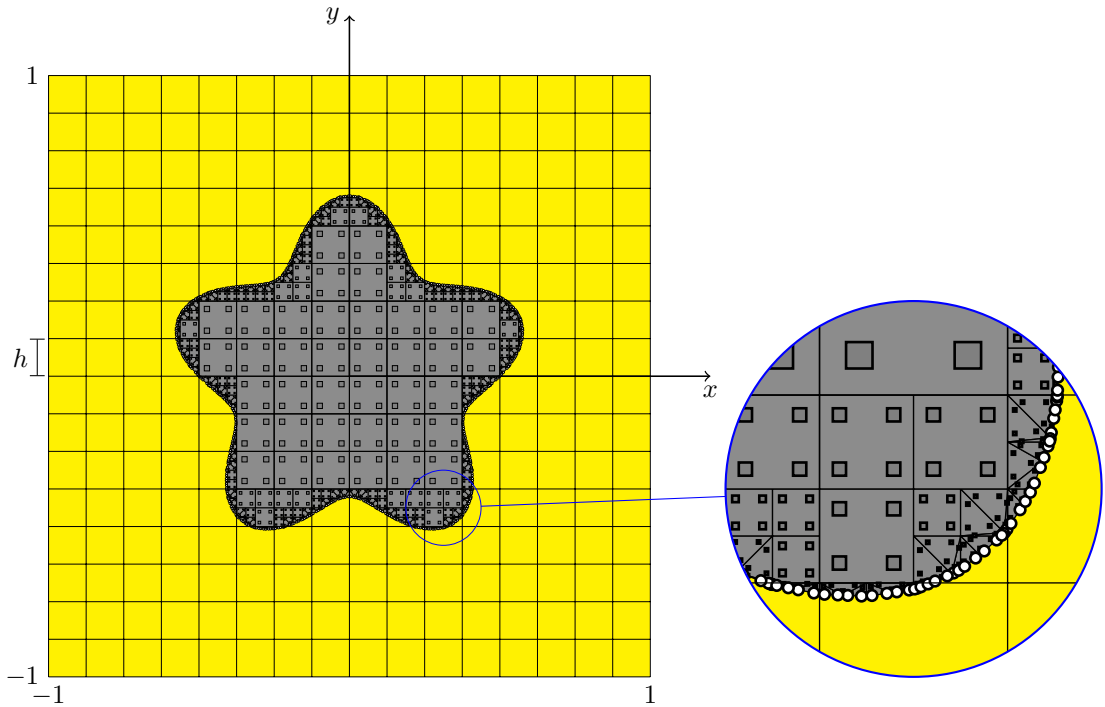

  \centering
  \begin{tikzpicture}
  \tikzmath{\wambient=4;
            \h=0.125*\wambient;
            }

  \begin{scope}[spy using outlines={circle, magnification=5, size=10*\h cm, connect spies, blue}]

    % embedding domain 
    \fill[yellow, fill opacity=1] (-\wambient,-\wambient) rectangle (\wambient,\wambient);
    
    % subcells (from data)
    \begin{scope}[scale=\wambient,line join=bevel,ultra thin,fill=black!45,draw=black]
      \input{integration_subcells.dat}
    \end{scope}

%     % elements (from data)
%     \begin{scope}[scale=\wambient,line join=bevel]
%       \input{integration_elements.dat}
%     \end{scope}

    % ipoints interior (from data)
    \begin{scope}[scale=\wambient,line width=0.2pt,inner sep=0pt,minimum size=1pt,rectangle,draw=black!100,fill=black!50]
      \input{integration_points.dat}
    \end{scope}

    % ipoints boundary (from data)
     \begin{scope}[scale=\wambient,line width=0.2pt,inner sep=0pt,minimum size=0.8pt,circle,draw=black!100,fill=white!100]
      \input{integration_boundary.dat}
     \end{scope}

    % grid
    \begin{scope}
      \draw[step=\h,black,ultra thin] (-\wambient,-\wambient) grid (\wambient,\wambient);
      \draw[|-|] (-\wambient-0.3*\h,0*\h) -- node[left]{$h$} ++(0,\h);
      \node[left] at (-\wambient,-\wambient) {$-1$};
      \node[left] at (-\wambient,\wambient) {$1$};
      \node[below] at (-\wambient,-\wambient) {$-1$};
      \node[below] at (\wambient,-\wambient) {$1$};
    \end{scope}

    % axes
    \draw[->,semithick] (0,0) -- (1.2*\wambient,0) node[below] {$x$} coordinate(x axis);
    \draw[->,semithick] (0,0) -- (0,1.2*\wambient) node[left] {$y$} coordinate(y axis);
    
    % integration point zoom
    \spy on ({(2.5*\h},{(-3.5*\h}) in node(zoom) at (15*\h,-3*\h);

  \end{scope}

\end{tikzpicture}
  \caption{The physical domain $\Omega$ defined by the level set function $\psi(x,y) = 0.5 + 0.1 {\rm sin}(5\theta) - r(x,y)$, with $\theta = {\rm arctan2}(y,x)$, and $r(x,y)^2 = x^2 + y^2$ (inspired by \cite{Lehrenfeld2016}). The physical domain is encapsulated by the embedding domain $(-1,1)^2$ on which the grid with grid size $h=\frac18$ is posed. Quadrature is performed by the integration procedure that is outlined in \cite{Verhoosel2015}. This procedure recursively bisects cut elements, until a maximum integration depth is reached. The bisected elements are then triangulated, to obtain integration subcells. Standard Gaussian integration points are applied on these subcells to evaluate volumetric integrals (gray squares). Boundary integrals are computed by Gaussian quadrature on the approximate boundary, formed by the edges of the integration subcells (white circles).\label{fig:domainStar}}
\end{figure}

We consider problems on a two-dimensional or three-dimensional domain $\Omega\subset\mathbb{R}^d$ ($d\in\{2,3\}$), that is referred to as the physical domain. The physical domain is encapsulated by a fictitious extension, to obtain an embedding domain of simple shape, as illustrated in Figure~\ref{fig:domainStar}. Because of its simple shape, it is trivial to generate a structured, e.g., tensor product, mesh for the embedding domain. These structured meshes render immersed finite elements ideally suitable for geometric multigrid approaches. The set of elements that intersects the physical domain is referred to as the background grid, and can serve as a substructure to construct different types of basis functions. We denote this mesh by $\mathcal{T}_h$, where $h$ refers to the grid size, and the basis functions that are supported on the physical domain span the approximation space $\mathcal{V}_h$. A non-trivial aspect of immersed finite element methods is the integration over cut elements. Herein we employ the procedure as outlined in \cite{Verhoosel2015}, which is illustrated in Figure~\ref{fig:domainStar}.

The numerical examples in this contribution consider problems in linear elasticity:%pertain to the linear elasticity problem:
\begin{equation}\left\{ \begin{array}{c l}
  \mbox{div}\left( \boldsymbol{\sigma} \right) + \boldsymbol{f} = \boldsymbol{0} & \mbox{in } \Omega, \\
  \boldsymbol{u} = \boldsymbol{g}^D & \mbox{on } \Gamma^D, \\
  \boldsymbol{\sigma} \boldsymbol{n} = \boldsymbol{g}^N & \mbox{on } \Gamma^N,                         
 \end{array}\right.\label{eq:strongMulti}\end{equation}
with Cauchy stress tensor $\boldsymbol{\sigma} = \boldsymbol{\sigma}(\boldsymbol{u}) =\lambda \mbox{div}\left(\boldsymbol{u}\right)\mathbf{I} + 2 \mu \nabla^s \boldsymbol{u}$, Lam\'{e} parameters $\lambda$ and $\mu$, $\nabla^s$ denoting the symmetric gradient operator, $\mathbf{n}$ the exterior unit normal vector and $\Gamma^D \cup \Gamma^N = \partial \Omega$ complementary parts of the boundary on which Dirichlet and Neumann conditions are prescribed. The formulations in this section are restricted to pure Dirichlet or Neumann boundary conditions, but can easily be modified to mixed boundary conditions with a prescribed normal displacement and tangential traction or Robin-type conditions. We consider approximations of \eqref{eq:strongMulti} based on a symmetric and coercive variational form:
\begin{equation}\left\{ \begin{array}{l}
\mbox{Find } \boldsymbol{u}_h \in \mathcal{V}_h \mbox{ such that for all } \boldsymbol{v}_h \in \mathcal{V}_h\mbox{:} \\[.5em]
\qquad a_h\left(\boldsymbol{v}_h,\boldsymbol{u}_h\right) = b_h\left(v_h\right),
\end{array}\right.
\label{eq:weakMulti}
\end{equation}
with the bilinear and linear operators defined as:
\begin{equation}\begin{aligned}
a_h\left(\boldsymbol{v}_h,\boldsymbol{u}_h\right) & = \int_\Omega \nabla^s \boldsymbol{v}_h : \boldsymbol{\sigma}(\boldsymbol{u}_h) \,\mbox{d}V + \int_{\Gamma^D} \left( \lambda \beta_h^\lambda (\boldsymbol{v}_h \cdot \boldsymbol{n}) (\boldsymbol{u}_h \cdot \boldsymbol{n}) + 2 \mu \beta_h^\mu \boldsymbol{v}_h \cdot \boldsymbol{u}_h \right) \,\mbox{d}S, \\
b_h\left(\boldsymbol{v}_h\right) &= \int_\Omega \boldsymbol{v}_h \cdot \mathbf{f} \,\mbox{d}V + \int_{\Gamma^D} \left( \lambda \beta_h^\lambda (\boldsymbol{v}_h \cdot \boldsymbol{n}) (\boldsymbol{g}^D \cdot \boldsymbol{n}) + 2 \mu \beta_h^\mu \boldsymbol{v}_h \cdot \boldsymbol{g}^D \right) \,\mbox{d}S + \int_{\Gamma^N} \mathbf{v}_h \cdot \mathbf{g}^N \,\mbox{d}S.
\end{aligned}\label{eq:operatorsMulti}\end{equation}
The weak formulation in \eqref{eq:weakMulti} and \eqref{eq:operatorsMulti} imposes the Dirichlet conditions by the penalty method with parameters $\beta_h^\lambda$ and $\beta_h^\mu$ that are chosen inversely proportional to the grid size. The analysis of the conditioning of immersed finite elements in \cite{SIPIC} conveys that cut-element-specific conditioning problems are not essentially affected by the type of weak enforcement of the boundary conditions, provided that with a Nitsche-type enforcement of Dirichlet conditions \cite{Nitsche1971} coercivity of the variational form is retained \cite{Embar2010}. Therefore all computations in this contribution apply the penalty method, which bypasses the computation of local stabilization parameters. It should be noted that, due to the penalty parameters, the operators $a_h\left(\cdot,\cdot\right)$ and $b_h\left(\cdot\right)$ are mesh-size dependent. This is relevant in regard of multigrid techniques, in which the same problem is considered on different meshes. The coarse problems in this contribution inherit the weak formulation from the finest grid, without modification of the penalty parameters, see Remark~\ref{rem:gridParamCorr} for details.

\begin{figure}[pt]
  \centering
  \begin{subfigure}{0.4\textwidth}
   \begin{center}
    \includegraphics[height=4.5cm]{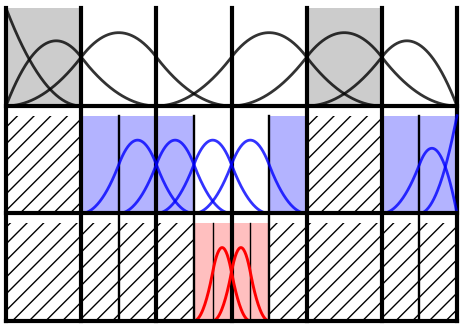}
    \caption{Hierarchical (non-truncated) B-splines\label{fig:HB_basis}}
   \end{center}
  \end{subfigure}
  \hfill
  \begin{subfigure}{0.4\textwidth}
   \begin{center}
    \includegraphics[height=4.5cm]{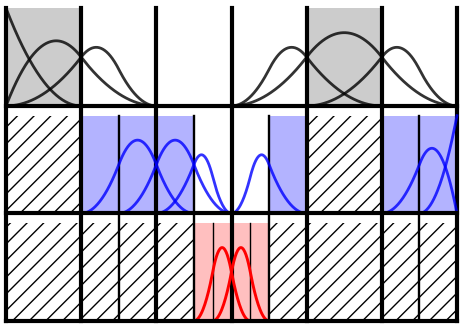}
    \caption{Truncated hierarchical B-splines\label{fig:THB_basis}}
   \end{center}
  \end{subfigure}
  \hspace{.2cm}
  \begin{minipage}{.08\textwidth}
   \begin{center}
    $m=0$ \\
    \vspace{1cm}
    $m=1$ \\
    \vspace{1cm}
    $m=2$ \\
    \vspace{.6cm}
   \end{center}
  \end{minipage}
  \caption{Illustration of an open, hierarchically refined, quadratic B-spline basis. The local refinement level is denoted by $m \leq M$, with $M$ the number of local refinement levels. Active elements $K_i^m$ of level $m\in\{0,1,2\}$ are indicated in gray, blue, and red, respectively. The set of all active elements of level $m$ is denoted by $\mathcal{T}_h^m$, and the union of these sets forms the hierarchically refined mesh $\mathcal{T}_h = \cup_{m=0}^{m=M} \mathcal{T}_h^m$. Refined elements, i.e., that have active elements at a finer level underneath, are shown blank and inactive elements, i.e., that have an active element at a coarser level above, are hatched. Figure~\ref{fig:HB_basis} shows the active standard (non-truncated) hierarchical B-spline basis functions. A basis function of level $m$ is active under the condition that: \emph{i)} it is supported on at least one active element of level $m$, and \emph{ii)} it is not supported on inactive elements of level $m$, i.e., the entire support is of local refinement level $\geq m$. Figure~\ref{fig:THB_basis} shows the truncated basis, which is obtained by truncating the basis functions with respect to active basis functions of finer levels.\label{fig:THB}}
\end{figure}

An advantageous property of immersed methods is that different approximation spaces can relatively easily be employed, by virtue of the regularity of the underlying mesh. We herein consider uniform discretizations with both traditional Lagrange basis functions and B-splines \cite{Prenter1975}, and locally refined discretizations with truncated hierarchical B-splines (THB-splines) \cite{Gianelli2012}. The construction of THB-splines is illustrated in Figure~\ref{fig:THB}. THB-splines are posed on a hierarchy of meshes, each of which has a certain number of active B-splines. The truncated basis is then obtained by truncating active B-splines with respect to the active B-splines on the finer grids in which they are nested. We refer the reader to \cite{Gianelli2012,Bracco2019} for details regarding the construction of THB-splines. Truncating the basis functions reduces the computational cost, because the truncated basis functions have smaller supports than the standard (non-truncated) hierarchical basis functions. This results in a sparser system matrix. In the context of the multigrid preconditioner developed herein, THB-splines have a particularly important advantage concerning the computational cost. As will be elaborated in Section~\ref{sec:smoothers}, the block selection for the Schwarz-type smoother in the developed preconditioner yields smaller blocks with THB-splines than it would with a non-truncated basis, and results in a nearly diagonal treatment of untrimmed basis functions. Furthermore, it is demonstrated in \cite{Hofreither2016} that mesh-fitting systems with truncated bases are generally better conditioned than systems with non-truncated bases, and that multigrid methods with a diagonal smoother are effective for systems with THB-splines. 

\section{Multigrid methods for immersed finite element methods}\label{sec:multigrid}

This section presents the developed geometric multigrid preconditioner for immersed finite element methods. Section~\ref{sec:conditioningIntro} discusses the conditioning effects of both cut elements and the mesh size of the background grid. In Section~\ref{sec:algorithm} the employed preconditioning algorithm is introduced. Section~\ref{sec:coarsening} considers the prolongation and restriction operators, specifically in the context of locally refined grids, and effective smoothers for immersed methods are discussed in Section~\ref{sec:smoothers}.

\subsection{Conditioning aspects of immersed finite element methods}\label{sec:conditioningIntro}

Introducing the basis $\{\boldsymbol{\phi}_i\}_{i=1}^n$ for the approximation space $\mathcal{V}_h$, the variational formulation in \eqref{eq:weakMulti} leads to the linear system:
\begin{equation}
 \mathbf{A}\mathbf{x} = \mathbf{b},
\label{eq:systemMulti}\end{equation}
with symmetric positive definite (SPD) matrix $\mathbf{A}$, with $A_{ij} = a_h\left( \boldsymbol{\phi}_i, \boldsymbol{\phi}_j \right)$, solution vector $\mathbf{x}$, such that $\mathbf{u}_h = \sum_{i=1}^n x_i \boldsymbol{\phi}_i$, and right hand side vector $\mathbf{b}$, with $b_i = b_h\left(\boldsymbol{\phi}_i\right)$. Immersed finite element methods generally lead to systems of equations that are fundamentally difficult to solve. Without dedicated treatment of the cut-element-specific ill-conditioning of systems derived from immersed finite element methods, in general the convergence of iterative solvers is severely retarded \cite{SIPIC,Prenter2019,John}. An important indicator of the feasibility of iterative solution procedures for a linear system as in \eqref{eq:systemMulti} is the condition number, $\kappa\left(\mathbf{A}\right)$. For the symmetric positive definite (SPD) systems considered in this contribution, the condition number is equal to the quotient of the largest to the smallest eigenvalue. It should be noted that the convergence of iterative solution methods is not merely dependent on the condition number, but rather depends on the entire spectrum of the system matrix, i.e., on the distribution of the complete set of eigenvalues. Systems with well-clustered spectra generally lead to faster convergence than systems with eigenvalues that are spread out. For a detailed discussion about the aspects affecting the performance of iterative solvers, the reader is directed to \cite{Saad,Greenbaum1997}.

\begin{figure}[pt]
   \begin{center}
   \begin{subfigure}{.49\textwidth}
   \includegraphics[height=5.1cm]{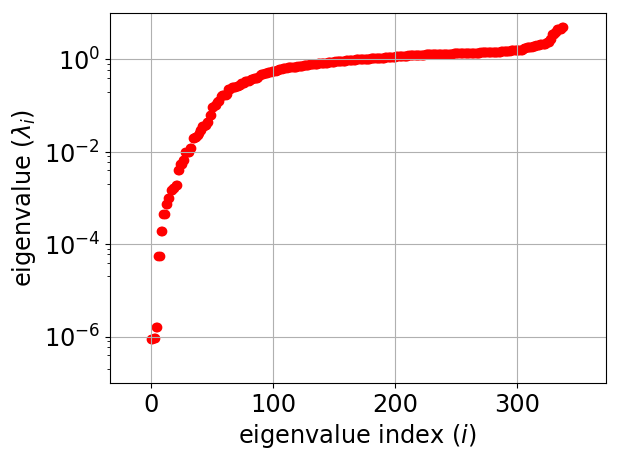}
   \caption{Eigenvalue spectrum\label{fig:JacobiSpectrum}}
  \end{subfigure}
  \end{center}
  \begin{subfigure}{.49\textwidth}
   \includegraphics[height=5cm]{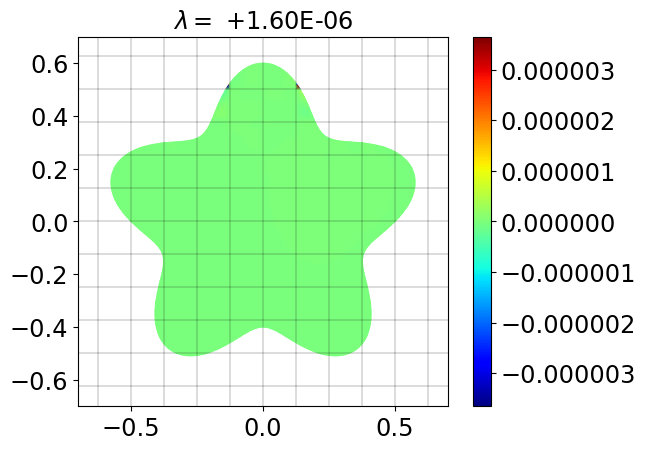}
   \caption{Smallest mode\label{fig:JacobiSmall}}
  \end{subfigure}
  \hspace{.5cm}
  \begin{subfigure}{.49\textwidth}
   \includegraphics[height=5cm]{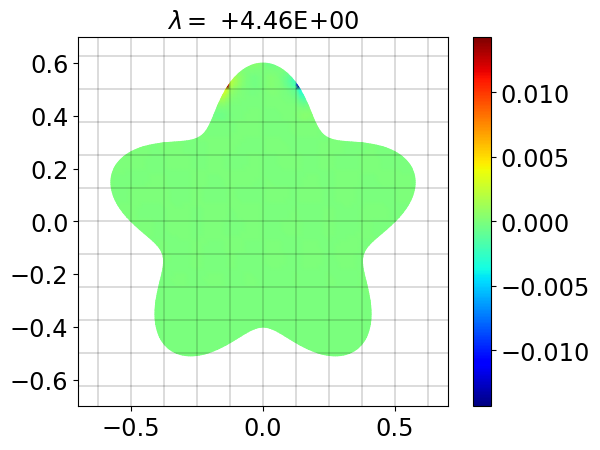}
   \caption{Largest mode\label{fig:JacobiLarge}}
  \end{subfigure}
  \begin{subfigure}{.49\textwidth}
   \vspace{.2cm}
   \includegraphics[height=5cm]{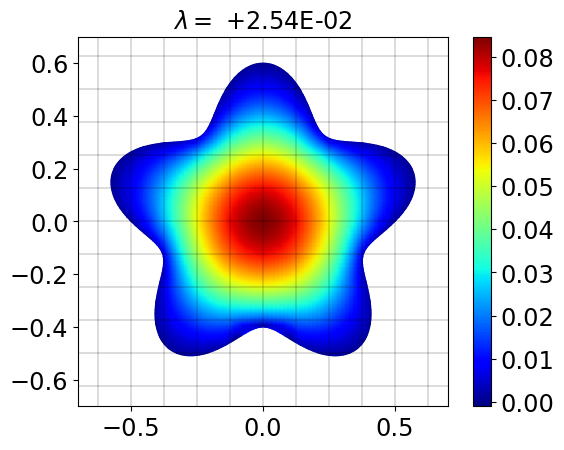}
   \caption{Smooth mode\label{fig:JacobiSmooth}}
  \end{subfigure}
  \hspace{.5cm}
  \begin{subfigure}{.49\textwidth}
   \vspace{.2cm}
   \includegraphics[height=5cm]{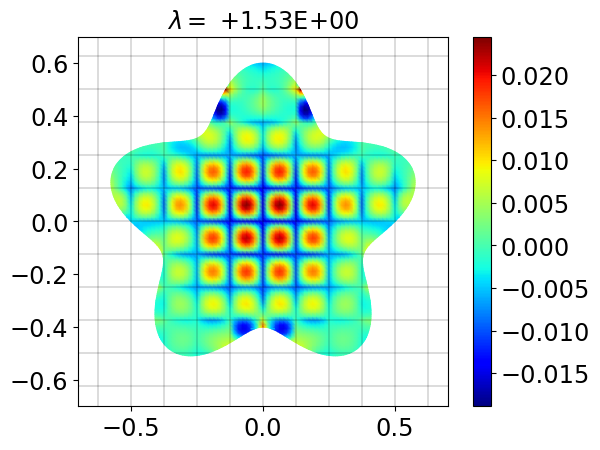}
   \caption{Oscillatory mode\label{fig:JacobiOscil}}
  \end{subfigure}
 \caption{Typical spectrum and characteristic eigenmodes of an immersed system that is preconditioned with Jacobi.\label{fig:Jacobi}}
\end{figure}

%Because all eigenvalues of SPD systems are real and positive, the spectra of the systems considered in this contribution can be visualized as in Figure~\ref{fig:JacobiSpectrum}. In this figure, the horizontal axis corresponds to the indices of eigenmodes, and the vertical axis corresponds to the eigenvalues. 
%The spectrum in Figure~\ref{fig:JacobiSpectrum} corresponds to a linear system derived from the Poisson problem on the geometry in Figure~\ref{fig:domainStar}, preconditioned with a Jacobi preconditioner. This system employs quadratic Lagrange basis functions on a mesh with grid size $h=\frac18$, and imposes Dirichlet conditions on the entire boundary by a penalty method with parameter $\frac2h$. Characteristic eigenfunctions corresponding to eigenvalues in this spectrum are illustrated in Figures~\ref{fig:JacobiSmall}--\ref{fig:JacobiOscil}.

To elucidate the general characteristics of spectra emanating from immersed finite element methods, Figure~\ref{fig:Jacobi} displays the spectrum of an immersed approximation of the Laplace operator on the star-shaped domain from Figure~\ref{fig:domainStar}, as well as 4 characteristic eigenmodes. The immersed approximation space is composed of quadratic Lagrange basis functions on a mesh with $h=\frac18$. Dirichlet boundary conditions are imposed by means of a penalty method with parameter $\frac2h$. Because symmetric positive definite (SPD) systems only posses real and positive eigenvalues, their spectra can be conveniently represented by plotting the eigenvalues $\lambda_i\in\mathbb{R}^+$ \emph{vs.}\ their index $i\in\mathbb{N}$, as in Figure~\ref{fig:JacobiSpectrum}. Figure~\ref{fig:JacobiSmall} displays a very small eigenmode that is only supported on a small cut element, which is exemplary for eigenmodes that are common in immersed finite element methods. As described in detail in \cite{SIPIC}, basis functions on small cut elements can become almost linearly dependent, which yields very small eigenvalues and is not repaired by Jacobi preconditioning. It can be observed that this eigenmode is very similar to the largest eigenmode of the system, plotted in Figure~\ref{fig:JacobiLarge}, as both eigenmodes consist of essentially the same basis functions. However, while the almost linearly dependent basis functions are subtracted from each other such that these cancel out in Figure~\ref{fig:JacobiSmall}, the opposite happens in Figure~\ref{fig:JacobiLarge}. Note that with the quadratic Lagrange basis, there are 4 almost linearly dependent basis functions that are only supported on the small cut elements, such that the largest eigenvalue is bounded from below by approximately $4$. Based on an analysis of the typical small eigenmodes in immersed finite element methods, an estimate of the condition number of SPD systems for second-order PDEs is derived in \cite{SIPIC}:
\begin{equation}
 \kappa\left(\mathbf{A}\right) = \mathcal{O} \left( \eta^{-\left(2p+1-2/d\right)} \right),
\end{equation}
with $p$ the polynomial degree of the approximation space, $d$ the number of dimensions and $\eta$ the smallest volume fraction, defined as the smallest relative intersection of an element with the physical domain:
\begin{equation}
 \eta = \min_{K_i \in \mathcal{T}_h | K_i \cap \Omega \neq \emptyset} \frac{|K_i \cap \Omega|}{|K_i|}.
\end{equation}
Since cut elements can be arbitrarily small, systems derived from immersed finite element formulations can be arbitrarily ill-conditioned. As a result, iterative solvers are generally ineffective in case that no dedicated treatment for the cut-element-induced conditioning problem is applied, see e.g., the introduction in Section~\ref{sec:introduction}. Figures~\ref{fig:JacobiSmooth} and \ref{fig:JacobiOscil} portray characteristic eigenfunctions that are not exclusive to immersed finite element methods. The smooth eigenmode in Figure~\ref{fig:JacobiSmooth} is very similar to the usual smallest eigenmode in mesh-fitting finite elements, and is in close correspondence with the smallest analytical eigenmode of the considered PDE. Figure~\ref{fig:JacobiOscil} plots the oscillatory eigenfunction with the highest frequency that can be captured by the grid. This eigenmode is similar to the usual largest eigenmode in mesh-fitting systems. The ratio between the largest and smallest eigenvalue in mesh-fitting systems, and with that the condition number, is therefore mesh-size dependent. It can be shown that for second-order PDEs it holds that \cite{Johnson}:
\begin{equation}
 \kappa\left(\mathbf{A}\right) = \mathcal{O}\left( h^{-2} \right).
\label{eq:kappaMeshSize}\end{equation}
This deteriorates the conditioning for very fine meshes, which also retards the convergence of iterative solvers. In particular, smooth eigenmodes as in Figure~\ref{fig:JacobiSmooth} -- which yield small eigenvalues of $\mathcal{O}(h^2)$ with Jacobi preconditioning -- converge slowly. For fixed point iteration methods a convergence rate of $1-\mathcal{O}(h^2)$ can be derived \cite{Briggs2000}, and also plain Krylov subspace methods generally require $\mathcal{O}\left( h^{-1} \right)$ iterations \cite{Saad}.

Multigrid methods are often applied to resolve the mesh-dependence of the condition number according to \eqref{eq:kappaMeshSize} and the corresponding slow convergence, and provide a conditioning and convergence rate that is independent of the mesh size. %A common method to resolve the mesh-dependence of the condition number according to \eqref{eq:kappaMeshSize}, and the corresponding slow convergence, is multigrid (preconditioning), which provides a convergence rate (conditioning) that is independent of the mesh size.
Due to the smoothness of small eigenmodes in mesh-fitting systems, these modes can be effectively approximated on a coarser grid. The concept of multigrid methods is to combine iterations on a fine grid -- by which the oscillatory modes with large eigenvalues as in Figure~\ref{fig:JacobiOscil} converge quickly -- with coarse grid corrections in which the approximation of the solution is amended by the solution on a coarser grid. These coarse grid corrections enhance the convergence of the smooth eigenmodes as in Figure~\ref{fig:JacobiSmooth}, such that a convergence behavior is obtained that is independent of the grid size. The reader is directed to \cite{Hackbusch1985,Wesseling1992,Briggs2000,Brandt2011} for reference works on multigrid techniques. In \cite{Prenter2019} an additive Schwarz preconditioner is presented that is tailored to immersed finite element methods, and resolves the conditioning problems related to cut elements. Furthermore, it is observed that the systems treated with this preconditioner behave similarly with respect to the grid size as mesh-fitting systems, both in terms of the condition number and in terms of the number of iterations with Krylov subspace methods. This results in a computational cost for solving immersed systems that -- while scaling better with the size of the system than direct solvers -- is suboptimal, in the sense that it is not yet linear with the number of degrees of freedom. In the following sections we therefore incorporate aspects of this preconditioner in a multigrid framework, to obtain a solution method that is robust to cut elements and independent of the size of the system.

\subsection{Multigrid V-cycle algorithm}\label{sec:algorithm}

\begin{figure}[pt]
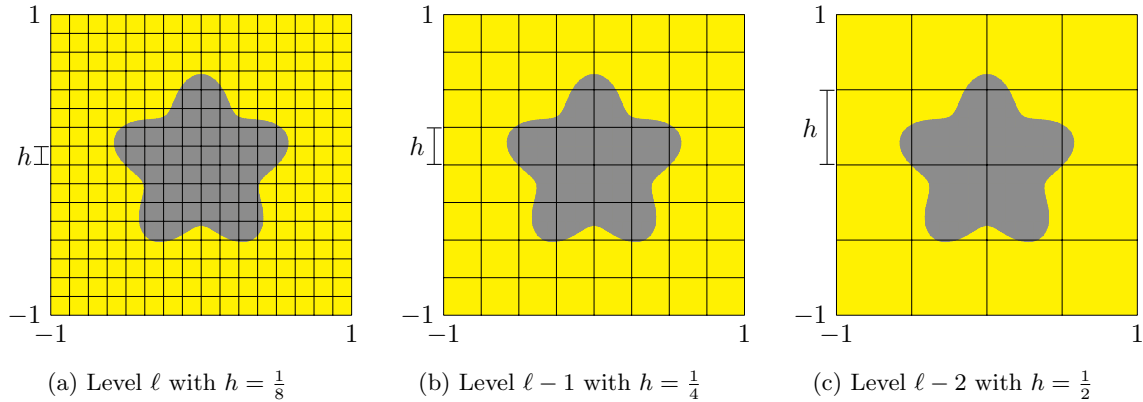

 \begin{center}
  \begin{subfigure}{.3\textwidth}
   \begin{center}
    \begin{tikzpicture}
  \tikzmath{\wambient=2;
            \h=0.125*\wambient;
            }

  \begin{scope}

    % embedding domain 
    \fill[yellow, fill opacity=1] (-\wambient,-\wambient) rectangle (\wambient,\wambient);
    
    % subcells (from data)
    \begin{scope}[scale=\wambient]
      \input{integration_subcells_Frits.dat}
    \end{scope}

    % grid
    \begin{scope}
      \draw[step=\h,black,ultra thin] (-\wambient,-\wambient) grid (\wambient,\wambient);
      \draw[|-|] (-\wambient-.5*\h,0*\h) -- node[left]{$h$} ++(0,\h);
      \node[left] at (-\wambient,-\wambient) {$-1$};
      \node[left] at (-\wambient,\wambient) {$1$};
      \node[below] at (-\wambient,-\wambient) {$-1$};
      \node[below] at (\wambient,-\wambient) {$1$};
    \end{scope}

  \end{scope}

\end{tikzpicture}
    \vspace{-.4cm}
    \caption{Level $\ell$ with $h=\frac18$\label{fig:domain8}}
   \end{center}
  \end{subfigure}
  \hfill
  \begin{subfigure}{.3\textwidth}
   \begin{center}
    \begin{tikzpicture}
  \tikzmath{\wambient=2;
            \h=0.25*\wambient;
            }

  \begin{scope}

    % embedding domain 
    \fill[yellow, fill opacity=1] (-\wambient,-\wambient) rectangle (\wambient,\wambient);
    
    % subcells (from data)
    \begin{scope}[scale=\wambient]
      \input{integration_subcells_Frits.dat}
    \end{scope}

    % grid
    \begin{scope}
      \draw[step=\h,black,ultra thin] (-\wambient,-\wambient) grid (\wambient,\wambient);
      \draw[|-|] (-\wambient-0.25*\h,0*\h) -- node[left]{$h$} ++(0,\h);
      \node[left] at (-\wambient,-\wambient) {$-1$};
      \node[left] at (-\wambient,\wambient) {$1$};
      \node[below] at (-\wambient,-\wambient) {$-1$};
      \node[below] at (\wambient,-\wambient) {$1$};
    \end{scope}

  \end{scope}

\end{tikzpicture}
    \vspace{-.4cm}
    \caption{Level $\ell-1$ with $h=\frac14$\label{fig:domain4}}
   \end{center}
  \end{subfigure}
  \hfill
  \begin{subfigure}{.3\textwidth}
   \begin{center}
    \begin{tikzpicture}
  \tikzmath{\wambient=2;
            \h=0.5*\wambient;
            }

  \begin{scope}

    % embedding domain 
    \fill[yellow, fill opacity=1] (-\wambient,-\wambient) rectangle (\wambient,\wambient);
    
    % subcells (from data)
    \begin{scope}[scale=\wambient]
      \input{integration_subcells_Frits.dat}
    \end{scope}

    % grid
    \begin{scope}
      \draw[step=\h,black,ultra thin] (-\wambient,-\wambient) grid (\wambient,\wambient);
      \draw[|-|] (-\wambient-.125*\h,0*\h) -- node[left]{$h$} ++(0,\h);
      \node[left] at (-\wambient,-\wambient) {$-1$};
      \node[left] at (-\wambient,\wambient) {$1$};
      \node[below] at (-\wambient,-\wambient) {$-1$};
      \node[below] at (\wambient,-\wambient) {$1$};
      
    \end{scope}

  \end{scope}

\end{tikzpicture}
    \vspace{-.4cm}
    \caption{Level $\ell-2$ with $h=\frac12$\label{fig:domain2}}
   \end{center}
  \end{subfigure}
  \caption{A hierarchy of nested discretizations on which the problem in \eqref{eq:weakMulti} can be solved.\label{fig:coarsening}}
 \end{center}
\end{figure}

We present the geometric multigrid method in the context of the correction scheme (CS), but the analyses, algorithms and results extend mutatis mutandis to the full approximation scheme (FAS). We consider an algebraic system of the form \eqref{eq:systemMulti}. Let $\mathbf{x}'$ denote an approximation to the solution of \eqref{eq:systemMulti} and let $\mathbf{r}$ denote the corresponding residual according to:
\begin{equation}
 \mathbf{r} = \mathbf{b} - \mathbf{A}\mathbf{x}'.
\label{eq:residualMulti}\end{equation}
To obtain the solution to \eqref{eq:systemMulti}, the approximation $\mathbf{x}'$ must be corrected as $\mathbf{x}'+\bar{\mathbf{x}}$ with:
\begin{equation}
 \bar{\mathbf{x}} = \mathbf{A}^{-1} \mathbf{r}.
\label{eq:updateMulti}\end{equation}
Multigrid methods are based on the notion that if $\mathbf{A}$ derives from an elliptic operator, then the inverse $\mathbf{A}^{-1}$ in \eqref{eq:updateMulti} can be efficiently approximated by a combination of fixed point iterations and a coarse grid correction. The fixed point iterations efficiently reduce the oscillatory components of the error, and are therefore commonly referred to as smoothing \cite{Briggs2000}. The coarse grid correction again involves the (approximate) inverse of a matrix, analogous to \eqref{eq:updateMulti}, but corresponding to a coarser grid. This inverse can likewise be approximated by smoothing operations and a coarse grid correction on an even coarser grid, see Figure~\ref{fig:coarsening}. In multigrid methods, the smoothing operations and coarse grid correction are therefore applied recursively, until a grid is reached that is sufficiently coarse to enable direct inversion at a negligible computational expense. Our analyses are based on the notion that one entire multigrid cycle generates an approximation to the inverse of matrix $\mathbf{A}$. The result of the cycle is the vector $\tilde{\mathbf{x}}$, that approximates the vector $\bar{\mathbf{x}}$ in \eqref{eq:updateMulti}.

In this contribution we consider the multigrid V-cycle as outlined in Algorithm~\ref{alg:MGV}. The input parameters of this algorithm are the number of (remaining) recursive multigrid levels $\ell$, see Figure~\ref{fig:coarsening}, and the residual of the linear system $\mathbf{r}_\ell$ at the current level. It is emphasized that the number of multigrid levels follows a different convention than the number of hierarchical refinement levels. The coarsest multigrid level in the V-cycle is $\ell=1$, while the coarsest hierarchical refinement level is $m=0$, see Figure~\ref{fig:THB}. The output of the algorithm is the approximate solution $\tilde{\mathbf{x}}_\ell$ to $\mathbf{A}_\ell^{-1}\mathbf{r}_\ell$. The algorithm first initializes the vector $\tilde{\mathbf{x}}_\ell$ in line 2, and then in line 4 conducts a pre-smoothing step\footnote{Note that in this contribution we only apply a single pre-smoothing and post-smoothing operation in each cycle, and do not consider the possibility of multiple smoothing operations.}, i.e., fixed point iteration, with $\gamma>0$ denoting a relaxation parameter for stability and $\mathbf{M}_\ell^{-1}$ an approximate factorization or inverse of $\mathbf{A}_\ell$, e.g., Jacobi or Gauss-Seidel. Subsequently, smooth components of the error -- which were not effectively reduced by the aforementioned smoothing operation -- are treated by applying a coarse grid correction with coarser level $\ell-1$ in lines 7 to 9. The recursive nature of multigrid is implemented in Algorithm~\ref{alg:MGV}, by applying the V-cycle also to obtain an approximate solution to the coarse grid correction problem. The direct solver in line 14 is only applied at the coarsest level $\ell=1$. Note that the residual in the input of the algorithm should be interpreted as in \eqref{eq:residualMulti} only on the finest level, and is a restriction of a finer residual in the recursive applications of the algorithm in line 8. To enforce symmetry of the linear operator induced by the algorithm, the post-smoothing operation in line 12 is performed with the adjoint of $\mathbf{M}_\ell^{-1}$. This is relevant for Gauss-Seidel-type smoothers, and is realized by a reverse sweep. The approximation of the solution is returned in line 16. The coarse grid correction is treated in more detail in Section~\ref{sec:coarsening}, and the smoothing operations are discussed in detail in Section~\ref{sec:smoothers}.

\begin{center}\begin{minipage}[pt]{.95\textwidth}
 \begin{algorithm}[H]
  \caption{V-cycle($\ell,\mathbf{r}_\ell$)\label{alg:MGV}}
  \vspace{.2cm}
  \uIf{$\ell > 1$}
  {
    \vspace{.2cm}
    $\tilde{\mathbf{x}}_\ell = \mathbf{0}$ \hspace{.1cm} \textcolor{blue}{\# initialize approximation of $\mathbf{A}_\ell^{-1}\mathbf{r}_\ell$} \\
    \vspace{.2cm}
    \textcolor{blue}{\# pre-smooth} \\
    $\tilde{\mathbf{x}}_\ell = \tilde{\mathbf{x}}_\ell + \gamma \mathbf{M}_\ell^{-1} \mathbf{r}_\ell$ \hspace{.1cm} \textcolor{blue}{\# smooth} \\
    $\mathbf{r}_\ell = \mathbf{r}_\ell - \mathbf{A}_\ell\tilde{\mathbf{x}}_\ell$ \hspace{.1cm} \textcolor{blue}{\# update residual} \\
    \vspace{.2cm}
    \textcolor{blue}{\# coarse grid correction} \\
    $\mathbf{r}_{\ell-1} = \mathbf{R}_{\ell} \mathbf{r}_\ell$ \hspace{.1cm} \textcolor{blue}{\# restrict current residual to coarser grid} \\
    $\tilde{\mathbf{x}}_{\ell-1} = \mbox{V-cycle(}\ell-1,\mathbf{r}_{\ell-1}\mbox{)}$ \hspace{.1cm} \textcolor{blue}{\# compute coarse grid correction} \\
    $\tilde{\mathbf{x}}_\ell = \tilde{\mathbf{x}}_\ell +  \mathbf{R}_{\ell}^{\rm T} \tilde{\mathbf{x}}_{\ell-1}$ \hspace{.1cm} \textcolor{blue}{\# prolongate coarse grid correction} \\
    $\mathbf{r}_\ell = \mathbf{r}_\ell - \mathbf{A}_\ell\tilde{\mathbf{x}}_\ell$ \hspace{.1cm} \textcolor{blue}{\# update residual} \\
    \vspace{.2cm}
    \textcolor{blue}{\# post-smooth} \\
    $\tilde{\mathbf{x}}_\ell = \tilde{\mathbf{x}}_\ell + \gamma\mathbf{M}_\ell^{-{\rm T}} \mathbf{r}_\ell$ \hspace{.1cm} \textcolor{blue}{\# smooth} \\
    \vspace{.2cm}
  }
  \ElseIf{$\ell==1$}
  {
    \vspace{.2cm}
    $\tilde{\mathbf{x}}_\ell = \mathbf{A}_\ell^{-1} \mathbf{r}_\ell$ \hspace{.1cm} \textcolor{blue}{\# direct solve at coarsest level} \\ 
    \vspace{.2cm}
  }
  \vspace{.2cm}
  \Return{$\tilde{\mathbf{x}}_\ell$}
  \vspace{.2cm}
 \end{algorithm}
\end{minipage}\end{center}

The multigrid V-cycle in Algorithm~\ref{alg:MGV} can itself be applied as a solver, by iteratively performing the cycle to reduce the residual in every step, i.e., the approximation $\mathbf{x}'$ is simply updated by directly adding $\tilde{\mathbf{x}}$. This is generally true for multigrid cycles, e.g., also the W-cycle or FMG-cycle \cite{Briggs2000}. We consider the V-cycle, as this simple setup is already suitable to demonstrate the effectivity of the multigrid concept in immersed FEM. Additionally, when the pre-smoothing and post-smoothing operations are chosen such that these are adjoint, the V-cycle algorithm yields a symmetric positive definite (SPD) linear operator. This has the advantage that, instead of direct application of the V-cycle as a solver, it can also be employed as a preconditioner in a conjugate gradient (CG) algorithm. In the results presented in Section~\ref{sec:resultsMulti}, this multigrid-preconditioned CG-solver is applied.
The advantage of Krylov subspace solvers compared to multigrid as a standalone solver is that the convergence of Krylov methods is not purely governed by the smallest eigenmode in the system. Therefore, these are more robust to artifacts in the spectrum resulting from e.g., geometrical complexities such as the artificial coupling that is observed in Section~\ref{sec:3d} \cite{Trottenberg2000}.

\subsection{Restriction, prolongation and coarse grid correction}\label{sec:coarsening}

The restriction and prolongation operations to communicate between different grid sizes, such as those in Figure~\ref{fig:coarsening}, are essential aspects of the coarse grid correction in lines 7 to 10 of Algorithm~\ref{alg:MGV}. Under the usual assumption that the grid size is doubled in the mesh coarsening, the grid lines of the coarser level $\ell-1$ in Figure~\ref{fig:domain4} coincide with grid lines at the finer level $\ell$ in Figure~\ref{fig:domain8}. Therefore, the level $\ell-1$ space is nested in the level $\ell$ space and the basis functions on level $\ell-1$ can be represented identically by linear combinations of the basis functions on level $\ell$. Denoting by $\boldsymbol{\Phi}_{\ell}$ and $\boldsymbol{\Phi}_{\ell-1}$ vectors of basis functions on level $\ell$ and level $\ell-1$, respectively, there exists a matrix of coefficients $\mathbf{R}_\ell$ such that:
\begin{equation}
 \boldsymbol{\Phi}_{\ell-1} = \mathbf{R}_{\ell} \boldsymbol{\Phi}_\ell.
\label{eq:coarsening}\end{equation}
The matrix $\mathbf{R}_\ell$ defines the restriction operator. This restriction operator is employed in line 7 to restrict the residual at level $\ell$ to the level $\ell-1$ coarse mesh. In line 8 of Algorithm~\ref{alg:MGV}, the solution to the level $\ell-1$ problem is approximated by recursive application of the V-cycle algorithm, except at the coarsest level $\ell=1$ where a direct solution is carried out. The (approximate) solution to the level $\ell-1$ problem constitutes the coarse grid correction, which is prolongated and added to $\tilde{\mathbf{x}}_\ell$ at level $\ell$ in line 9. Let us note, that by the nesting of the approximation spaces at the different levels, the prolongation from level $\ell-1$ to level $\ell$ corresponds to injection. By virtue of relation \eqref{eq:coarsening} between the basis functions, we have the identities:
\begin{equation}
\tilde{\mathbf{x}}_{\ell-1}^T\boldsymbol{\Phi}_{\ell-1}= \tilde{\mathbf{x}}_{\ell-1}^T \left(\mathbf{R}_\ell\boldsymbol{\Phi}_\ell\right) = \left(\mathbf{R}_\ell^T\tilde{\mathbf{x}}_{\ell-1}\right)^T\boldsymbol{\Phi}_\ell
\end{equation}
In terms of the coefficients, the prolongation thus corresponds to the adjoint of the restriction operator. In line 10 the level $\ell$ residual is updated after the coarse grid correction.

\begin{remark}\label{rem:gridParamCorr}
 The coarse grid correction in Algorithm~\ref{alg:MGV} is performed purely algebraically. Therefore, the coarse problem inherits the weak formulation from the fine problem, without adapting the mesh-dependent parameters in the operators in \eqref{eq:operatorsMulti}. In essence, the finer level basis functions are replaced by the coarser level basis functions, such that the system matrix and residual at the coarser level can simply be obtained as $\mathbf{A}_{\ell-1} = \mathbf{R}_\ell\mathbf{A}_\ell\mathbf{R}_\ell^{\rm T}$ and $\mathbf{r}_{\ell-1} = \mathbf{R}_{\ell} \mathbf{r}_\ell$. An alternative approach is to adapt the operators in weak form \eqref{eq:weakMulti} to the coarser grid, as in e.g., \cite{Bramble1991,Gopalakrishnan2003}. The implementation of the algebraic approach in Algorithm~\ref{alg:MGV} is considerably simpler, however, and has been observed to adequately resolve the smooth eigenmodes in all considered cases. It should be noted that the extension of this simple algebraic coarsening approach to other grid-size dependencies is not verified in this contribution. Examples of other grid-size dependencies in the weak formulation are the stabilization parameter in Nitsche's method, see \cite{Nitsche1971,Embar2010} and for spectral effects specifically \cite{Harari2018}, and the ghost penalty, see \cite{Burman2010,BurmanHansbo2012}.
\end{remark}
% \begin{remark}\label{rem:gridParamCorr}
%  The coarse grid correction in Algorithm~\ref{alg:MGV} is performed purely algebraically. In essence, the finer level basis functions are replaced by the coarser level basis functions, without adapting the parameters in the operators in \eqref{eq:operatorsMulti}. The system matrix at the coarser level can therefore simply be obtained as $\mathbf{A}_{\ell-1} = \mathbf{R}_\ell\mathbf{A}_\ell\mathbf{R}_\ell^{\rm T}$. An alternative approach is to adapt the operators in weak form \eqref{eq:weakMulti} to the coarser grid. This is a relevant difference in the context of immersed finite elements, by virtue of the fact that the weak form is generally grid-size dependent, e.g., via mesh-dependent penalty parameters. The implementation of the algebraic approach in Algorithm~\ref{alg:MGV} is considerably simpler, and has been observed to adequately resolve the smooth eigenmodes in all considered cases. It should be noted, however, that the extension of this simple algebraic coarsening approach to other grid-size dependencies is not verified in this contribution. Examples of such grid-size dependencies are the stabilization parameter in Nitsche's method, see \cite{Nitsche1971,Embar2010} and for spectral effects specifically \cite{Harari2018}, and the ghost penalty, see \cite{Burman2010,BurmanHansbo2012}. 
% \end{remark}

\input{mesh.tex}

We have so far restricted ourselves to uniform meshes. For hierarchically refined meshes, as illustrated in Figure~\ref{fig:levelL}, the coarsening procedure in the multigrid method requires reconsideration. To enable an adequate approximation of functions that are smooth with respect to the local mesh width on level $\ell$, we apply a type of hierarchical derefinement to obtain a nested coarse space on level $\ell-1$. In practice, the mesh on level $\ell-1$ is the coarsest mesh for which the mesh of level $\ell$ can be obtained by a single level of hierarchical refinements, i.e., one uniform subdivision of a certain set of elements. The construction of such non-uniform coarser spaces is summarized in Algorithm~\ref{alg:coarsen}, and illustrated in Figure~\ref{fig:localCoarsening}. The algorithm denotes components of the mesh at level $\ell$ by: $K_i^m \in \mathcal{T}_h^m \subset \mathcal{T}_h$. As introduced in Figure~\ref{fig:THB}, $K_i^m$ denotes an active element at local refinement level $m$ with index $i$. The set of all active elements of local refinement level $m$ is denoted by $\mathcal{T}_h^m = \{K_i^m\}$, and $\mathcal{T}_h=\cup_{m=0}^{m=M} \mathcal{T}_h^m$ denotes the full mesh, with $M$ denoting the number of local refinement levels. Elements of the coarser mesh at level $\ell-1$ are denoted as $k_i^m \in \mathcal{T}_{2h}^m \subset \mathcal{T}_{2h}$. The subscript $2h$ for the coarse mesh is to be conceived of as a symbolic notation. Note that the elements in $\mathcal{T}_{2h}^{m+1}$ are of the same size as the elements in $\mathcal{T}_h^m$. In line 1 of the algorithm, the coarser mesh is initialized as a uniform mesh with only unrefined elements of local refinement level $m=0$, i.e., $\mathcal{T}_{2h} = \mathcal{T}_{2h}^0$. Line 3 initiates a loop over the refinement levels $0 \leq m \leq M-1$. Within this loop, the algorithm loops over all the elements $k_i^m \in \mathcal{T}_{2h}^m$, that are currently $m$ times refined. If $k_i^m$ intersects an element of $\mathcal{T}_h$ that is refined more than $m$ times, element $k_i^m$ is refined in lines 7 to 9. Note that line 7 abuses notation to simplify the expression, and that the refining of element $k_i^m$ implies removing $k_i^m$ from $\mathcal{T}_{2h}^m$ and adding the refinements to $\mathcal{T}_{2h}^{m+1}$. Note that this hierarchical derefinement procedure differs from the approach in \cite{Hofreither2016}, which employs an existing hierarchy of refined grids for the coarse grid corrections.

\begin{center}\begin{minipage}{.95\textwidth}
 \begin{algorithm}[H]
  \caption{Coarsen$\left(\mathcal{T}_h\right)$\label{alg:coarsen}}
  \vspace{.2cm}
  initialize $\mathcal{T}_{2h}$ \hspace{.1cm} \textcolor{blue}{\# initialize coarse mesh $\mathcal{T}_{2h}=\mathcal{T}_{2h}^0$ with uniform unrefined elements} \\
  \vspace{.2cm}
  \textcolor{blue}{\# loop over local refinement levels} \\
  \For{$m \in \{0,...,M-1\}$}
  {
    \vspace{.2cm}
    \textcolor{blue}{\# loop over coarse mesh elements of level $m$} \\
    \For{$k_i^m \in \mathcal{T}_{2h}^m$} 
    {
      \vspace{.2cm}
      \textcolor{blue}{\# check if $k_i^m$ intersects elements of fine mesh $\mathcal{T}_h$ of local refinement level $>m$} \\
      \If{$\left( \cup_{\tilde{m} = m+1}^{\tilde{m}=M} \mathcal{T}_h^{\tilde{m}} \right) \cap k_i^m \neq \emptyset$}
      {
        \vspace{.2cm}
        refine $k_i^m$ \\
        \vspace{.2cm}
      }
      \vspace{.2cm}
    }
    \vspace{.2cm}
  }
  \vspace{.2cm}
  \Return{$\mathcal{T}_{2h}$} \\
  \vspace{.2cm}
 \end{algorithm}
\end{minipage}\end{center}

\subsection{Smoothers for immersed finite element methods}\label{sec:smoothers}

In lines 4 and 12 of Algorithm~\ref{alg:MGV}, smoothing operations, i.e., fixed point iterations, are performed to resolve the components of the error that cannot be adequately captured by the coarse grid. Effective application of the multigrid algorithm requires that the eigenvalues of $\gamma\mathbf{M}^{-1}\mathbf{A}$ that correspond to non-smooth eigenfunctions are close to $1$. Note that in the formulations from here on, the subscripts indicating the level in the multigrid solver are omitted to simplify the notation. All these formulations are independent of the level $\ell$, however. Stability of fixed point iterations requires:
\begin{equation}
 0 \leq \lambda_{\rm min}\left( \gamma\mathbf{M}^{-1}\mathbf{A} \right) \leq \lambda_{\rm max}\left( \gamma\mathbf{M}^{-1}\mathbf{A} \right) \leq 2,
\label{eq:stabilityMulti}\end{equation}
with $\lambda_{\rm min}(\cdot)$ and $\lambda_{\rm max}(\cdot)$ denoting the smallest and largest eigenvalues, such that the spectral radius of the fixed point iteration is bounded:
\begin{equation}
 \rho\left( \mathbf{I} - \gamma\mathbf{M}^{-1}\mathbf{A} \right) < 1,
\end{equation}
with $\mathbf{I}$ denoting the identity matrix with the same size as system matrix $\mathbf{A}$. In the case that $\lambda_{\rm max}\left(\gamma\mathbf{M}^{-1}\mathbf{A}\right)>2$, the error component in the direction of the eigenvector corresponding to the largest eigenvalue will increase with smoothing, which may result in divergence. For this reason, smoothers such as Jacobi and additive Schwarz require a sufficiently small relaxation parameter $\gamma$. Smoothers such as Gauss-Seidel and multiplicative Schwarz are unconditionally stable and do not require relaxation, see e.g., \cite[Theorem 4.10 \& 14.9]{Saad}. As already mentioned in the description of Algorithm~\ref{alg:MGV} in Section~\ref{sec:algorithm}, it should be noted that symmetry of the linear operator that is induced by the V-cycle insists that the post-smoothing operation is the adjoint of the pre-smoothing operation. Furthermore, it should be mentioned that the computational efficiency can potentially be improved by performing multiple smoothing operations in each cycle. Such enhancements are, however, not considered in this contribution.

Jacobi iterations are not suitable as a smoother for immersed finite elements, which is illustrated by the example in Figure~\ref{fig:Jacobi}. The smallest eigenmode in Figure~\ref{fig:JacobiSmall} with a very small eigenvalue is barely affected by the smoothing, and cannot be captured on a coarser grid. Furthermore, the relatively large eigenmodes caused by almost linear dependencies as in Figure~\ref{fig:JacobiLarge} impose a small relaxation parameter, which further impairs the conditioning. The following subsections examine the suitability of a Gauss-Seidel smoother and Schwarz-type smoothers based on the preconditioner for immersed finite elements developed in \cite{Prenter2019}.

\subsubsection*{Gauss-Seidel}

\begin{figure}[pt]
  \begin{subfigure}{.49\textwidth}
   \vspace{.2cm}
   \includegraphics[height=5.1cm]{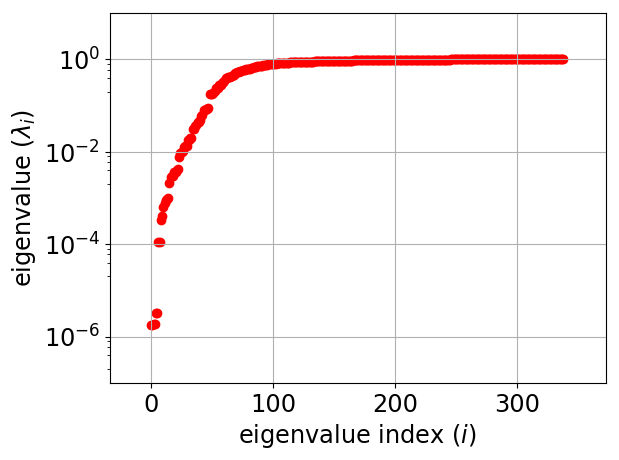}
   \caption{Eigenvalue spectrum\label{fig:GaussSeidelSpectrum}}
  \end{subfigure}
  \hspace{.5cm}
  \begin{subfigure}{.49\textwidth}
   \includegraphics[height=5cm]{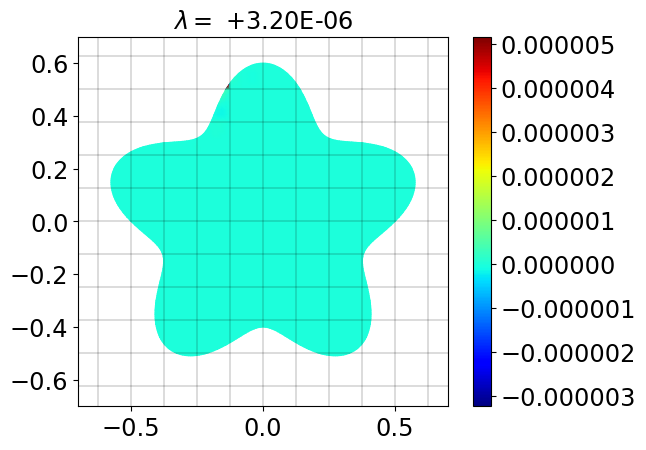}
   \vspace{.3cm}
   \caption{Very small mode\label{fig:GaussSeidelSmall}}
  \end{subfigure}
  \\
  \rule{0pt}{0pt}\hspace{.8cm}
  \begin{subfigure}{.49\textwidth}
   \vspace{.2cm}
   \includegraphics[height=5cm]{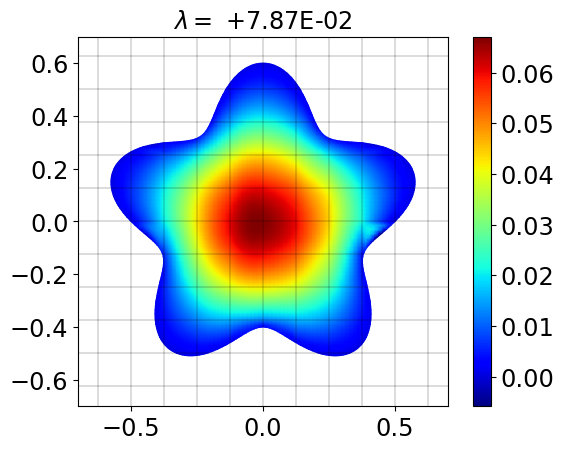}
   \caption{Smooth mode\label{fig:GaussSeidelSmooth}}
  \end{subfigure}
  \hspace{-.4cm}
  \begin{subfigure}{.49\textwidth}
   \vspace{.2cm}
   \includegraphics[height=5cm]{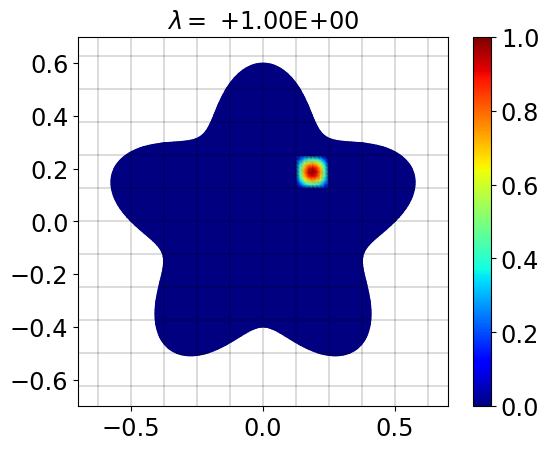}
   \caption{Element of largest eigenspace\label{fig:GaussSeidelLarge}}
  \end{subfigure}
 \caption{Typical spectrum and characteristic eigenmodes of an immersed system that is preconditioned by a double fixed point iteration with Gauss-Seidel.\label{fig:GaussSeidel}}
\end{figure}

A typical spectrum and characteristic eigenmodes with standard Gauss-Seidel preconditioning are shown in Figure~\ref{fig:GaussSeidel}. Similar to Figure~\ref{fig:Jacobi} for Jacobi preconditioning, these figures correspond to the Laplace operator on the geometry in Figure~\ref{fig:domainStar} with quadratic Lagrange basis functions and boundary conditions imposed by the penalty method. To obtain a symmetric preconditioner, a double fixed point iteration with adjoint Gauss-Seidel operations is applied in these figures:
\begin{equation}
 \left( \mathbf{I} - \mathbf{M}^{-{\rm T}} \mathbf{A} \right) \left( \mathbf{I} - \mathbf{M}^{-1} \mathbf{A} \right) \mathbf{y}_{\lambda_i} = \left( 1 - \lambda_i \right) \mathbf{y}_{\lambda_i},
\end{equation}
with $\left( \mathbf{I} - \mathbf{M}^{-1} \mathbf{A} \right)$ corresponding to the initial Gauss-Seidel sweep, $\left( \mathbf{I} - \mathbf{M}^{-{\rm T}} \mathbf{A} \right)$ corresponding to the reverse sweep, and $\mathbf{y}_{\lambda_i}$ denoting the eigenvector corresponding to the $i$th eigenvalue $\lambda_i$. Hence, Figure~\ref{fig:GaussSeidel} presents the eigenmodes of the system $\left( \mathbf{M}^{-1} + \mathbf{M}^{-{\rm T}} - \mathbf{M}^{-{\rm T}}\mathbf{A}\mathbf{M}^{-1} \right)\mathbf{A}$. As will follow in \eqref{eq:eigenMulti}, this is actually very similar to the V-cycle, except for the omission of the coarse grid correction. To reduce the computational cost, the Gauss-Seidel routine is implemented with a graph coloring algorithm\footnote{This algorithm divides the basis functions in sets, or colors, such that the supports of basis functions with the same color do not intersect. The union over the sets of all colors constitutes the full approximation space. Note that, for Lagrange basis functions on uniform grids, this requires $(p+1)^d$ different colors. As basis functions of the same color do not intersect, updating the approximation of the solution, $\tilde{\mathbf{x}}$, at an index corresponding to a certain color, does not affect the residual, $\mathbf{r}$, at the other indices corresponding to that color. This enables a Gauss-Seidel routine that sweeps over all indices of a certain color at once, instead of sequentially updating each index of the approximation of the solution and updating residual accordingly.} \cite{Adams2003}. In Figure~\ref{fig:GaussSeidelSmall} a very small eigenvalue is plotted, which is similar to the smallest eigenmode with Jacobi preconditioning in Figure~\ref{fig:JacobiSmall}. Note that the eigenvalues of these eigenfunctions differ by a factor of (approximately) $2$, which is an expected consequence of the double iteration with Gauss-Seidel versus the single iteration with Jacobi. In Figure~\ref{fig:GaussSeidelSmooth} it is shown that also the spectrum with Gauss-Seidel preconditioning contains an eigenmode that is in close correspondence with the smallest analytical eigenmode of the PDE, and is similar to the usual smallest eigenmode with mesh-fitting techniques. Since the colors are selected such that basis functions of the same color do not intersect, the unit vectors corresponding to basis functions with the last color of the graph coloring algorithm yield an eigenspace with an eigenvalue of exactly $1$. This can be observed in the spectrum in Figure~\ref{fig:GaussSeidelSpectrum} and is illustrated in Figure~\ref{fig:GaussSeidelLarge}.

\begin{figure}[pt]
  \begin{subfigure}{.49\textwidth}
   \vspace{.2cm}
   \includegraphics[height=5.1cm]{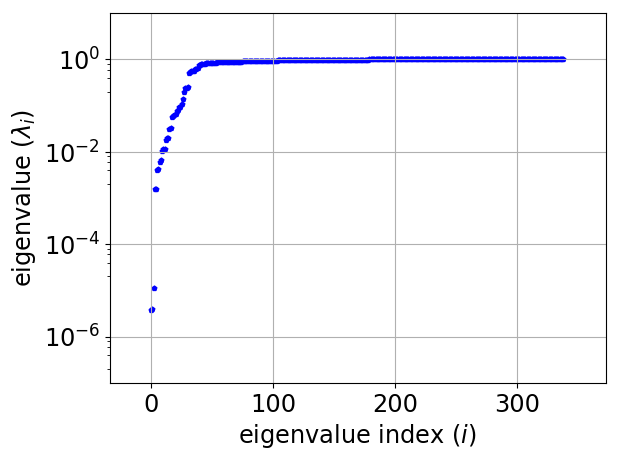}
   \caption{Eigenvalue spectrum\label{fig:GaussSeidelSpectrumMultigrid}}
  \end{subfigure}
  \hspace{.5cm}
  \begin{subfigure}{.49\textwidth}
   \includegraphics[height=5cm]{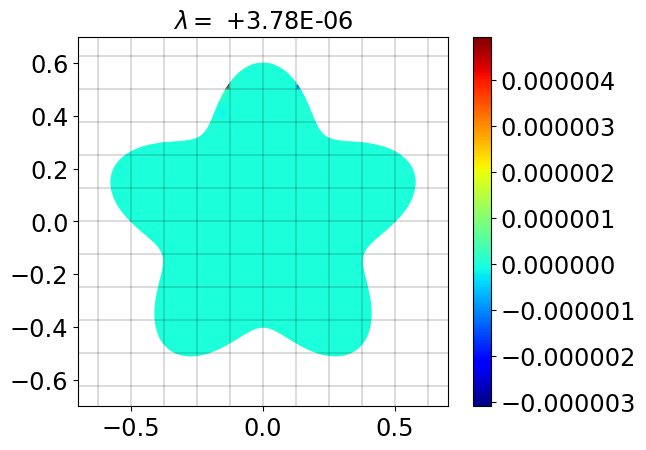}
   \vspace{.3cm}
   \caption{Smallest mode\label{fig:GaussSeidelSmallMultigrid}}
  \end{subfigure}
 \caption{Spectrum and a characteristic very small eigenmode of an immersed system that is preconditioned by a two-level V-cycle with Gauss-Seidel as smoother.\label{fig:GaussSeidelMultigrid}}
\end{figure}

Gauss-Seidel is not suitable as a smoother for immersed finite element methods, despite the unconditional stability by which -- in contrast to Jacobi -- it does not require more relaxation in an immersed setting than in a mesh-fitting setting. The problem that renders Gauss-Seidel ineffective for immersed finite element methods is the small eigenmode plotted in Figure~\ref{fig:GaussSeidelSmall}. Similar to Jacobi, the eigenvalue of this mode is too small to be adequately treated by the smoothing operations, and it is also not resolved by the coarse grid correction. This is shown in Figure~\ref{fig:GaussSeidelMultigrid}, which displays the spectrum of the same problem preconditioned by the V-cycle in Algorithm~\ref{alg:MGV} with $\ell=2$ levels and a single Gauss-Seidel sweep as smoother. This eigenvalue problem can be formulated as:
\begin{equation}
 \left( \mathbf{I} - \mathbf{M}_\ell^{-{\rm T}} \mathbf{A}_\ell \right) \left( \mathbf{I} - \mathbf{R}_\ell^{\rm T} \mathbf{A}_{\ell-1}^{-1} \mathbf{R}_\ell \mathbf{A}_\ell \right)\left( \mathbf{I} - \mathbf{M}_\ell^{-1} \mathbf{A}_\ell \right) \mathbf{y}_{\lambda_i} = \left( 1 - \lambda_i \right) \mathbf{y}_{\lambda_i},
\label{eq:eigenMulti}\end{equation}
or as:
\begin{equation}
\mbox{V-cycle}\big(\ell=2,\mathbf{r}_\ell=\mathbf{A}\mathbf{y}_{\lambda_i}\big) = \lambda_i \mathbf{y}_{\lambda_i}.
\end{equation}
While a comparison of the spectra in Figures~\ref{fig:GaussSeidelSpectrum} and \ref{fig:GaussSeidelSpectrumMultigrid}, i.e., without and with the coarse grid correction, reveals that several small eigenmodes are resolved by the coarse grid correction, Figures~\ref{fig:GaussSeidelSmall} and \ref{fig:GaussSeidelSmallMultigrid} demonstrate that both spectra contain approximately the same small eigenmode. The eigenvalues of this mode are barely affected by the coarse grid correction, with eigenvalue $3.20\cdot10^{-6}$ for the double fixed point iteration and eigenvalue $3.78\cdot10^{-6}$ for the full V-cycle.

\clearpage

\subsubsection*{Additive Schwarz}

The spectra and eigenfunctions with Jacobi and Gauss-Seidel in Figures~\ref{fig:Jacobi} and \ref{fig:GaussSeidel} clearly demonstrate that these are not robust to cut elements. This is consistent with the analysis of the conditioning problems in \cite{SIPIC}, which points out that diagonal preconditioners do not adequately mitigate the almost linear dependencies that occur in immersed finite element methods. In \cite{Prenter2019} it is derived that almost linearly dependent basis functions can be effectively treated collectively when these are inverted in a block manner by a Schwarz-type method. Based on the additive Schwarz lemma, \cite{Prenter2019} shows that additive Schwarz preconditioning is actually a very natural approach to resolve the small eigenmodes caused by almost linear dependencies. Furthermore, it is demonstrated that in terms of the condition number and the number of iterations in an iterative solution method, immersed methods with additive Schwarz preconditioning behave similar to mesh-fitting techniques. This section therefore examines the spectrum of immersed methods with additive Schwarz preconditioning, to establish the suitability as a smoother in a multigrid method.

In additive Schwarz preconditioning, a set of $N$ index blocks is selected, which correspond to sets of basis functions. For each index block $j\leq N$, the system matrix $\mathbf{A} \in \mathbb{R}^{n \times n}$ is restricted to the indices in the block, denoted by $\mathbf{A}_j \in \mathbb{R}^{n_j \times n_j}$. The block matrices are then inverted and prolongated to a matrix of size $n \times n$. The additive Schwarz preconditioner is obtained by summing these matrices:
\begin{equation}
 \mathbf{M}^{-1} = \sum_{j=1}^N \mathbf{P}_j \underbrace{\left(\mathbf{P}_j^{\rm T}\mathbf{A}\mathbf{P}_j \right)^{-1}}_{\mathbf{A}_j^{-1}}\mathbf{P}_j^{\rm T},
\end{equation}
with $\mathbf{P}_j\in \mathbb{R}^{n \times n_j}$ a matrix that prolongates a block vector $\mathbf{y}_j \in \mathbb{R}^{n_j}$ to a vector $\mathbf{P}_j \mathbf{y}_j = \mathbf{y} \in \mathbb{R}^{n}$ -- with nonzero entries only at the indices in block $j$ -- and the transpose of $\mathbf{P}_j$ a restriction operator that restricts a vector $\mathbf{z}\in\mathbb{R}^{n}$ to a block vector $\mathbf{P}_j^{\rm T} \mathbf{z} = \mathbf{z}_j \in \mathbb{R}^{n_j}$ -- containing only the indices in block $j$.

An essential aspect of additive Schwarz preconditioners is the choice of the index blocks. It is pointed out in \cite{Prenter2019} that almost linearly dependent basis functions are required to be in an index block together. Furthermore, it is demonstrated that devising a block for each cut element with all basis functions supported on it is an effective strategy to satisfy this requirement for uniform grids. As demonstrated in \cite{John}, however, it is not trivial to generalize this concept to locally refined meshes. Therefore this contribution applies an alternative strategy to select the Schwarz blocks based on so-called encapsulating supports, which is inspired by the Schwarz-type smoother developed for divergence-conforming discretizations in \cite{Coley2018}. Accordingly, for every basis function a block is devised, containing all the basis functions whose support completely lies inside the support of the basis function associated to the block, see Figure~\ref{fig:inclusionMulti}. Note that in this contribution the support of a basis function refers to the support within the physical domain. The block that is associated to function $\phi_j$ is defined as: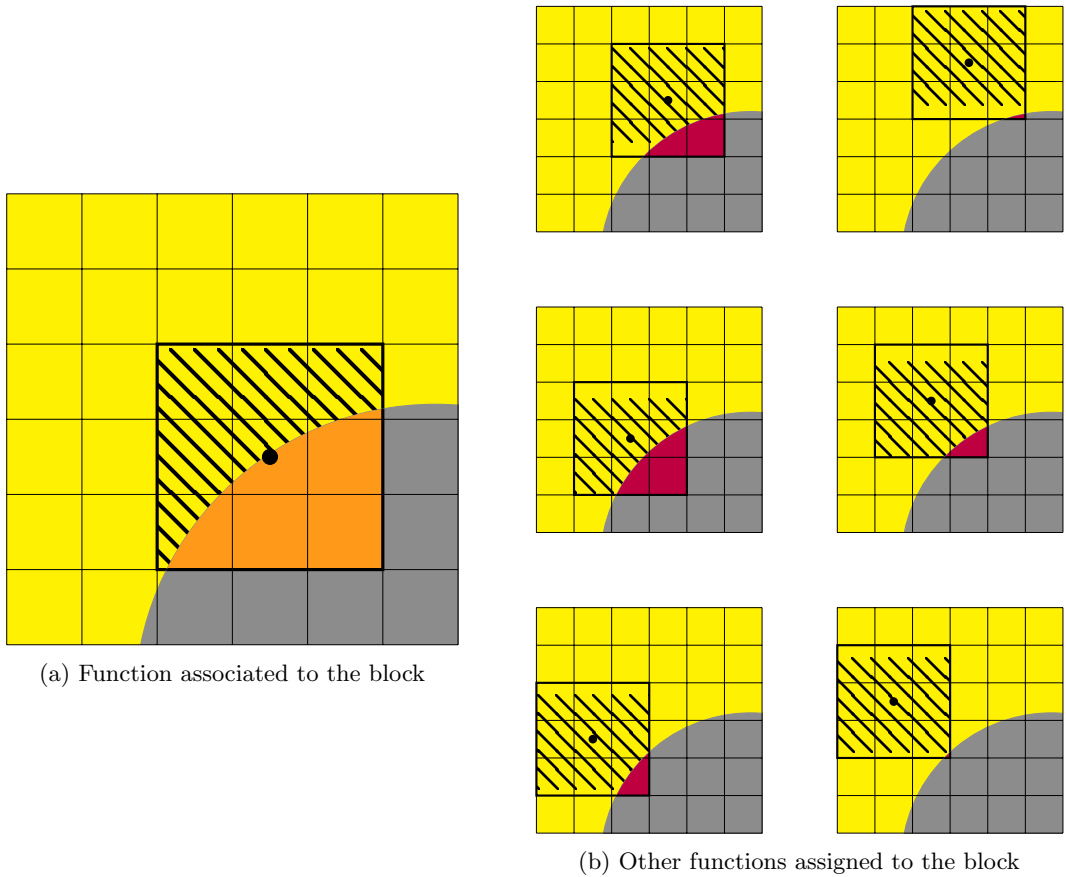
\begin{figure}[pt]
    \tikzset
    {
        hatch distance/.store in=\hatchdistance,
        hatch distance=8pt,
        hatch thickness/.store in=\hatchthickness,
        hatch thickness=1pt
    }

    \makeatletter
    \pgfdeclarepatternformonly[\hatchdistance,\hatchthickness]{flexible hatch}
    {\pgfqpoint{0pt}{0pt}}
    {\pgfqpoint{\hatchdistance}{-\hatchdistance}}
    {\pgfpoint{\hatchdistance-1pt}{\hatchdistance-1pt}}%
    {
        \pgfsetcolor{\tikz@pattern@color}
        \pgfsetlinewidth{\hatchthickness}
        \pgfpathmoveto{\pgfqpoint{0pt}{0pt}}
        \pgfpathlineto{\pgfqpoint{\hatchdistance}{-\hatchdistance}}
        \pgfusepath{stroke}
    }
    \makeatother
 \begin{center}
  \begin{subfigure}{.45\textwidth}
   \begin{center}
    \begin{tikzpicture}
     [ scale=1 ]
     \definecolor{myOrange}{rgb}{1.1,.6,.1}
          
     \fill[yellow, fill opacity=1] (0,0) rectangle (6,6);
     \draw[black, pattern=flexible hatch, hatch distance=10pt, hatch thickness=1.5pt, very thick] (2,1) rectangle (5,4);
     \begin{scope} 
      \clip (0,0) rectangle (6,6);
      \fill[black!45,fill opacity=1,draw=black!45] (5.7,-.8) circle (4);
     \end{scope}
     \begin{scope} 
      \clip (0,0) rectangle (6,6);
      \clip (5.7,-.8) circle (4);
      \draw[black, very thick, fill=myOrange, fill opacity=1] (2,1) rectangle (5,4);
     \end{scope}
     \fill[black,fill opacity=1,draw=black] (3.5,2.5) circle (.1);
     \draw[step=1,black,ultra thin] (0,0) grid (6,6);
    \end{tikzpicture}
    \caption{Function associated to the block\label{fig:associatedMulti}}
   \end{center}
  \end{subfigure}
  \hfill
  \begin{subfigure}{.54\textwidth}
   \begin{center}
    \begin{tikzpicture}
     [ scale=.5 ]
     \definecolor{myPurple}{rgb}{.5,0.0,5}     
     \def\dx{-4}
     \def\dy{8}
     \def\funcx{0}
     \def\funcy{1}
     \fill[yellow, fill opacity=1] (0+\dx,0+\dy) rectangle (6+\dx,6+\dy);
     \draw[black, pattern=flexible hatch, thick] (2+\dx+\funcx,1+\dy+\funcy) rectangle (5+\dx+\funcx,4+\dy+\funcy);
     \begin{scope} 
      \clip (0+\dx,0+\dy) rectangle (6+\dx,6+\dy);
      \fill[black!45,fill opacity=1,draw=black!45] (5.7+\dx,-.8+\dy) circle (4);
     \end{scope}
     \begin{scope} 
      \clip (0+\dx,0+\dy) rectangle (6+\dx,6+\dy);
      \clip (5.7+\dx,-.8+\dy) circle (4);
      \draw[black, thick, fill=purple, fill opacity=1] (2+\dx+\funcx,1+\dy+\funcy) rectangle (5+\dx+\funcx,4+\dy+\funcy);
     \end{scope}
     \fill[black,fill opacity=1,draw=black] (3.5+\dx+\funcx,2.5+\dy+\funcy) circle (.1);
     \draw[step=1,black,ultra thin] (0+\dx,0+\dy) grid (6+\dx,6+\dy);
     \def\dx{4}
     \def\dy{8}
     \def\funcx{0}
     \def\funcy{2}
     \fill[yellow, fill opacity=1] (0+\dx,0+\dy) rectangle (6+\dx,6+\dy);
     \draw[black, pattern=flexible hatch, thick] (2+\dx+\funcx,1+\dy+\funcy) rectangle (5+\dx+\funcx,4+\dy+\funcy);
     \begin{scope} 
      \clip (0+\dx,0+\dy) rectangle (6+\dx,6+\dy);
      \fill[black!45,fill opacity=1,draw=black!45] (5.7+\dx,-.8+\dy) circle (4);
     \end{scope}
     \begin{scope} 
      \clip (0+\dx,0+\dy) rectangle (6+\dx,6+\dy);
      \clip (5.7+\dx,-.8+\dy) circle (4);
      \draw[black, thick, fill=purple, fill opacity=1] (2+\dx+\funcx,1+\dy+\funcy) rectangle (5+\dx+\funcx,4+\dy+\funcy);
     \end{scope}
     \fill[black,fill opacity=1,draw=black] (3.5+\dx+\funcx,2.5+\dy+\funcy) circle (.1);
     \draw[step=1,black,ultra thin] (0+\dx,0+\dy) grid (6+\dx,6+\dy);
     \def\dx{-4}
     \def\dy{0}
     \def\funcx{-1}
     \def\funcy{0}
     \fill[yellow, fill opacity=1] (0+\dx,0+\dy) rectangle (6+\dx,6+\dy);
     \draw[black, pattern=flexible hatch, thick] (2+\dx+\funcx,1+\dy+\funcy) rectangle (5+\dx+\funcx,4+\dy+\funcy);
     \begin{scope} 
      \clip (0+\dx,0+\dy) rectangle (6+\dx,6+\dy);
      \fill[black!45,fill opacity=1,draw=black!45] (5.7+\dx,-.8+\dy) circle (4);
     \end{scope}
     \begin{scope} 
      \clip (0+\dx,0+\dy) rectangle (6+\dx,6+\dy);
      \clip (5.7+\dx,-.8+\dy) circle (4);
      \draw[black, thick, fill=purple, fill opacity=1] (2+\dx+\funcx,1+\dy+\funcy) rectangle (5+\dx+\funcx,4+\dy+\funcy);
     \end{scope}
     \fill[black,fill opacity=1,draw=black] (3.5+\dx+\funcx,2.5+\dy+\funcy) circle (.1);
     \draw[step=1,black,ultra thin] (0+\dx,0+\dy) grid (6+\dx,6+\dy);
     \def\dx{4}
     \def\dy{0}
     \def\funcx{-1}
     \def\funcy{1}
     \fill[yellow, fill opacity=1] (0+\dx,0+\dy) rectangle (6+\dx,6+\dy);
     \draw[black, pattern=flexible hatch, thick] (2+\dx+\funcx,1+\dy+\funcy) rectangle (5+\dx+\funcx,4+\dy+\funcy);
     \begin{scope} 
      \clip (0+\dx,0+\dy) rectangle (6+\dx,6+\dy);
      \fill[black!45,fill opacity=1,draw=black!45] (5.7+\dx,-.8+\dy) circle (4);
     \end{scope}
     \begin{scope} 
      \clip (0+\dx,0+\dy) rectangle (6+\dx,6+\dy);
      \clip (5.7+\dx,-.8+\dy) circle (4);
      \draw[black, thick, fill=purple, fill opacity=1] (2+\dx+\funcx,1+\dy+\funcy) rectangle (5+\dx+\funcx,4+\dy+\funcy);
     \end{scope}
     \fill[black,fill opacity=1,draw=black] (3.5+\dx+\funcx,2.5+\dy+\funcy) circle (.1);
     \draw[step=1,black,ultra thin] (0+\dx,0+\dy) grid (6+\dx,6+\dy);
     \def\dx{-4}
     \def\dy{-8}
     \def\funcx{-2}
     \def\funcy{0}
     \fill[yellow, fill opacity=1] (0+\dx,0+\dy) rectangle (6+\dx,6+\dy);
     \draw[black, pattern=flexible hatch, thick] (2+\dx+\funcx,1+\dy+\funcy) rectangle (5+\dx+\funcx,4+\dy+\funcy);
     \begin{scope} 
      \clip (0+\dx,0+\dy) rectangle (6+\dx,6+\dy);
      \fill[black!45,fill opacity=1,draw=black!45] (5.7+\dx,-.8+\dy) circle (4);
     \end{scope}
     \begin{scope} 
      \clip (0+\dx,0+\dy) rectangle (6+\dx,6+\dy);
      \clip (5.7+\dx,-.8+\dy) circle (4);
      \draw[black, thick, fill=purple, fill opacity=1] (2+\dx+\funcx,1+\dy+\funcy) rectangle (5+\dx+\funcx,4+\dy+\funcy);
     \end{scope}
     \fill[black,fill opacity=1,draw=black] (3.5+\dx+\funcx,2.5+\dy+\funcy) circle (.1);
     \draw[step=1,black,ultra thin] (0+\dx,0+\dy) grid (6+\dx,6+\dy);
     \def\dx{4}
     \def\dy{-8}
     \def\funcx{-2}
     \def\funcy{1}
     \fill[yellow, fill opacity=1] (0+\dx,0+\dy) rectangle (6+\dx,6+\dy);
     \draw[black, pattern=flexible hatch, thick] (2+\dx+\funcx,1+\dy+\funcy) rectangle (5+\dx+\funcx,4+\dy+\funcy);
     \begin{scope} 
      \clip (0+\dx,0+\dy) rectangle (6+\dx,6+\dy);
      \fill[black!45,fill opacity=1,draw=black!45] (5.7+\dx,-.8+\dy) circle (4);
     \end{scope}
     \begin{scope} 
      \clip (0+\dx,0+\dy) rectangle (6+\dx,6+\dy);
      \clip (5.7+\dx,-.8+\dy) circle (4);
      \draw[black, thick, fill=purple, fill opacity=1] (2+\dx+\funcx,1+\dy+\funcy) rectangle (5+\dx+\funcx,4+\dy+\funcy);
     \end{scope}
     \fill[black,fill opacity=1,draw=black] (3.5+\dx+\funcx,2.5+\dy+\funcy) circle (.1);
     \draw[step=1,black,ultra thin] (0+\dx,0+\dy) grid (6+\dx,6+\dy);
    \end{tikzpicture}
    \caption{Other functions assigned to the block\label{fig:assignedMulti}}
   \end{center}
  \end{subfigure}
  
  \caption{Greville abscissae and supports of quadratic B-splines assigned to a block. Figure~\ref{fig:associatedMulti} displays the support within the physical domain of the function associated to the block in orange. Figure~\ref{fig:assignedMulti} displays the supports of the other functions assigned to the block in purple. The fictitious support of the functions is not considered, but is hatched to increase the clarity of the figures. Note that the support of the functions in Figure~\ref{fig:assignedMulti} completely lies inside the support of the function in Figure~\ref{fig:associatedMulti}.\label{fig:inclusionMulti}}
 \end{center}
\end{figure}
\begin{equation}
 \left\{ \, \phi_k \, : \,  {\rm supp}_\Omega\left(\phi_k\right) \subseteq {\rm supp}_\Omega\left(\phi_j\right) \, \right\},
\end{equation}
with ${\rm supp}_\Omega(\phi_k)$ denoting the support of basis function $\phi_k$ within physical domain $\Omega$. For vector-valued problems, separate blocks are devised for basis functions describing different vectorial components of the solution, similar to the blocks in \cite{Prenter2019}. By construction, each block therefore contains the basis function associated to it, and for untrimmed (truncated hierarchical) B-splines this approach yields an approximately diagonal preconditioner. Let us note here the importance of the truncation of the basis functions in the locally refined approximations. As non-truncated hierarchical bases yield a very large number of basis functions with overlapping supports, this correspondingly results in very large blocks in the additive Schwarz preconditioner, which leads to significant computational costs. Trimmed basis functions on small cut elements -- which can be almost linearly dependent -- are also assigned to blocks associated to other functions. This satisfies the requirement formulated in \cite{Prenter2019} that almost linearly dependent basis functions need to be in a block together, and therefore resolves the small eigenmodes that are characteristic for immersed finite element methods. This strategy to select the Schwarz blocks is directly applicable to both B-splines on uniform grids and truncated hierarchical B-splines on non-uniform grids, in contrast to the element-wise strategy in \cite{Prenter2019}. It is, however, not natural to directly apply this strategy to Lagrange basis functions, as these are not all supported on the same number of elements. For the Lagrange bases, blocks are therefore only devised for the nodal basis functions\footnote{In this contribution different Lagrange basis functions are indicated as nodal, edge, face and volume functions. With $d$ denoting the number of dimensions: \emph{nodal} functions attain the value $1$ at a vertex of the grid and span $2^d$ elements, \emph{edge} functions attain the value $1$ on an edge and span $2^{d-1}$ elements, \emph{face} functions attain the value $1$ on a face and span $2^{d-2}$ elements, and \emph{volume} or \emph{element internal} functions attain the value $1$ inside an element and span $1$ element. Note that face functions are not considered in the two-dimensional example.}. This implies that with Lagrange bases a block is devised for every cluster of $2^d$ elements. Note that this yields approximately the same number of blocks for uniform Lagrange bases as for uniform B-spline bases, but, in contrast, does not reduce to a purely diagonal treatment of untrimmed basis functions for Lagrange bases. It should be mentioned that the block selection described here is not the only possible and effective method to select blocks for immersed finite elements, and interested readers are directed to \cite{BeiraoDaVeiga2012,BeiraoDaVeiga2013Schwarz,delaRiva2018} for a study of suitable block selections in isogeometric analysis and to the reference works \cite{Smith1996,Toselli2005} for considerations regarding the block selections with traditional finite element bases.

\begin{figure}[pt]
  \begin{subfigure}{.49\textwidth}
   \vspace{.2cm}
   \includegraphics[height=5.1cm]{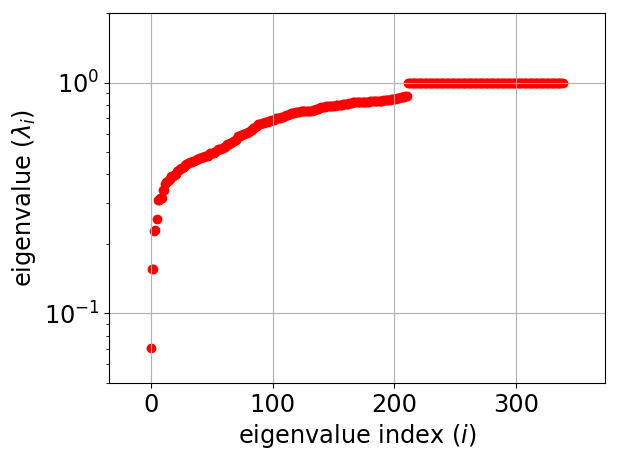}
   \caption{Eigenvalue spectrum\label{fig:AdditiveSchwarzSpectrum}}
  \end{subfigure}
  \hspace{.35cm}
  \begin{subfigure}{.49\textwidth}
   \includegraphics[height=5cm]{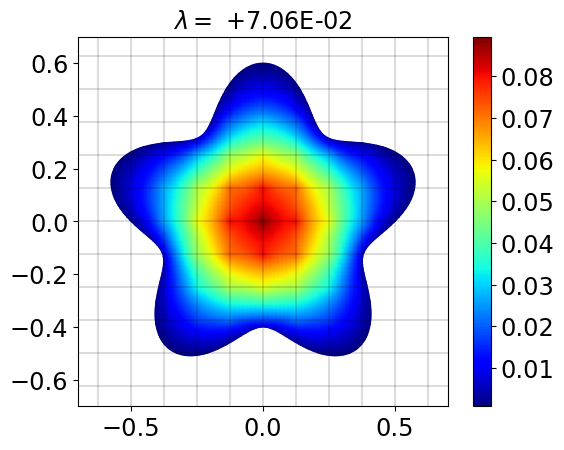}
   \vspace{.3cm}
   \caption{Smooth mode (smallest)\label{fig:AdditiveSchwarzSmooth}}
  \end{subfigure}
  \\
  \rule{0pt}{0pt}\hspace{.8cm}
  \begin{subfigure}{.49\textwidth}
   \vspace{.2cm}
   \includegraphics[height=5cm]{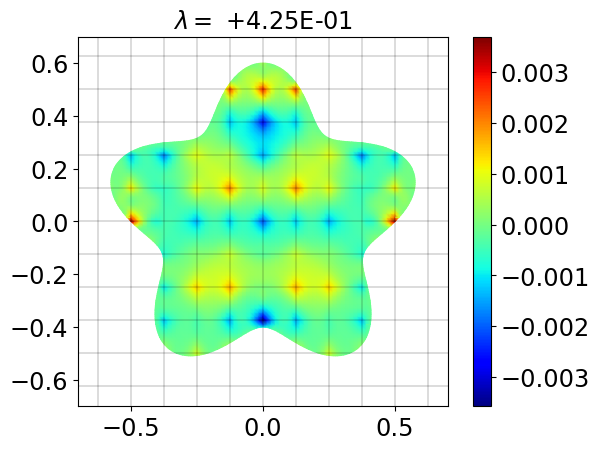}
   \caption{Nodal mode\label{fig:AdditiveSchwarzNodal}}
  \end{subfigure}
  \hspace{-.5cm}
  \begin{subfigure}{.49\textwidth}
   \vspace{.2cm}
   \includegraphics[height=5cm]{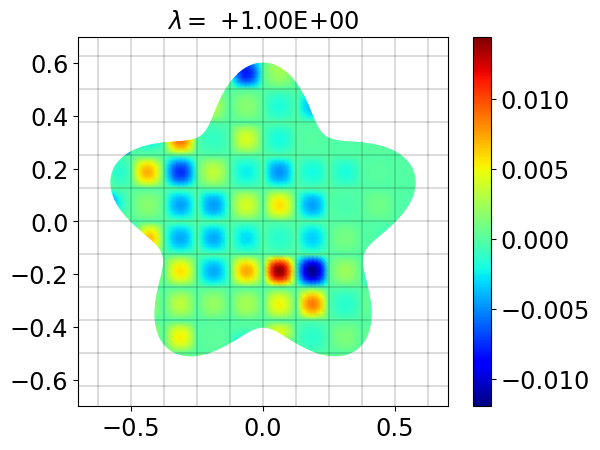}
   \caption{Volume mode (largest)\label{fig:AdditiveSchwarzLarge}}
  \end{subfigure}
 \caption{Typical spectrum and characteristic eigenmodes of an immersed system that is preconditioned by a double fixed point iteration with additive Schwarz.\label{fig:AdditiveSchwarz}}
\end{figure}

Efficiency and stability of additive Schwarz as a smoother requires adequate selection of the relaxation parameter $\gamma$. For the smoothing operations to efficiently reduce the error components in the directions of modes with eigenfunctions that can not be adequately captured on coarser grids, it is required that the relaxed eigenvalues $\gamma\lambda_i\left(\mathbf{M}^{-1}\mathbf{A}\right)$ corresponding to such eigenmodes are close to $1$. The requirement with regard to stability is formulated in \eqref{eq:stabilityMulti}, and states that $\lambda_{\rm min}\left(\mathbf{M}^{-1}\mathbf{A}\right) \geq 0$ and that $\gamma \lambda_{\rm max}\left(\mathbf{M}^{-1}\mathbf{A}\right) \leq 2$. The positivity of the eigenmodes follows from the symmetric positive definiteness of both $\mathbf{M}^{-1}$ and $\mathbf{A}$. Since the eigenvalues of $\mathbf{M}^{-1}\mathbf{A}$ coincide with the eigenvalues of $\mathbf{M}^{-\frac12}\mathbf{A}\mathbf{M}^{-\frac12}$, the eigenmodes can be bounded from above by:
\begin{subequations}\begin{align}
 \lambda_{\rm max} \left(\mathbf{M}^{-1}\mathbf{A}\right) & = \lambda_{\rm max} \left(\mathbf{M}^{-\frac12}\mathbf{A}\mathbf{M}^{-\frac12}\right) = \max\limits_{\mathbf{y}} \frac{\mathbf{y}^{\rm T}\mathbf{M}^{-\frac12}\mathbf{A}\mathbf{M}^{-\frac12}\mathbf{y}}{\mathbf{y}^{\rm T}\mathbf{y}} = \max\limits_{\mathbf{z}} \frac{\mathbf{z}^{\rm T}\mathbf{A}\mathbf{z}}{\mathbf{z}^{\rm T}\mathbf{M}\mathbf{z}} \\[.0em]
 & = \max\limits_{\mathbf{z}} \frac{\mathbf{z}^{\rm T}\mathbf{A}\mathbf{z}}{\min\limits_{\sum_{j=1}^{N} \mathbf{P}_j \mathbf{z}_j=\mathbf{z}} \sum_{j=1}^N \mathbf{z}_j^{\rm T}\mathbf{P}_j^{\rm T}\mathbf{A}\mathbf{P}_j\mathbf{z}_j} \label{eq:applyASLemma} \\[.3em]
 & = \max\limits_{\{\mathbf{z}_j\}_{j=1}^N} \frac{\left(\sum_{j=1}^N \mathbf{P}_j^{\rm T}\mathbf{z}_j^{\rm T}\right)\mathbf{A}\left(\sum_{j=1}^N \mathbf{P}_j\mathbf{z}_j\right)}{\sum_{j=1}^N \mathbf{z}_j^{\rm T}\mathbf{P}_j^{\rm T}\mathbf{A}\mathbf{P}_j\mathbf{z}_j} \\[.3em]
% & = \max\limits_{\{\mathbf{z}_j\}_{j=1}^N} \frac{\sum_{K_i} \left(\sum_{j=1}^N \mathbf{P}_j^{\rm T}\mathbf{z}_j^{\rm T}\right)\mathbf{A}^{K_i} \left(\sum_{j=1}^N \mathbf{P}_j\mathbf{z}_j\right)}{\sum_{K_i}\sum_{j=1}^N \mathbf{z}_j^{\rm T}\mathbf{P}_j^{\rm T}\mathbf{A}^{K_i}\mathbf{P}_j\mathbf{z}_j} \\[.3em]
 & \leq \max\limits_{K_i} \max\limits_{\{\mathbf{z}_j\}_{j=1}^{N_{K_i}}} \frac{\left(\sum_{j=1}^{N_{K_i}} \mathbf{P}_j^{\rm T}\mathbf{z}_j^{\rm T}\right)\mathbf{A}^{K_i} \left(\sum_{j=1}^{N_{K_i}} \mathbf{P}_j\mathbf{z}_j\right)}{\sum_{j=1}^{N_{K_i}} \mathbf{z}_j^{\rm T}\mathbf{P}_j^{\rm T}\mathbf{A}^{K_i}\mathbf{P}_j\mathbf{z}_j} \label{eq:applySummation} \\[.3em]
 & \leq \max\limits_{K_i} \max\limits_{\{\mathbf{z}_j\}_{j=1}^{N_{K_i}}} \frac{N_{K_i} \sum_{j=1}^{N_{K_i}} \mathbf{z}_j^{\rm T}\mathbf{P}_j^{\rm T}\mathbf{A}^{K_i}\mathbf{P}_j\mathbf{z}_j}{\sum_{j=1}^{N_{K_i}} \mathbf{z}_j^{\rm T}\mathbf{P}_j^{\rm T}\mathbf{A}^{K_i}\mathbf{P}_j\mathbf{z}_j} = \max\limits_{K_i} N_{K_i}, \label{eq:applyCS}
\end{align}\end{subequations}
with $\mathbf{A}^{K_i}\in\mathbb{R}^{n \times n}$ the part of system matrix $\mathbf{A}$ that results from the integration over element $K_i$, and $N_{K_i}$ denoting the number of blocks containing basis functions that are supported on this element. Note that in \eqref{eq:applyASLemma} the additive Schwarz lemma is applied, see e.g., \cite{Smith1996,Toselli2005}. In \eqref{eq:applySummation} the maximal quotient with system matrix $\mathbf{A}$ is replaced by the maximum over the quotients with the element contributions $\mathbf{A}^{K_i}$, and in \eqref{eq:applyCS} the Cauchy-Schwarz inequality is applied. For Lagrange basis functions this bound is observed to be considerably sharp, as volume basis functions are contained in $N_{K_i}$ blocks and form eigenfunctions similar to the one in Figure~\ref{fig:AdditiveSchwarzLarge}. These functions are not captured in the coarse grid correction, such that efficient smoothing requires $\gamma=N_{K_i}^{-1}=4^{-1}$ for a two-dimensional Lagrange basis. This relaxation parameter yields $\gamma \lambda_{\rm max}\left(\mathbf{M}^{-1}\mathbf{A}\right) \leq 1 < 2$, such that also the stability condition is satisfied.

\begin{figure} [pt]
   \begin{subfigure}{.49\textwidth}
   \vspace{.2cm}
   \includegraphics[height=5.1cm]{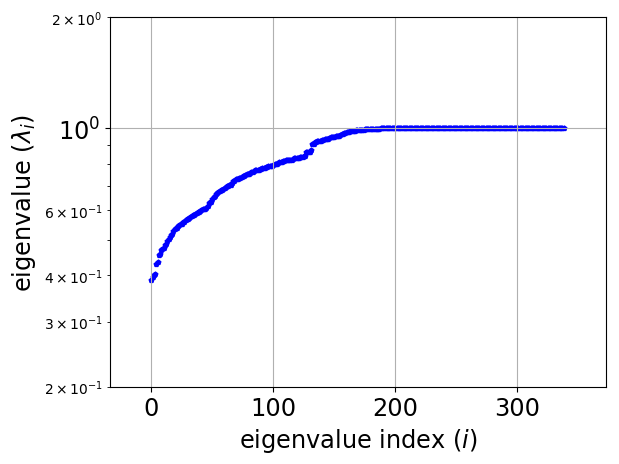}
   \caption{Eigenvalue spectrum\label{fig:AdditiveSchwarzSpectrumMultigrid}}
  \end{subfigure}
  \hfill
  \begin{subfigure}{.49\textwidth}
   \includegraphics[height=5cm]{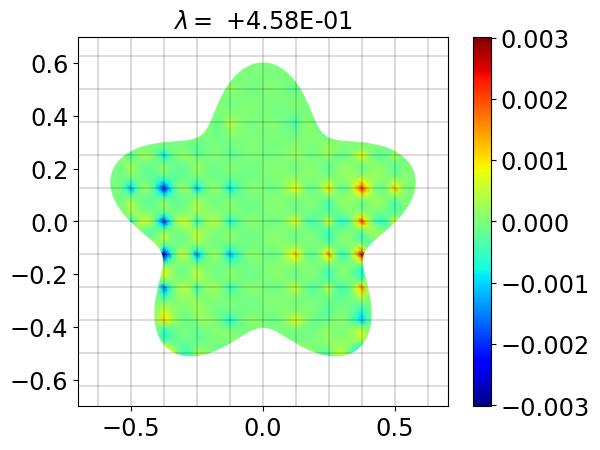}
   \vspace{.3cm}
   \caption{Nodal mode\label{fig:AdditiveSchwarzNodalMultigrid}}
  \end{subfigure}
 \caption{Spectrum and eigenmode of an immersed system that is preconditioned by a two-level V-cycle with additive Schwarz as smoother.\label{fig:AdditiveSchwarzMultigrid}}
\end{figure}

Figure~\ref{fig:AdditiveSchwarz} presents the spectrum and characteristic eigenmodes of the Laplace operator on the geometry in Figure~\ref{fig:domainStar} with a quadratic Lagrange basis that is preconditioned by a double fixed point iteration with additive Schwarz and the relaxation parameter $\gamma=\frac14$. The double fixed point iteration is applied such that later on these results can easily be related to the results with the V-cycle with additive Schwarz smoothing in Figure~\ref{fig:AdditiveSchwarzMultigrid}. Figure~\ref{fig:AdditiveSchwarzSmooth} shows that the smallest eigenvalue in the system with additive Schwarz does not correspond to an eigenfunction on a small cut element. Instead, the smallest eigenmode is similar to the usual smallest eigenmode with mesh-fitting methods, which can be considered as the smoothest possible mode satisfying the boundary conditions. Figure~\ref{fig:AdditiveSchwarzNodal} plots an example of an eigenfunction that almost entirely consists of nodal basis functions. The spectrum contains multiple of such modes, and these are important because these are the smallest eigenmodes that cannot be adequately captured by the coarse grid correction. Therefore, these nodal modes form the bottleneck for the condition number of immersed systems preconditioned by the V-cycle with additive Schwarz smoothing, as will follow in Figure~\ref{fig:AdditiveSchwarzMultigrid}. In fact, the eigenvalue of such modes can be clarified. Since nodal basis functions are in only $1$ index block, a fixed point iteration reduces the contribution of such functions by approximately a factor $1-\gamma=\frac34$. The double fixed point iteration therefore results in an eigenvalue of approximately $1-(1-\gamma)^2=\frac7{16}\approx0.438$. The largest eigenmodes in the spectrum have an eigenvalue of approximately $1$, and correspond to eigenfunctions consisting almost entirely out of basis functions that are only supported on 1 element and therefore are contained in 4 index blocks. On the interior this only involves volume basis functions, resulting in the volume mode in Figure~\ref{fig:AdditiveSchwarzLarge}. By virtue of the eigenvalue of approximately $1$, error components in the direction of the eigenvectors of these modes are effectively eliminated by the additive Schwarz smoother.

Based on the smallest eigenmodes in the spectrum with additive Schwarz, this technique is suitable as a smoother in a multigrid method for immersed finite elements. The results of this smoother in a two-level V-cycle are presented in Figure~\ref{fig:AdditiveSchwarzMultigrid}. It can be observed from the spectra in Figures~\ref{fig:AdditiveSchwarzSpectrum} and \ref{fig:AdditiveSchwarzSpectrumMultigrid} -- i.e., respectively without the coarse grid correction and with the coarse grid correction -- that the smooth eigenmodes are effectively resolved, such that a method is obtained that is robust to both cut elements and the grid size. It is noteworthy that the smallest eigenvalues with the multigrid method correspond to nodal modes, see Figure~\ref{fig:AdditiveSchwarzNodalMultigrid}. The limited effectiveness of the multigrid procedure for these modes derives from the fact that these modes are relatively insensitive to the smoothing operations, see Figure~\ref{fig:AdditiveSchwarzNodal}, and that these modes cannot be adequately captured on a coarser grid. With multigrid as a standalone solver, these modes would yield a convergence rate between $0.5$ and $0.6$. This rate can be improved by applying the multigrid cycle as a preconditioner in a Krylov subspace solver.

\begin{figure}[pt]
  \centering
  \begin{subfigure}{.49\textwidth}
   \centering
   \includegraphics[height=5.1cm]{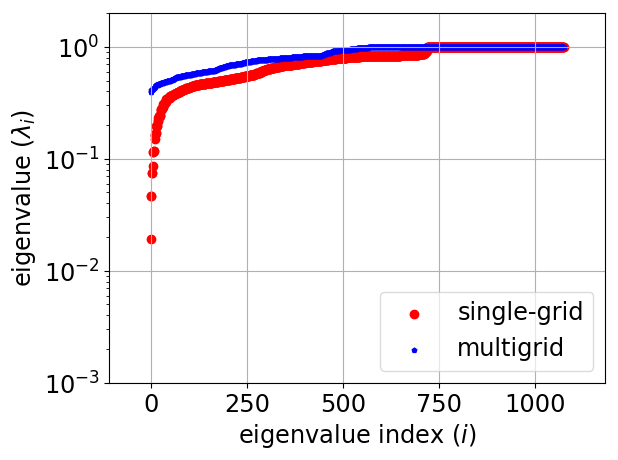}
   \caption{$32\times32$ elements\label{fig:AdditiveSchwarzSpectra32}}
  \end{subfigure}
  \hfill
  \begin{subfigure}{.49\textwidth}
   \centering
   \includegraphics[height=5.1cm]{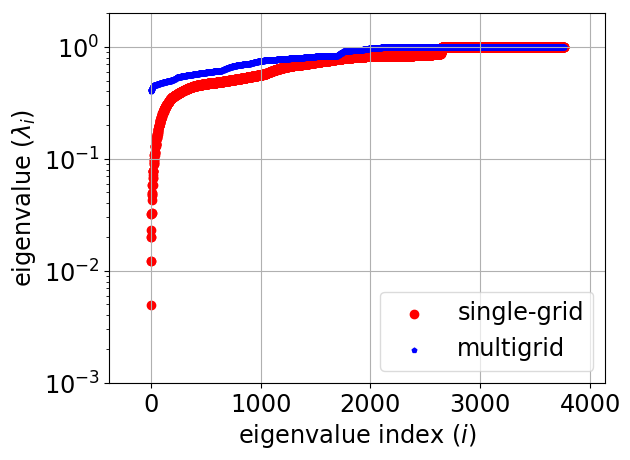}
   \caption{$64\times64$ elements\label{fig:AdditiveSchwarzSpectra64}}
  \end{subfigure}
 \caption{Spectra an immersed systems that are preconditioned by a double fixed point iteration with additive Schwarz (single-grid) and by a two-level V-cycle with additive Schwarz as smoother (multigrid).\label{fig:AdditiveSchwarzGridsize}}
\end{figure}

Figure~\ref{fig:AdditiveSchwarzGridsize} presents the spectra of the same problem with finer grids, preconditioned by a double fixed point iteration with additive Schwarz and preconditioned by the full V-cycle with additive Schwarz smoothing. Comparing these, and also the spectra with a grid of $16\times16$ elements in Figures~\ref{fig:AdditiveSchwarzSpectrum} and \ref{fig:AdditiveSchwarzSpectrumMultigrid}, conveys that systems without the coarse grid correction closely follow the grid-size dependence of the conditioning formulated in \eqref{eq:kappaMeshSize}, and that with the full V-cycle a conditioning is obtained that is independent of the grid size, with for this problem eigenvalues $0.4 < \lambda_{\rm min} < 0.5$ and $\lambda_{\rm max} = 1$.

While smoothing with additive Schwarz results in a conditioning that is independent of the grid size, the method is affected by the small relaxation parameter. This will be even more severe in three dimensions and for B-spline bases, which are supported on more elements and therefore require smaller and degree-dependent relaxation parameters, i.e., for B-splines $N^{K_i}=(p+1)^d$ with $p$ the spline degree and $d$ the number of dimensions. Therefore, the next section considers multiplicative Schwarz as a smoother, which is unconditionally stable \cite{Saad} and thereby circumvents the small relaxation parameter required for stability.

\clearpage

\subsubsection*{Multiplicative Schwarz}

Multiplicative Schwarz can be considered as the block equivalent of Gauss-Seidel. While with additive Schwarz all the locally inverted block matrices are applied to the same residual -- similar to the diagonal elements in Jacobi -- with multiplicative Schwarz the residual is updated after each block -- similar to the update of the residual after each diagonal element in Gauss-Seidel. Multiplicative Schwarz can be formulated by initializing the zero vector $\tilde{\mathbf{y}}^0 = \mathbf{0}$, defining the initial residual $\tilde{\mathbf{r}}^0 = \mathbf{r}$ and looping over the blocks $j \leq N$:
\begin{equation}\begin{aligned}
 \delta \tilde{\mathbf{y}}^j &= \mathbf{P}_j \mathbf{A}_j^{-1} \mathbf{P}_j^{\rm T} \tilde{\mathbf{r}}^{j-1}, \\
 \tilde{\mathbf{y}}^j &= \tilde{\mathbf{y}}^{j-1} + \delta \tilde{\mathbf{y}}^j, \\
 \tilde{\mathbf{r}}^j &= \tilde{\mathbf{r}}^{j-1} - \mathbf{A} \delta \tilde{\mathbf{y}}^j.
\end{aligned}\label{eq:multiplicativeSchwarz1}\end{equation}
The linear operator of multiplicative Schwarz is then defined as $\mathbf{M}^{-1}\mathbf{r} = \tilde{\mathbf{y}}^N$. Similar to Gauss-Seidel, multiplicative Schwarz does not require stabilization \cite{Saad}, which is the most important motivation to examine multiplicative Schwarz as a smoother in multigrid methods for immersed finite elements.

The linear operator induced by a fixed point iteration with multiplicative Schwarz is not symmetric. Therefore, it is important that in the post-smoothing the direction in the loop over the index blocks is reversed, to restore symmetry of the linear operator induced by the V-cycle. As can be observed in \eqref{eq:multiplicativeSchwarz1}, the application of multiplicative Schwarz requires updating the residual at every step, which is computationally expensive and impedes parallelization. Therefore, similar to standard Gauss-Seidel, the index blocks are ordered by a graph coloring algorithm \cite{Adams2003}. To this end, the blocks are divided in $C$ colors indicated by $c \leq C$. The basis functions in a block of a certain color do not intersect functions in a different block with the same color. Therefore, blocks of the same color are not affected by each other's update of the residual. As a result, the residual only needs to be updated after executing all blocks of a certain color. This reduces the computational cost and enables parallelization of the routine. The multiplicative Schwarz procedure with graph coloring can be formulated by initializing the zero vector $\tilde{\mathbf{z}}^0 = \mathbf{0}$, defining the initial residual again as $\tilde{\mathbf{r}}^0 = \mathbf{r}$ and looping over the colors $c \leq C$:
\begin{equation}\begin{aligned}
 \delta \tilde{\mathbf{z}}^c & = \sum_{j\in \mathcal{J}_c} \mathbf{P}_j \mathbf{A}_j^{-1} \mathbf{P}_j^{\rm T} \tilde{\mathbf{r}}^{c-1}, \\
 \tilde{\mathbf{y}}^c &= \tilde{\mathbf{y}}^{c-1} + \delta \tilde{\mathbf{y}}^c, \\
 \tilde{\mathbf{r}}^c &= \tilde{\mathbf{r}}^{c-1} - \mathbf{A} \delta \tilde{\mathbf{y}}^c,
\end{aligned}\end{equation}
with $\mathcal{J}_c$ denoting the set of blocks with color $c$. The linear operator of multiplicative Schwarz with graph coloring is defined as $\mathbf{M}^{-1}\mathbf{r} = \tilde{\mathbf{z}}^C$. The required minimal number of colors can be deduced by considering the overlap between the supports of basis functions. Because the supports of identically-colored basis functions are not allowed to intersect, for scalar problems and uniform grids the required number of blocks equals the number of elements on which basis functions are supported. For Lagrange bases, uniform B-spline bases and truncated hierarchical B-spline bases, the employed number of colors therefore amounts to $2^d$, $(p+1)^d$, and $M(p+1)^d$, respectively. For hierarchical bases the number of colors required for a uniform grid is multiplied by the number of hierarchical levels and for vector-valued problems the number of colors is multiplied by the number of components of the vector.

\begin{figure}
  \centering
  \begin{subfigure}{.49\textwidth}
   \centering
   \includegraphics[height=5.1cm]{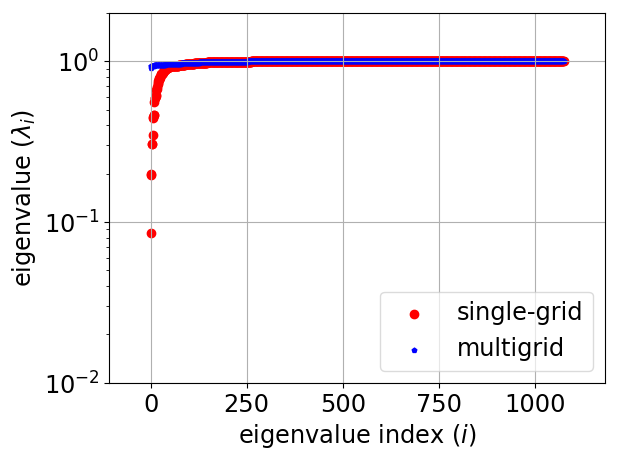}
   \caption{$32\times32$ elements\label{fig:MultiplicativeSchwarzSpectra32}}
  \end{subfigure}
  \hfill
  \begin{subfigure}{.49\textwidth}
   \centering
   \includegraphics[height=5.1cm]{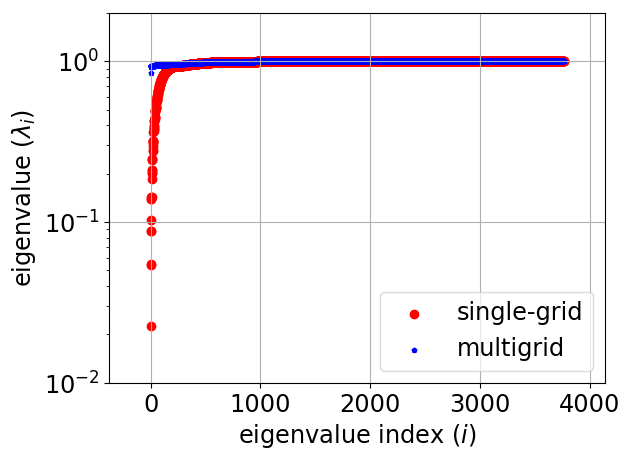}
   \caption{$64\times64$ elements\label{fig:MultiplicativeSchwarzSpectra64}}
  \end{subfigure}
 \caption{Spectra of immersed systems that are preconditioned by a double fixed point iteration with multiplicative Schwarz (single-grid) and by a two-level V-cycle with multiplicative Schwarz as smoother (multigrid).\label{fig:MultiplicativeSchwarzGridsize}}
 \vspace{-3mm}
\end{figure}

Figure~\ref{fig:MultiplicativeSchwarzGridsize} presents typical spectra of immersed systems, preconditioned by a symmetric double fixed point iteration with multiplicative Schwarz and preconditioned by a two-level V-cycle with multiplicative Schwarz smoothing. These results again pertain to the Laplace operator on the domain in Figure~\ref{fig:domainStar} with quadratic Lagrange basis functions. The single-grid results without the coarse grid correction show a largest eigenmode of $1$ and a mesh-dependent smallest eigenmode of order $\mathcal{O}(h^2)$, such that the condition number follows the relation for mesh-fitting systems in \eqref{eq:kappaMeshSize}. Similar to the results with additive Schwarz in Figure~\ref{fig:AdditiveSchwarzGridsize}, the results with the two-level V-cycle show a spectrum that is robust to both cut elements and the mesh size of the background grid. The conditioning with multiplicative Schwarz is much better and even nearly optimal, however, because it is not compromised by the small relaxation parameter that was required for the stability of additive Schwarz. Therefore, we conclude that multiplicative Schwarz is the most suitable smoothing procedure for multigrid methods for immersed systems. In view of its superior smoothing properties over additive Schwarz, in the numerical examples in Section~\ref{sec:resultsMulti} we restrict our considerations to multiplicative Schwarz.

\begin{remark}\label{rem:basisComparison}
The resulting preconditioning techniques for Lagrange bases and B-splines have much in common, but also small differences are to be noted. As mentioned previously, Lagrange bases require a considerably smaller number of colors in the implementation. Also, the solver for the Lagrange bases requires approximately $2.5$ times fewer iterations, as can be observed in the results in Section~\ref{sec:resultsMulti}. While both bases require approximately the same number of Schwarz blocks on the same grid, an advantage of the B-splines is that untrimmed basis functions result in a diagonal treatment, while untrimmed Lagrange basis functions still result in blocks with size $(2p-1)^d$. Also the considerably smaller number of degrees of freedom with B-splines is an aspect to consider. The methods for both bases, however, can be further optimized in regard of the computational efficiency. Therefore, we would like to emphasize that this contribution is intended to demonstrate the feasibility of multigrid methods for immersed finite elements, and not as a qualitative comparison between the different bases. Furthermore, while the computation time is linear with the number of degrees of freedom for both bases, the computation time is highly dependent on the implementation. For that reason, the computation time to solve the systems in Section~\ref{sec:resultsMulti} is not presented.
\end{remark}

\begin{remark}
 Similar to the Schwarz-type methods for immersed finite elements presented in \cite{Prenter2019} and \cite{John}, it is possible that block matrices contain eigenvalues of the order of the machine precision. Directly inverting such block matrices is unstable, as due to round-off errors these inverses can be inaccurate and can even contain negative eigenvalues for positive definite systems. As described in detail in \cite[Remark 3.3]{Prenter2019}, in case a block matrix contains an eigenvalue that is a factor $10^{16}$ smaller than the largest diagonal entry of the matrix, the basis function that is dominant in the corresponding eigenvector is removed from the block. Because this only pertains to basis functions with extremely small contributions, this does not affect the convergence of the iterative solver or the accuracy of the solution.
\end{remark}

\begin{remark}\label{rem:coarseGS}
 In principle, the ill-conditioning effects of small cut elements only need to be resolved on the finest grid, as the coarser levels are only required to resolve smooth components of the error. Diagonal smoothers therefore suffice on the coarser levels, and Schwarz-type smoothing is only required on the finest level. Applying Schwarz-type smoothing on all levels, however, retains the natural recursive character of the multigrid algorithm. For conciseness, we have therefore opted not to include results with diagonal smoothing on the coarser levels in the numerical results in Section~\ref{sec:resultsMulti}. % In principle, the ill-conditioning effects of small cut elements only need to be resolved on the finest grid, as the coarser levels are only required to resolve smooth components of the error. Diagonal smoothers therefore suffice on the coarser levels, such that Schwarz-type smoothing on these levels is actually not necessary. Applying Schwarz-type smoothing on all levels $\ell>1$ retains the natural recursive character of the algorithm, however, and for conciseness we have not included results with diagonal smoothing on the coarser levels in the numerical results in Section~\ref{sec:resultsMulti}. In principle, the ill-conditioning effects of small cut elements only need to be resolved on the finest grid, as the coarser levels are only required to be stable and to resolve smooth components of the error. Gauss-Seidel iterations satisfy both these requirements, such that Schwarz-type smoothing on these levels is actually not necessary. Applying multiplicative Schwarz on all levels $\ell>1$ is the most elegant solution, however, and for conciseness we have not included results with Gauss-Seidel on the coarser levels in the numerical results in Section~\ref{sec:resultsMulti}.
\end{remark}

\section{Numerical examples}\label{sec:resultsMulti}

In this section we assess the developed geometric multigrid preconditioning technique for immersed finite element methods on a range of numerical examples. In all these examples the pre-smoothing and post-smoothing operations consist of a single multiplicative Schwarz sweep as described in Section~\ref{sec:smoothers}. First, Section~\ref{sec:3d} considers three-dimensional elasticity problems on uniform grids with both B-splines and Lagrange basis functions. The goal of these simulations is to demonstrate the performance of the preconditioner on increasingly complex geometries. An aspect that is not covered in this section is the treatment of local mesh refinements. Therefore the preconditioner is applied to a level set based topology optimization problem with truncated hierarchical B-splines in Section~\ref{sec:topology}. All simulations are performed with quadratic discretizations. This is sufficient to demonstrate the difference between Lagrange basis functions and B-splines, and quadratic function spaces already result in severely ill-conditioned systems in case a dedicated treatment is not applied. The assembly of quadratic systems is considerably less expensive than systems of third degree or higher, however, enabling finer grids to be assessed.

\subsection{Linear elasticity problems}\label{sec:3d}

This section assesses the developed preconditioning technique in a conjugate gradient solver and discusses the obtained results. All problems are discretized with both quadratic Lagrange basis functions and quadratic B-splines on uniform grids. We first consider the deformation of a tooth-shaped domain subject to a distributed boundary traction. Next, an apple-shaped geometry containing two geometrically singular points is presented, which is deformed by a gravitational load. The third example considers a complex geometry of a $\mu$CT-scanned trabecular bone specimen, that is compressed by a prescribed displacement at the boundary. In this third example an effect regarding the number of levels in the preconditioner is observed. This effect is further investigated in the specifically designed test case in the last example, which is posed on a geometry in the shape of a triple helix.

\subsubsection*{Tooth-shaped geometry subject to a distributed boundary traction}

\begin{figure}[pt]
  \centering
  \begin{subfigure}{.49\textwidth}
   \centering
   \includegraphics[width=.9\textwidth]{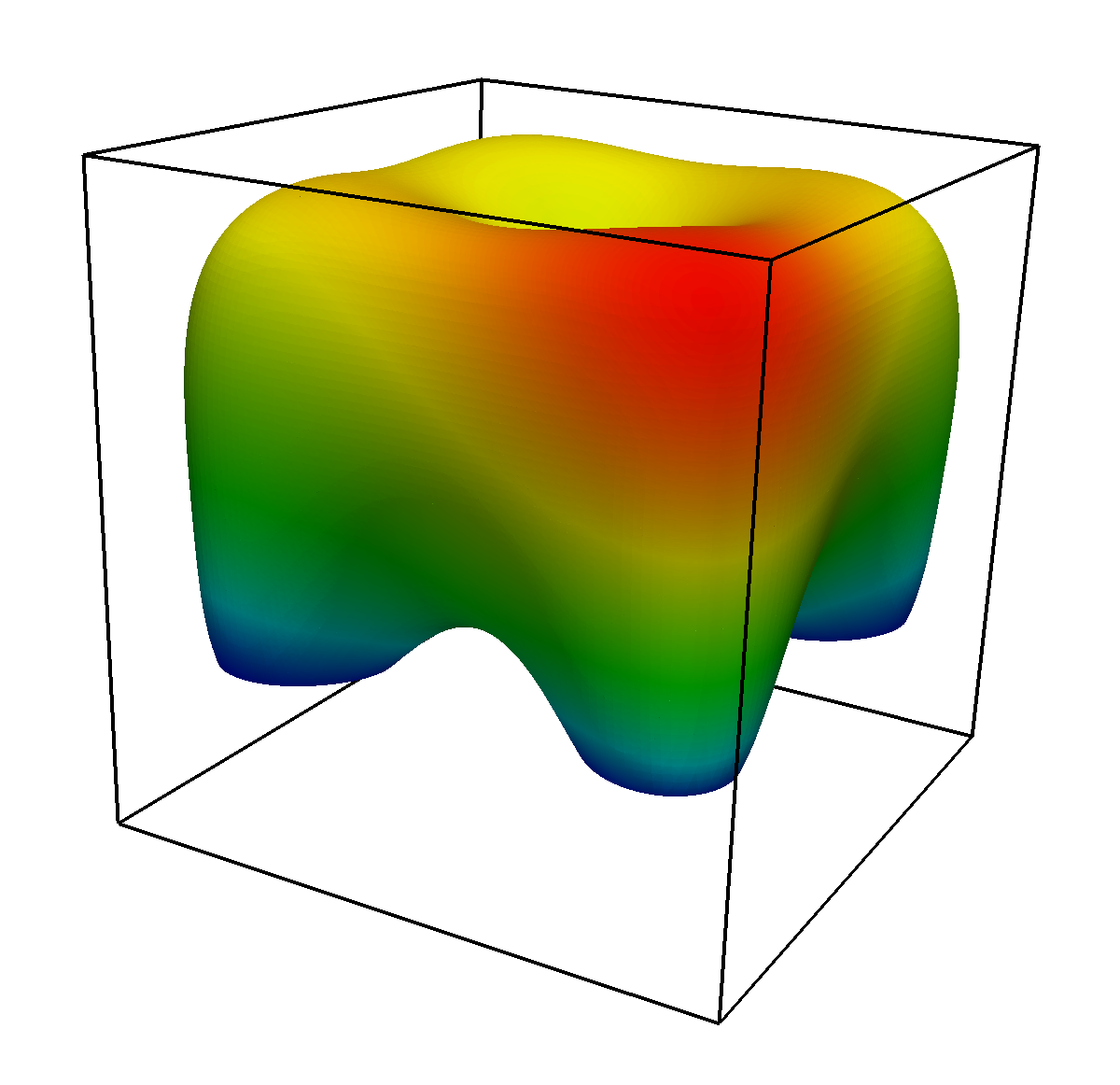} \\ 
   \includegraphics[width=.9\textwidth]{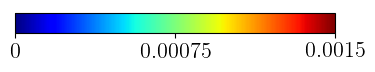}
   \caption{Displacement magnitude\label{fig:ToothDisplacement}}
  \end{subfigure}
  \hfill
  \begin{subfigure}{.49\textwidth}
   \centering
   \includegraphics[width=.9\textwidth]{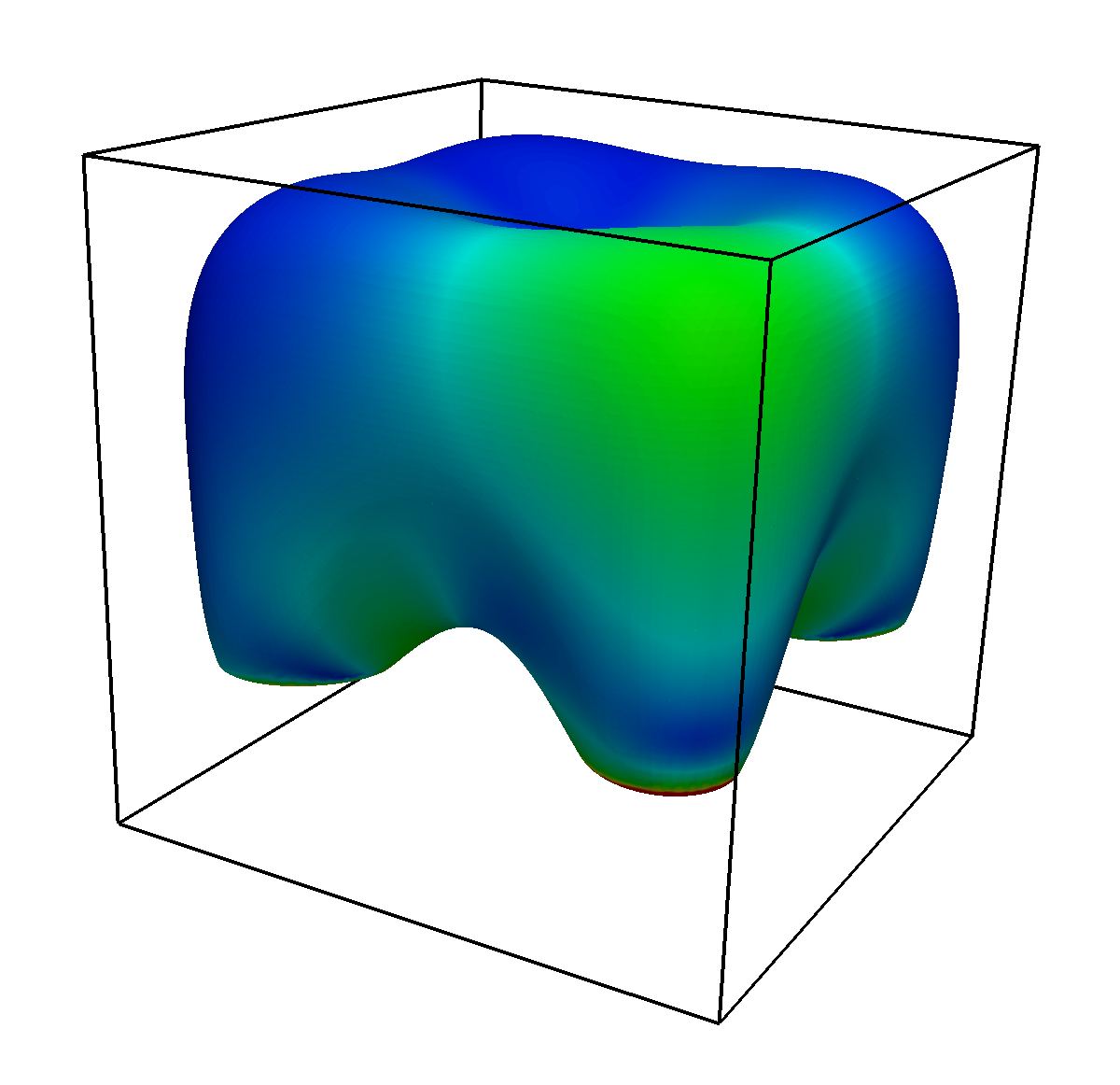} \\ 
   \includegraphics[width=.9\textwidth]{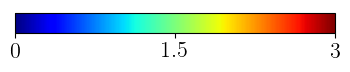}
   \caption{Frobenius norm of the stress tensor\label{fig:ToothStress}}
  \end{subfigure}
 \caption{Displacements and stresses of the tooth-shaped geometry. The displayed results are computed on a discretization with $80^3$ elements and quadratic B-spline basis functions.\label{fig:ToothSolution}}
\end{figure}

This first example considers a dimensionless problem posed on an immersed geometry with the shape of a tooth, see Figure~\ref{fig:ToothSolution}. The embedding domain is the cube $(-2,2)^3$, and the tooth-shaped geometry is obtained by trimming with the level set function:
\begin{equation}
 \psi_1(x,y,z) = 16 \left( 1 - \sum_{i=0}^6 e^{-\rho_i(x,y,z)} \right) - x^4 - y^4 - z^4,
\end{equation}
with:
\begin{equation*}\begin{array}{l c c c c c}
\rho_0(x,y,z) = & x^2 & + & y^2 & + & (z-2)^2, \\
\rho_1(x,y,z) = & (x-2)^2 & + & (y-2)^2 & + & \left(\frac{z+2}{2}\right)^2, \\
\rho_2(x,y,z) = & (x-2)^2 & + & (y+2)^2 & + & \left(\frac{z+2}{2}\right)^2, \\
\rho_3(x,y,z) = & (x+2)^2 & + & (y+2)^2 & + & \left(\frac{z+2}{2}\right)^2, \\
\rho_4(x,y,z) = & (x+2)^2 & + & (y-2)^2 & + & \left(\frac{z+2}{2}\right)^2, \\
\rho_5(x,y,z) = & x^2 & + & \left(\frac{y}{2}\right)^2 & + & \left(\frac{z+2}{2}\right)^2, \\
\rho_6(x,y,z) = & \left(\frac{x}{2}\right)^2 & + & y^2 & + & \left(\frac{z+2}{2}\right)^2.
\end{array}\end{equation*}
In this level set $\rho_0(x,y,z)$ creates the dent in the surface of the tooth, $\rho_1(x,y,z)$, $\rho_2(x,y,z)$, $\rho_3(x,y,z)$, and $\rho_4(x,y,z)$ create the dents in the sides of the tooth and the roots are obtained by $\rho_5(x,y,z)$ and $\rho_6(x,y,z)$. The tips of the roots below $z=-1$ are trimmed by a second trimming operation with the level set function:
\begin{equation}
 \psi_2(x,y,z) = z+1.
\end{equation}
This creates four surfaces on which homogeneous Dirichlet conditions are applied. On the rest of the boundary a normal traction is applied with the magnitude:
\begin{equation}
 \mathbf{g}^N(x,y,z) = - \mathbf{n} e^{ -\frac14 \left( (x-2)^2 + (y-2)^2 + (z-2)^2 \right) },
\end{equation}
which concentrates around the corner of the tooth. The Lam\'{e} parameters are set to $\lambda=\mu=10^3$, and the Dirichlet conditions are enforced by the penalty method with penalty parameters $\beta_h^\lambda=\beta_h^\mu= \frac2h$. The resulting displacements and stresses are shown in Figure~\ref{fig:ToothSolution}.

To demonstrate the robustness of the preconditioning technique to the number of elements, we discretize the embedding domain with grids of $20^3$, $40^3$, and $80^3$ elements. This results in, respectively, $129\cdot10^3$, $874\cdot10^3$, and $6.43\cdot10^6$ degrees of freedom (DOFs) with quadratic Lagrange basis functions and $21.6\cdot10^3$, $129\cdot10^3$, and $878\cdot10^3$ DOFs with quadratic B-splines. The integration depth as defined in Section~\ref{sec:methodMulti} is set to $2$ for the grid with $20^3$ elements, $1$ for the grid with $40^3$ elements, and $0$ for the grid with $80^3$ elements. This implies that the grid with $20^3$ elements is first partitioned by $2$ consecutive bisectioning operations before it is triangulated, the grid with $40^3$ elements is first partitioned by $1$ bisectioning operation before it is triangulated, and cut elements in the grid with $80^3$ elements are directly triangulated. Note that this results in identical integrated geometries for all grid sizes. The multigrid preconditioner is applied with $2$ and $3$ levels as defined in Section~\ref{sec:multigrid} for the grid with $20^3$ elements, $2$, $3$ and $4$ levels for the grid with $40^3$ elements and $2$, $3$, $4$ and $5$ levels for the grid with $80^3$ elements. These numbers of levels are chosen such that for all grid sizes, the largest number of levels results in a direct solution in the preconditioner for a system derived from $4^3$ elements.

\begin{figure}[pt]
  \centering
  \begin{subfigure}{.49\textwidth}
   \centering
   \includegraphics[width=.9\textwidth]{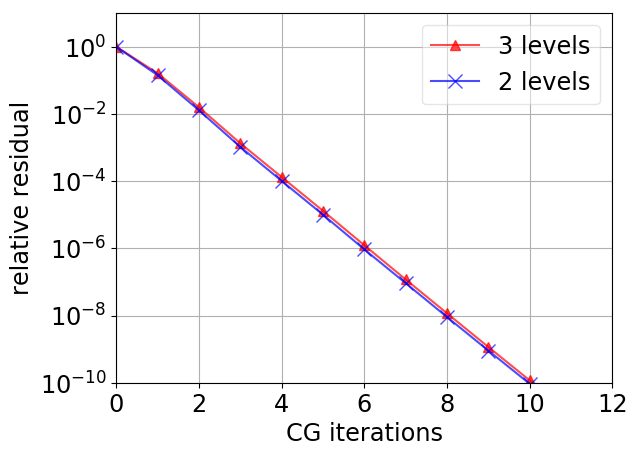}
   \caption{Lagrange basis with $20^3$ elements\label{fig:Tooth20Lagrange}}
  \end{subfigure}
  \hfill
  \begin{subfigure}{.49\textwidth}
   \centering
   \includegraphics[width=.9\textwidth]{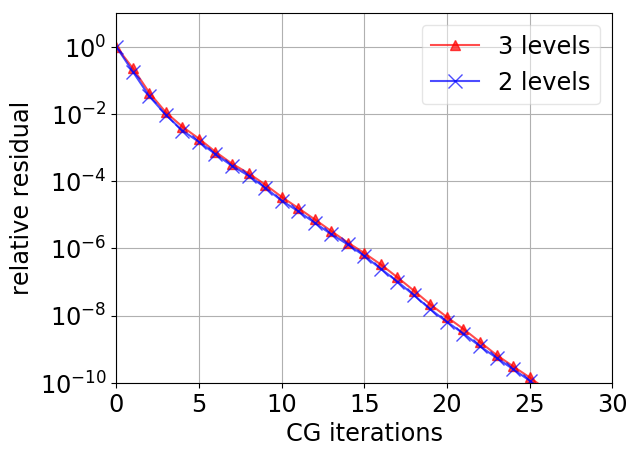}
   \caption{B-splines with $20^3$ elements\label{fig:Tooth20Splines}}
  \end{subfigure}
  \begin{subfigure}{.49\textwidth}
   \centering
   \includegraphics[width=.9\textwidth]{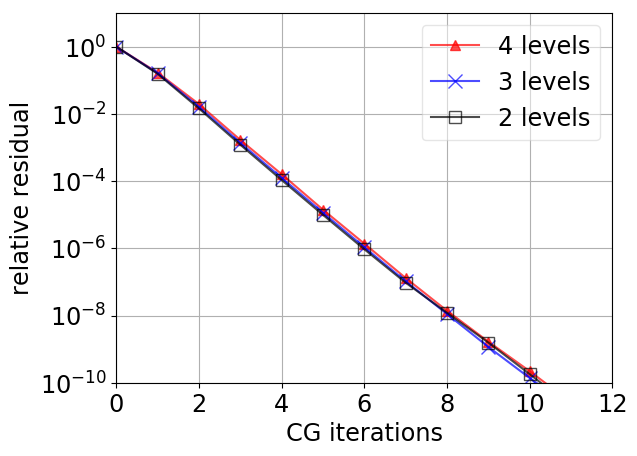}
   \caption{Lagrange basis with $40^3$ elements\label{fig:Tooth40Lagrange}}
  \end{subfigure}
  \hfill
  \begin{subfigure}{.49\textwidth}
   \centering
   \includegraphics[width=.9\textwidth]{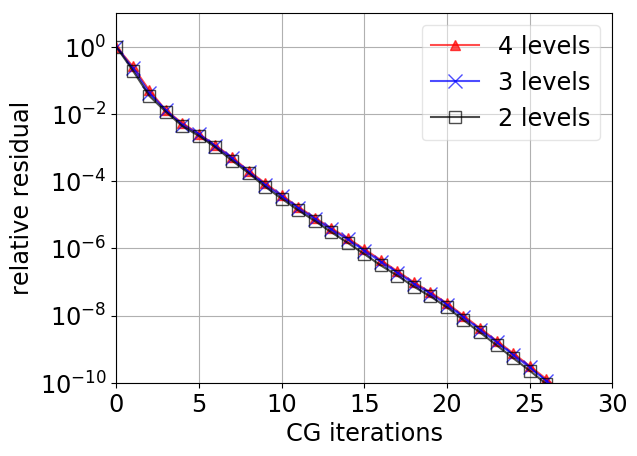}
   \caption{B-splines with $40^3$ elements\label{fig:Tooth40Splines}}
  \end{subfigure}
  \begin{subfigure}{.49\textwidth}
   \centering
   \includegraphics[width=.9\textwidth]{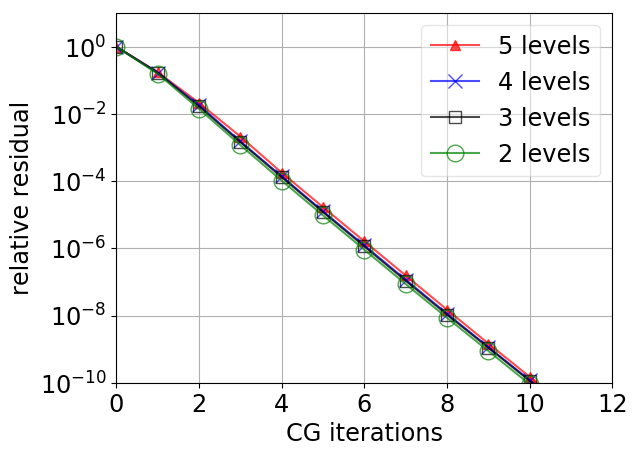}
   \caption{Lagrange basis with $80^3$ elements\label{fig:Tooth80Lagrange}}
  \end{subfigure}
  \hfill
  \begin{subfigure}{.49\textwidth}
   \centering
   \includegraphics[width=.9\textwidth]{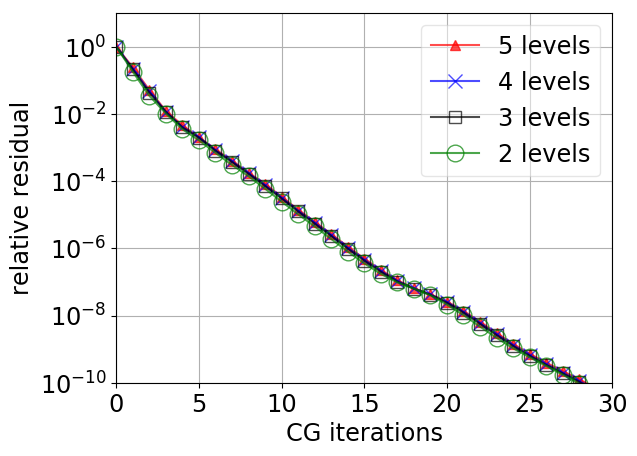}
   \caption{B-splines with $80^3$ elements\label{fig:Tooth80Splines}}
  \end{subfigure}
 \caption{Convergence plots of the problems on the tooth-shaped geometry.\label{fig:ToothResults}}
\end{figure}

The convergence of the multigrid-preconditioned conjugate gradient solver is plotted in Figure~\ref{fig:ToothResults}. The results demonstrate a convergence behavior that is virtually independent of the number of elements. Also, the number of iterations is nearly independent of the number of levels in the preconditioner. This can be attributed to the simple compact shape of the geometry, which is resolved well even on the coarsest grid of $4^3$ elements, such that with all numbers of levels the coarse grid corrections provide an effective approximation of the smooth eigenmodes. Finally, it can be observed that the systems with B-splines require more iterations than the systems with Lagrange basis functions. In this regard it should be noted, however, that for untrimmed basis functions the algorithm to select the Schwarz blocks yields a diagonal smoother for B-splines and a block smoother for Lagrange basis functions, \emph{cf.}\ the discussion on the computational cost of Lagrange and B-spline spaces in Remark~\ref{rem:basisComparison} that also considers the difference in the number of colors and the number of DOFs.

\subsubsection*{Apple-shaped geometry subject to a gravitational load}

\begin{figure}[pt]
  \centering
  \begin{subfigure}{.49\textwidth}
   \centering
   \includegraphics[width=.9\textwidth]{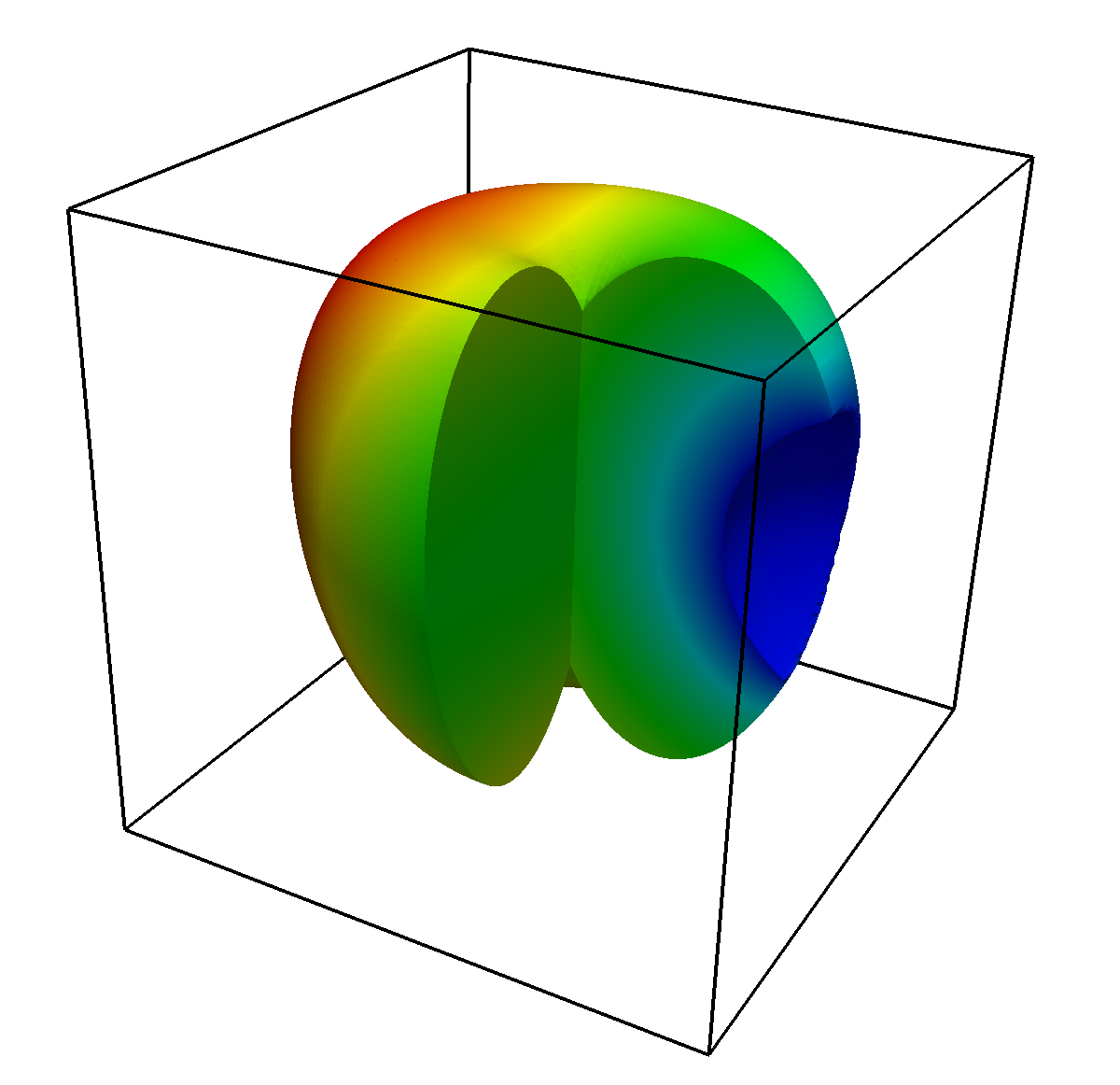} \\ 
   \includegraphics[width=.9\textwidth]{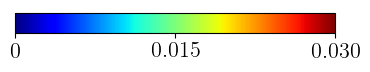}
   \caption{Displacement magnitude\label{fig:AppleDisplacement}}
  \end{subfigure}
  \hfill
  \begin{subfigure}{.49\textwidth}
   \centering
   \includegraphics[width=.9\textwidth]{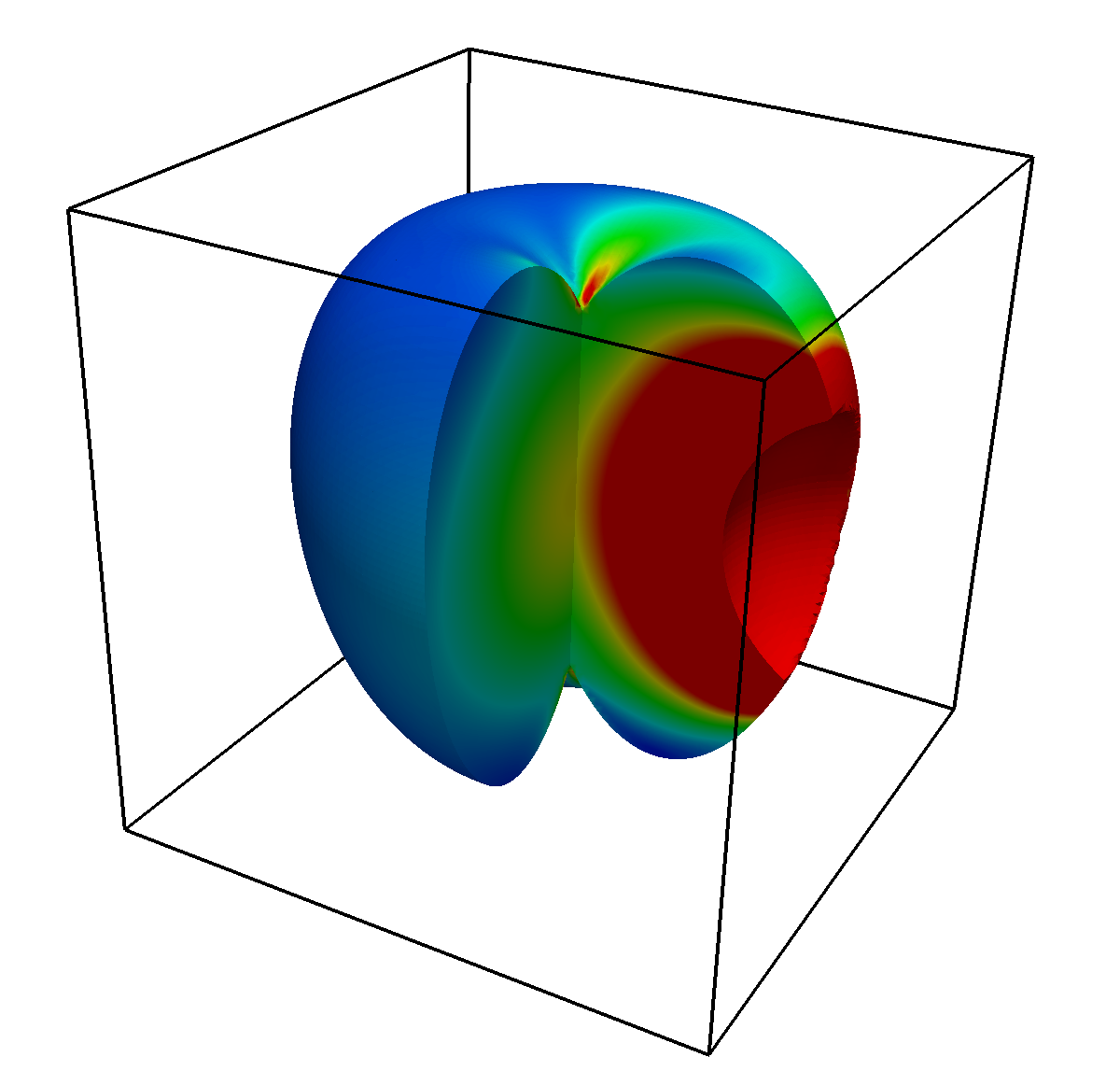} \\ 
   \includegraphics[width=.9\textwidth]{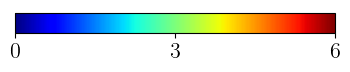}
   \caption{Frobenius norm of the stress tensor\label{fig:AppleStress}}
  \end{subfigure}
 \caption{Displacements and stresses in the geometry with the shape of an apple. The results are obtained with a B-spline basis on a grid of $80^3$ elements.\label{fig:AppleSolution}}
\end{figure}

This second example is designed to establish the suitability of the preconditioner for non-convex geometries with sharp reentrant corners where stress singularities are to be expected, such that the effectivity of multigrid methods is not generally evident \cite{Yserentant1986}. A dimensionless problem is posed on the geometry with the shape of an apple, see Figure~\ref{fig:AppleSolution}. The embedding domain is again the cube $(-2,2)^3$, and the shape of the apple is derived through a trimming operation with the level set function:
\begin{equation}
 \psi_1(x,y,z) = 1 - \left( \frac{z}{1.7} \right)^2 - \left( \frac{r(x,y)}{1+z/17} - 0.7 \right)^2,
\end{equation}
with $r(x,y)^2 = x^2 + y^2$. A second level set function models a bite being taken out of the apple:
\begin{equation}
 \psi_2(x,y,z) = \left( x - 2 \right)^2 + y^2 + z^2 - 1.
\end{equation}
Homogeneous Dirichlet conditions are imposed on the surface of the bite, and a volumetric load $\mathbf{f}=(0,0,-1)$ is applied to model gravity acting on the apple. The Lam\'{e} parameters are again set to $\lambda=\mu=10^3$ and the Dirichlet conditions are again enforced by the penalty method with parameters $\beta_h^\lambda=\beta_h^\mu= \frac2h$. Figure~\ref{fig:AppleSolution} displays the resulting displacements and stresses.

\begin{figure}[pt]
  \centering
  \begin{subfigure}{.49\textwidth}
   \centering
   \includegraphics[width=.9\textwidth]{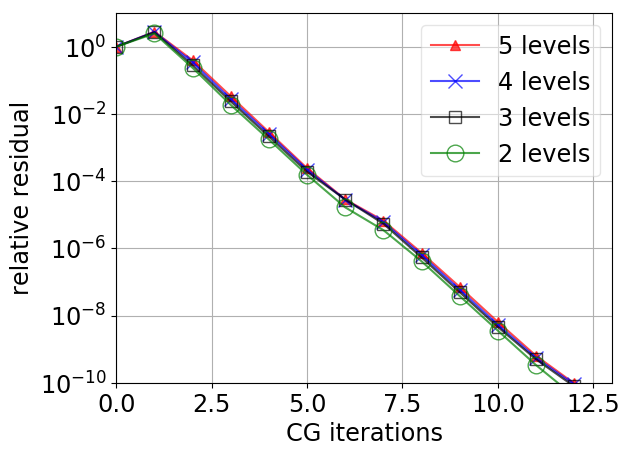}
   \caption{Lagrange\label{fig:AppleLagrange}}
  \end{subfigure}
  \hfill
  \begin{subfigure}{.49\textwidth}
   \centering
   \includegraphics[width=.9\textwidth]{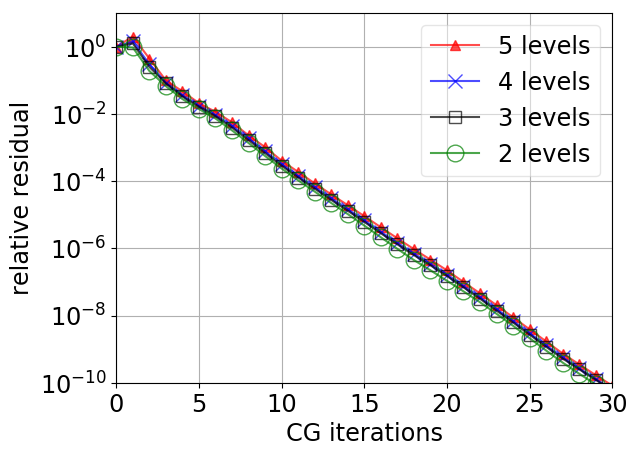}
   \caption{B-splines\label{fig:AppleSplines}}
  \end{subfigure}
 \caption{Convergence of the problem on the apple-shaped geometry with a grid of $80^3$ elements.\label{fig:AppleResults}}
\end{figure}

The problem is discretized with quadratic Lagrange basis functions and quadratic B-splines on a background grid with $80^3$ elements. This yields $4.87\cdot10^6$ DOFs supported in the physical domain with the Lagrange basis and $667\cdot10^3$ DOFs with the B-splines. The integration depth on the cut elements is set to $0$, implying that the the cut elements are directly triangulated and not first partitioned with a bisectioning operation. The convergence of the multigrid-preconditioned conjugate gradient solver with $2$, $3$, $4$, and $5$ levels for both bases is shown in Figure~\ref{fig:AppleResults}. The obtained numbers of iterations do not show an effect of the cusps in the geometry, and are similar to those for the tooth-shaped geometry. Also, the convergence is again virtually independent of the number of levels in the multigrid cycle.

\clearpage

\subsubsection*{Trabecular bone specimen loaded in compression}

\begin{figure}[pt]
 \centering
 \begin{subfigure}{.49\textwidth}
  \centering
  \includegraphics[width=.9\textwidth]{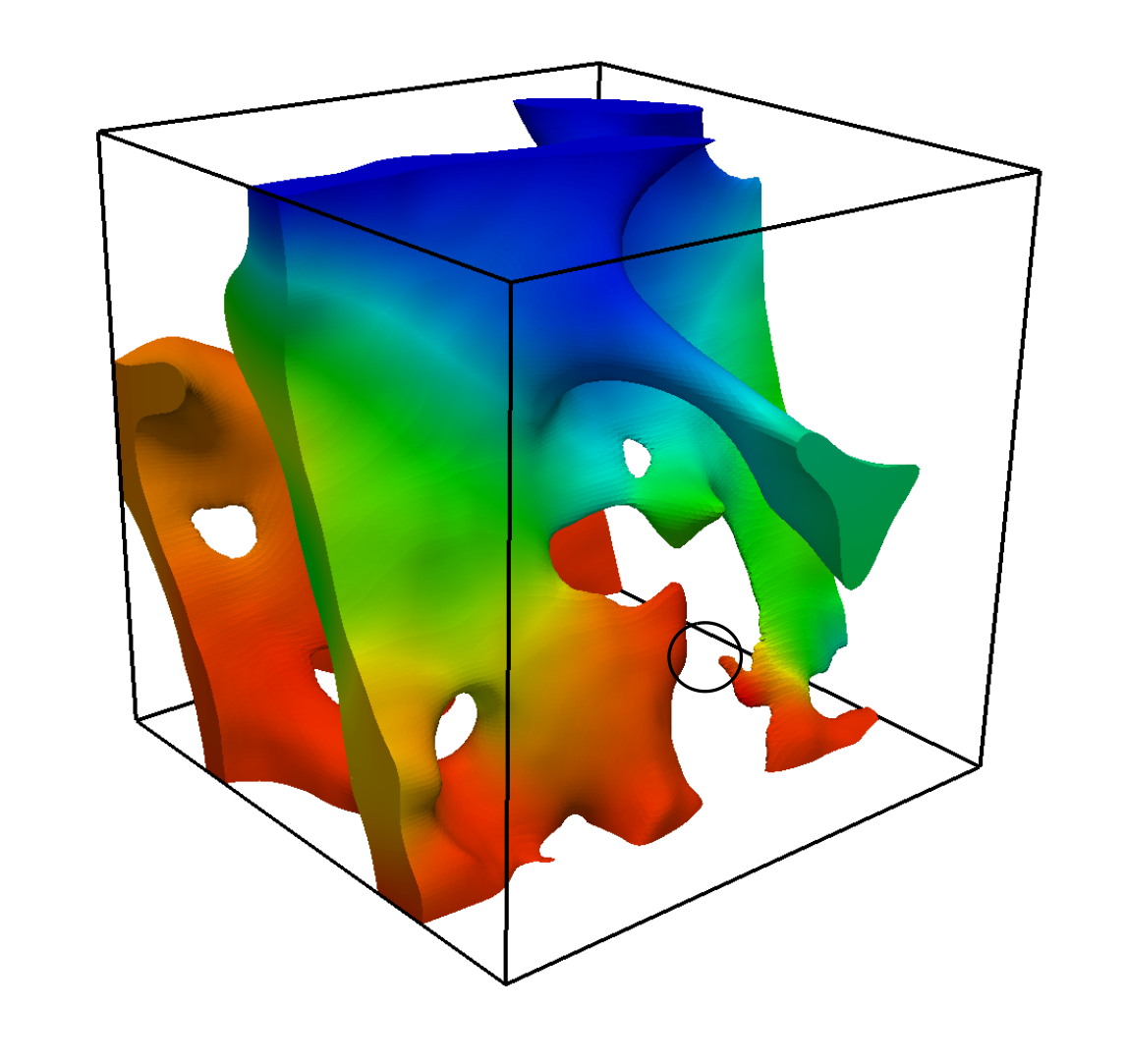} \\ 
  \includegraphics[width=.9\textwidth]{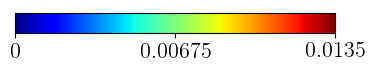}
  \caption{Displacement magnitude [mm]\label{fig:TrabecularDisplacement}}
 \end{subfigure}
 \hfill
 \begin{subfigure}{.49\textwidth}
  \centering
  \includegraphics[width=.9\textwidth]{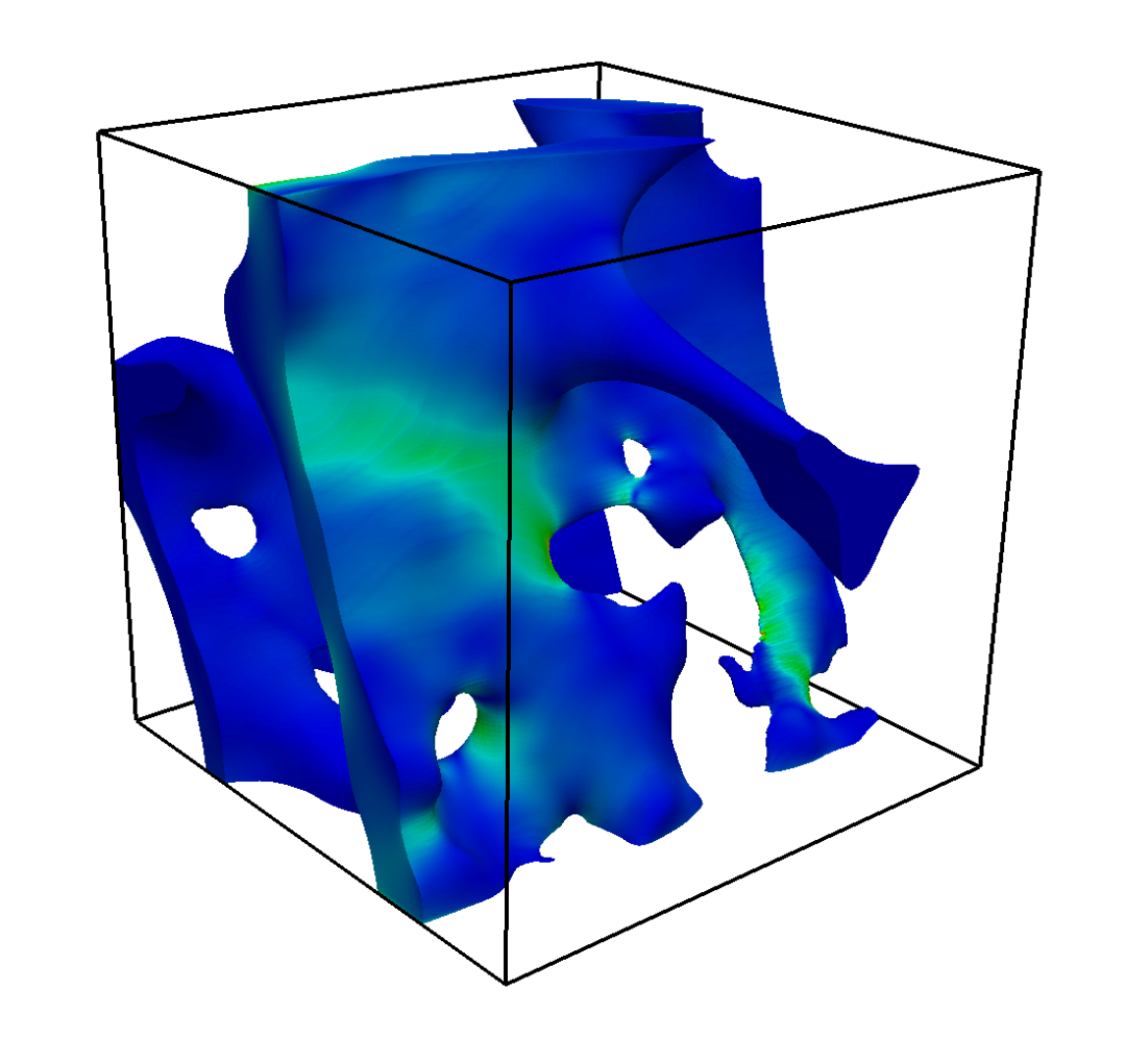} \\ 
  \includegraphics[width=.9\textwidth]{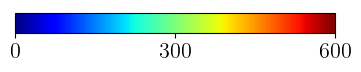}
  \caption{Frobenius norm of the stress tensor [MPa]\label{fig:TrabecularStress}}
 \end{subfigure}
 \caption{Solution of the elasticity problem on the three-dimensional trabecular bone geometry.\label{fig:TrabecularSolution}}
\end{figure}

This third three-dimensional test case considers the challenging geometry of a $\mu$CT-scanned trabecular bone specimen. This geometry was first presented in \cite{Verhoosel2015}, and is displayed in Figure~\ref{fig:TrabecularSolution}, together with the embedding domain of dimension $(0,1.28)^3$ mm$^3$. A linear elastic material model is employed with Young's modulus $E=10$ GPa and Poisson's ratio $\nu=0.3$. The specimen is compressed with an average uniaxial compressive strain of $1\%$, by imposing a homogeneous Dirichlet condition at the top boundary, and at the bottom boundary prescribing a normal displacement of $0.0128$ mm while constraining the tangential displacement. These boundary conditions are weakly enforced by the penalty method with penalty parameters $\beta_h^\lambda=\beta_h^\mu=\frac2h$.

We consider different grids on the embedding domain with $32^3$, $64^3$, and $128^3$ elements. The linear systems derived from these grids contain, respectively, $182\cdot10^3$, $1.03\cdot10^6$, and $6.65\cdot10^6$ DOFs with the quadratic Lagrange bases and $39.9\cdot10^3$, $189\cdot10^3$, and $1.05\cdot10^6$ DOFs with the quadratic B-splines. The integration depth is set to $2$ for the grid with $32^3$ elements, $1$ for the grid with $64^3$ elements, and $0$ for the grid with $128^3$ elements. The multigrid cycle by which the systems are preconditioned applies $2$ and $3$ levels for the discretizations with $32^3$ elements, $2$, $3$ and $4$ levels for the discretizations with $64^3$ elements and $2$, $3$, $4$ and $5$ levels for the discretizations with $128^3$ elements. With these numbers of levels, the preconditioner with the largest number of levels applies a direct solver to a system derived from $8^3$ elements for all three grid sizes, similar to the first example.

\begin{figure}[pt]
  \centering
  \begin{subfigure}{.49\textwidth}
   \centering
   \includegraphics[width=.9\textwidth]{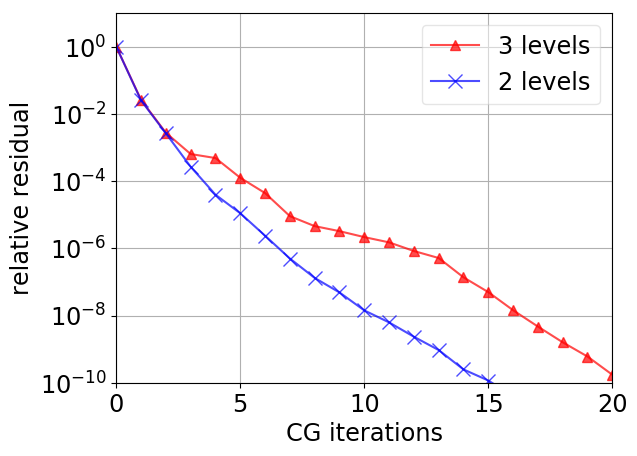}
   \caption{Lagrange basis with $32^3$ elements\label{fig:Trabecular32Lagrange}}
  \end{subfigure}
  \hfill
  \begin{subfigure}{.49\textwidth}
   \centering
   \includegraphics[width=.9\textwidth]{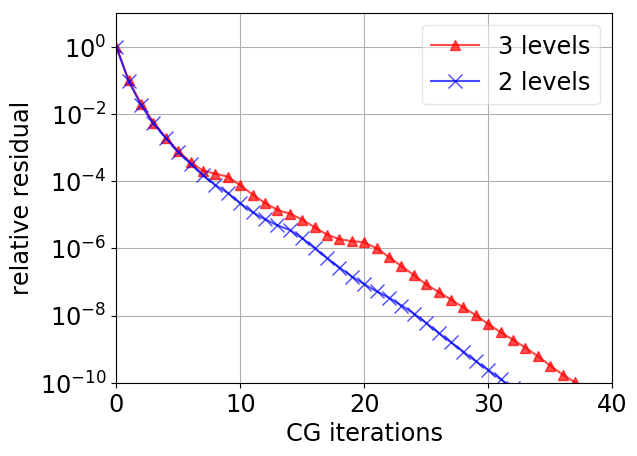}
   \caption{B-splines with $32^3$ elements\label{fig:Trabecular32Splines}}
  \end{subfigure}
  \begin{subfigure}{.49\textwidth}
   \centering
   \includegraphics[width=.9\textwidth]{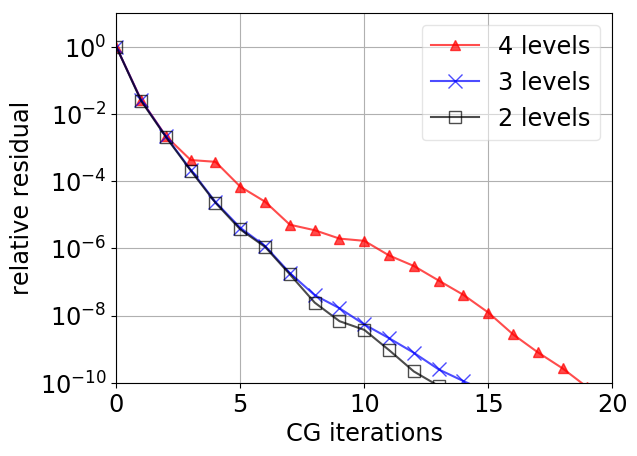}
   \caption{Lagrange basis with $64^3$ elements\label{fig:Trabecular64Lagrange}}
  \end{subfigure}
  \hfill
  \begin{subfigure}{.49\textwidth}
   \centering
   \includegraphics[width=.9\textwidth]{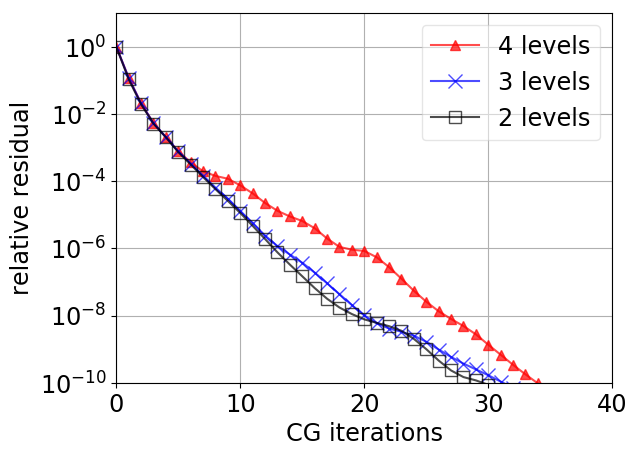}
   \caption{B-splines with $64^3$ elements\label{fig:Trabecular64Splines}}
  \end{subfigure}
  \begin{subfigure}{.49\textwidth}
   \centering
   \includegraphics[width=.9\textwidth]{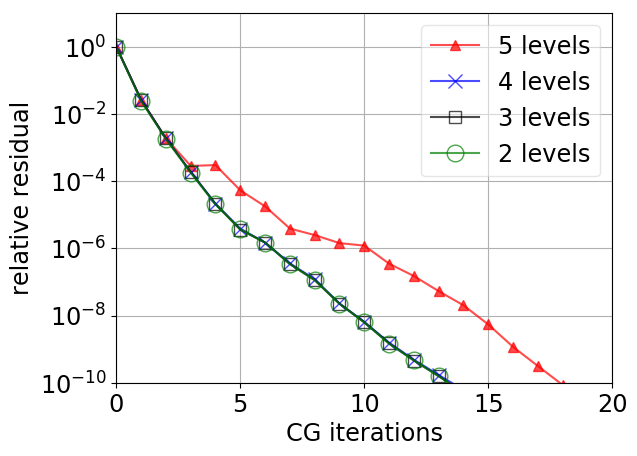}
   \caption{Lagrange basis with $128^3$ elements\label{fig:Trabecular128Lagrange}}
  \end{subfigure}
  \hfill
  \begin{subfigure}{.49\textwidth}
   \centering
   \includegraphics[width=.9\textwidth]{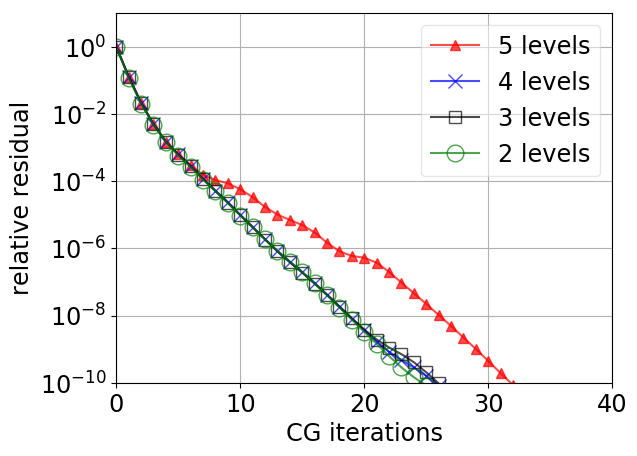}
   \caption{B-splines with $128^3$ elements\label{fig:Trabecular128Splines}}
  \end{subfigure}
 \caption{Convergence plots of the test case on the trabecular bone specimen.\label{fig:TrabecularResults}}
\end{figure}

\begin{figure}[pt]
 \centering
 \begin{subfigure}{.49\textwidth}
  \centering
  \includegraphics[width=.9\textwidth]{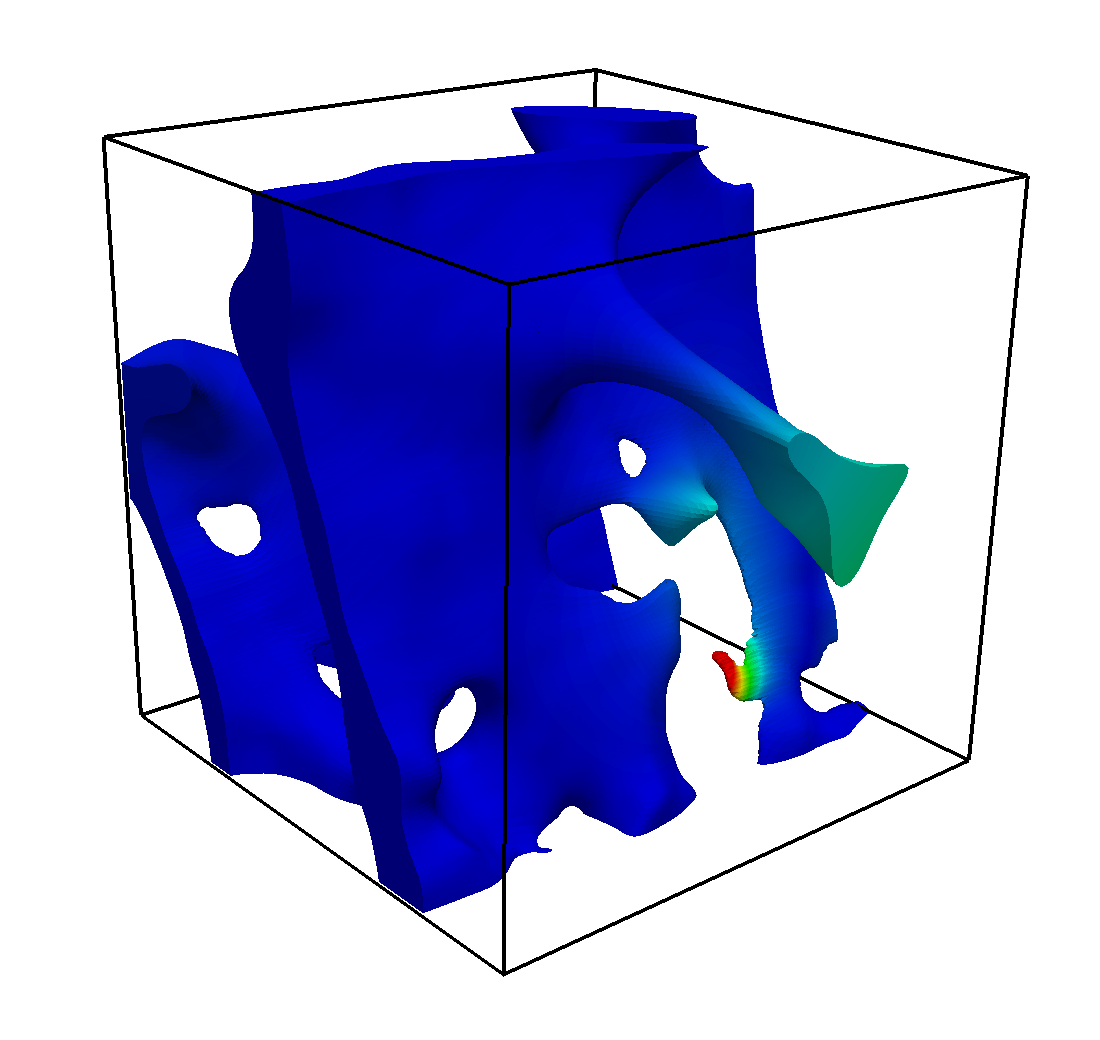} \\ 
  \includegraphics[width=.9\textwidth]{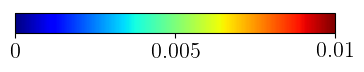}
  \caption{Displacement magnitude [mm]\label{fig:TrabecularSmallDisplacement}}
 \end{subfigure}
 \hfill
 \begin{subfigure}{.49\textwidth}
  \centering
  \includegraphics[width=.9\textwidth]{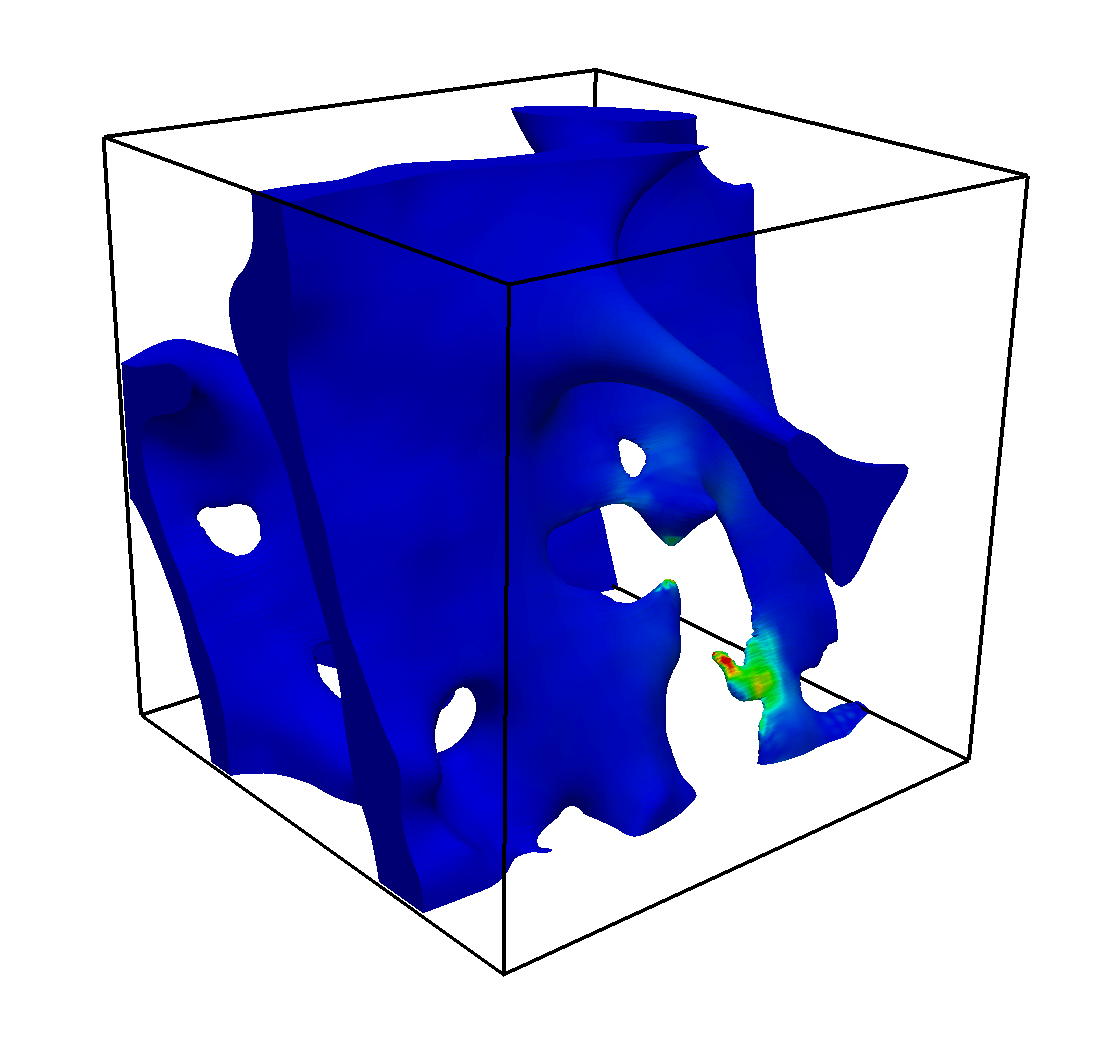} \\ 
  \includegraphics[width=.9\textwidth]{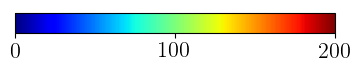}
  \caption{Frobenius norm of the stress tensor [MPa]\label{fig:TrabecularSmallStress}}
 \end{subfigure}
 \caption{Displacement and stress of the smallest eigenmode of the three-level preconditioned system with quadratic B-splines, formed on a grid of $32^3$ elements. The eigenmode shows an artificial stress state caused by a nonphysical coupling between disconnected parts of the geometry and the resulting deformation.\label{fig:TrabecularSmall}}
\end{figure}

The convergence of the systems with the different numbers of levels in the preconditioner is shown in Figure~\ref{fig:TrabecularResults}. It can be observed that the convergence of the systems with $32^3$ elements and $2$ levels, $64^3$ elements with $2$ and $3$ levels, and $128^3$ elements with $2$, $3$, and $4$ levels is very similar to that in the previous examples. For this test case, slightly more iterations are required. This is conjecturally connected with the complicated multiscale geometry, by which smaller geometric features are less adequately represented on coarse meshes. The convergence with the largest number of levels in the preconditioner, i.e., the preconditioners in which a direct solver is applied to a system with $8^3$ elements, is significantly slower. This is caused by underresolution and corresponding nonphysical behavior in the systems with $8^3$ elements, which was also observed in \cite{Verhoosel2015}. Both with the Lagrange and the B-spline bases, these coarse systems contain basis functions with a disjoint support in the physical domain, i.e., basis functions that cover the gap between disconnected parts of the geometry as e.g., in the circle in Figure~\ref{fig:TrabecularDisplacement}. Physically, disconnected parts of the geometry should be able to move freely with respect to each other. When a basis function is supported on both parts, however, these disconnected parts are artificially coupled in the numerical model. This results in very different stiffness properties between the systems with a mesh of $8^3$ elements and systems with meshes of $16^3$ elements and finer. Therefore, the coarse grid corrections with $8^3$ elements do not approximate the smooth eigenmodes of the finer meshes adequately, which retards the convergence of the multigrid-preconditioned iterative solver. A similar effect can be observed in nearly incompressible elasticity, where the convergence deteriorates when the coarse grids experience volumetric locking, as in e.g., \cite{Wieners2000}. To illustrate the artificial coupling, the eigenfunction with the smallest eigenvalue in the three-level preconditioned system with B-splines on a grid of $32^3$ elements is shown in Figure~\ref{fig:TrabecularSmall}. This function was obtained by $100$ iterations in a power algorithm, and it is clearly observable that this smallest preconditioned eigenmode contains peaks in the displacement and stress fields at points where basis functions of the coarse grid are supported on disconnected parts of the geometry. Such an artificial-coupling mode does not occur if the system is preconditioned with two levels, i.e., with a direct solver applied to the system with $16^3$ elements. While the artificial-coupling effect only moderately increases the number of iterations in this test case, and the iterative solver still converges in an acceptable number of iterations, it can not be ruled out that it can potentially hinder the convergence in other test cases more severely. Therefore, this effect is investigated in more detail in the next test case, which presents a geometry that is specifically designed to study this effect.

\subsubsection*{Triple helix loaded in compression}

\begin{figure}[pt]
 \centering
 \begin{subfigure}{.49\textwidth}
  \centering
  \includegraphics[width=.9\textwidth]{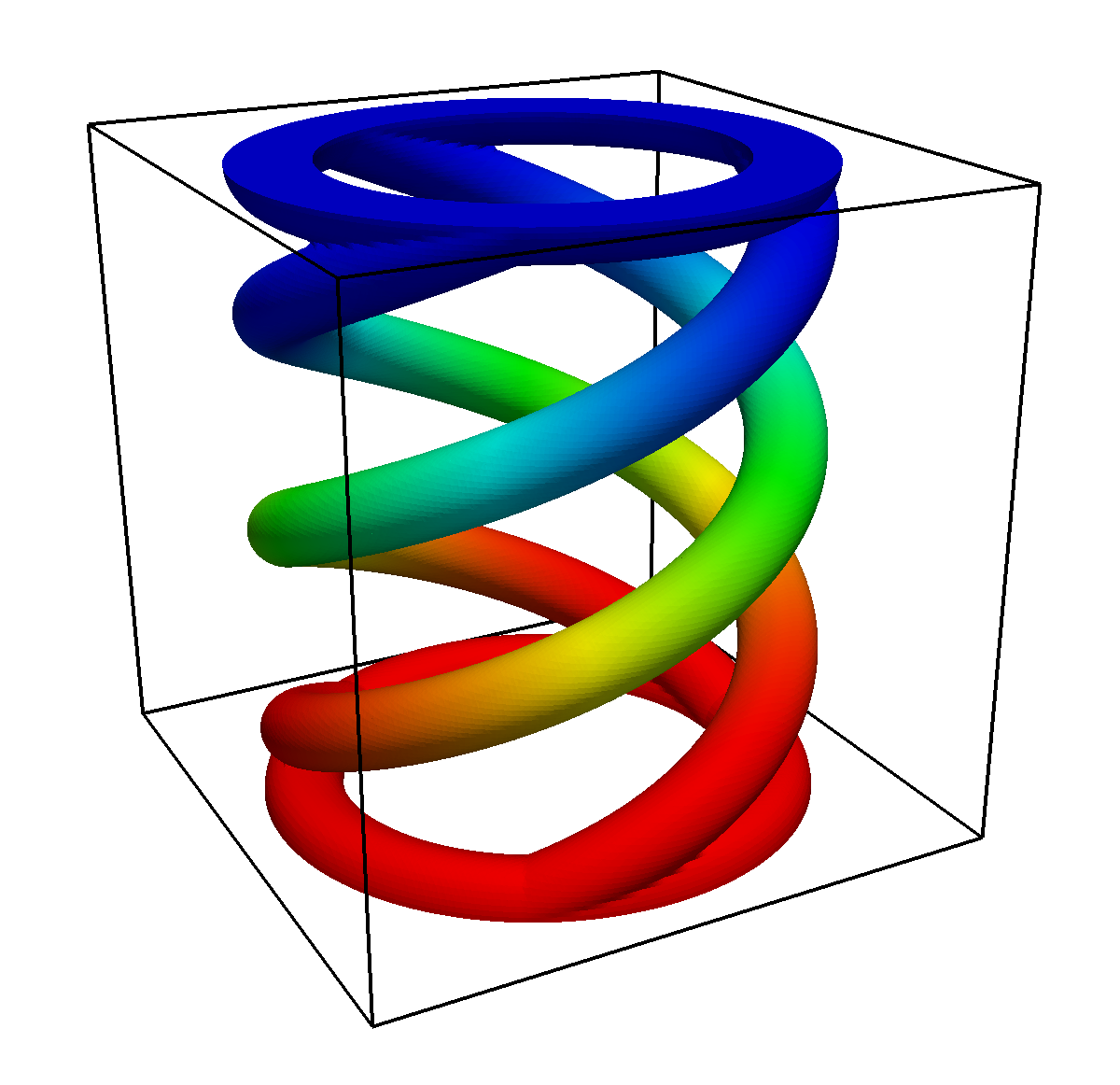} \\ 
  \includegraphics[width=.9\textwidth]{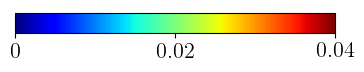}
  \caption{Displacement magnitude\label{fig:SpiralDisplacement}}
 \end{subfigure}
 \hfill
 \begin{subfigure}{.49\textwidth}
  \centering
  \includegraphics[width=.9\textwidth]{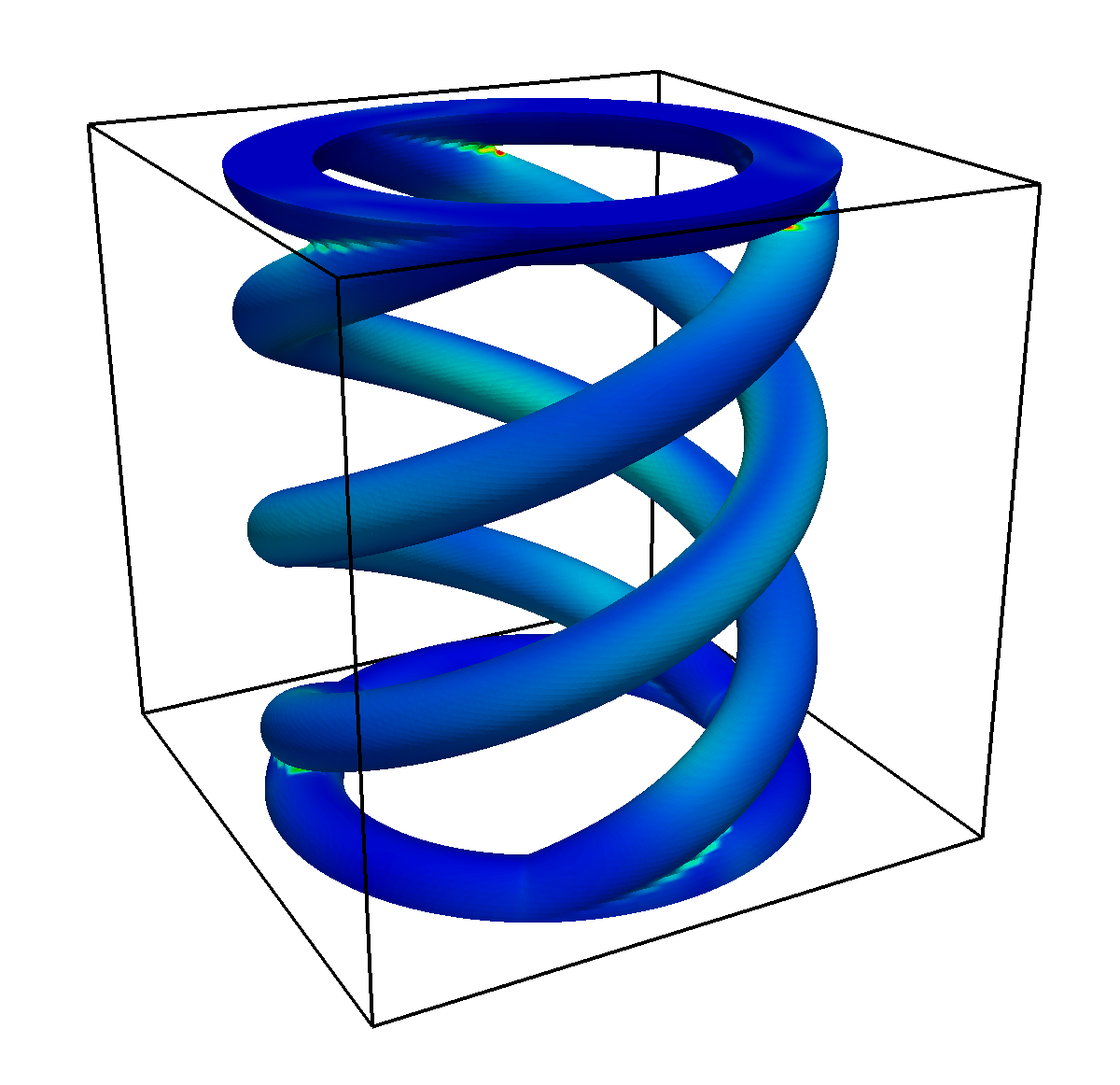} \\ 
  \includegraphics[width=.9\textwidth]{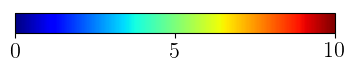}
  \caption{Frobenius norm of the stress tensor\label{fig:SpiralStress}}
 \end{subfigure}
 \caption{Solution of the elasticity problem on the triple-helix-geometry with B-spline basis functions on a grid of $128^3$ elements.\label{fig:SpiralSolution}}
 \vspace{-3mm}
\end{figure}

This final three-dimensional example is specifically designed to investigate the effect observed in the trabecular bone specimen, by which disconnected parts of the domain are artificially coupled on coarse grids. We consider a dimensionless problem on the geometry in Figure~\ref{fig:SpiralSolution}, which consists of a triple helix that is connected by a half ring at the top and bottom boundary. The embedding domain again consists of the cube $(-2,2)^3$. The helixes are obtained by rotations of the level set function:
\begin{equation}
 \psi_1(x,y,z) = r_{\rm inner}^2 - \left(r(x,y) - R_{\rm outer}\right)^2 - \left(z - 2 \frac{{\rm arctan2}(y,x)}{\pi}\right)^2,
\end{equation}
with $r_{\rm inner} = 0.25$, $R_{\rm outer} = 1.5$, and $r(x,y)^2 = x^2 + y^2$. For the half ring at the top and bottom boundary the level set function:
\begin{equation}
 \psi_2^\pm(x,y,z) = r_{\rm inner}^2 - \left(r - R_{\rm outer}\right)^2 - \left(z \pm 2\right)^2,
\end{equation}
is applied. The Lam\'{e} parameters are again set to $\lambda=\mu=10^3$. At the top boundary homogeneous Dirichlet conditions are imposed, and at the bottom boundary a normal displacement of $0.04$ prescribed. It should be noted that, as opposed to the bottom boundary condition on the trabecular bone specimen, the tangential displacement at the bottom boundary is not constrained. The boundary conditions are again weakly imposed by the penalty method with parameters $\beta_h^\lambda=\beta_h^\mu=\frac2h$. The solution with B-splines on a grid with $128^3$ elements is shown in Figure~\ref{fig:SpiralSolution}.

\renewcommand{\floatpagefraction}{1}

\begin{figure}[pt]
  \vspace{-3mm}
  \centering
  \begin{subfigure}{.49\textwidth}
   \centering
   \includegraphics[width=.9\textwidth]{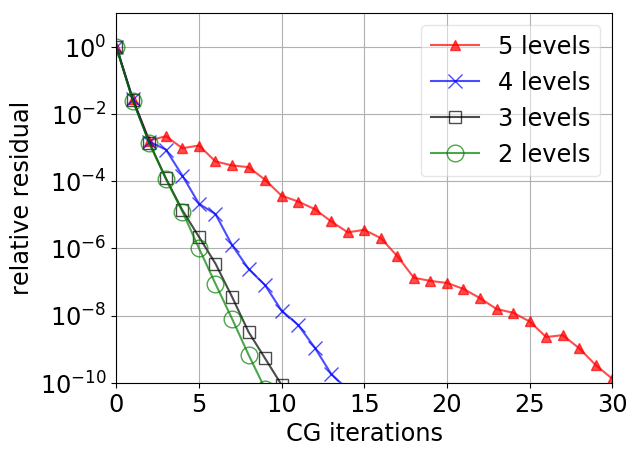}
 \vspace{-0mm}
   \caption{Lagrange\label{fig:SpiralLagrange}}
 \vspace{-0mm}
  \end{subfigure}
  \hfill
  \begin{subfigure}{.49\textwidth}
   \centering
   \includegraphics[width=.9\textwidth]{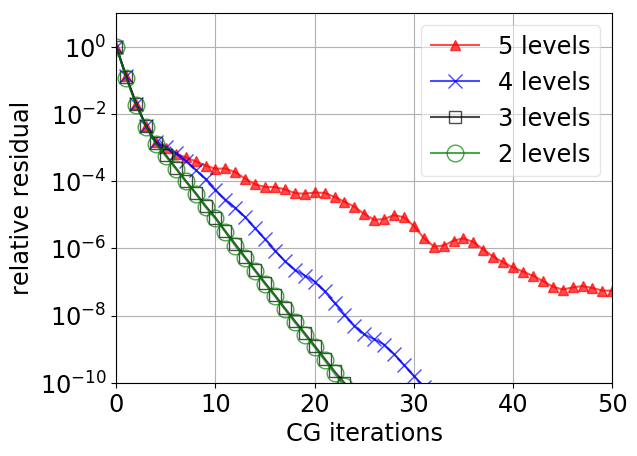}
 \vspace{-0mm}
   \caption{B-splines\label{fig:SpiralSplines}}
 \vspace{-0mm}
  \end{subfigure}
 \caption{Convergence of the triple helix geometry with $128^3$ elements.\label{fig:SpiralResults}}
 \vspace{-1mm}
\end{figure}

\begin{figure}[pt]
 \centering
 \begin{subfigure}{.49\textwidth}
  \centering
  \includegraphics[width=.9\textwidth]{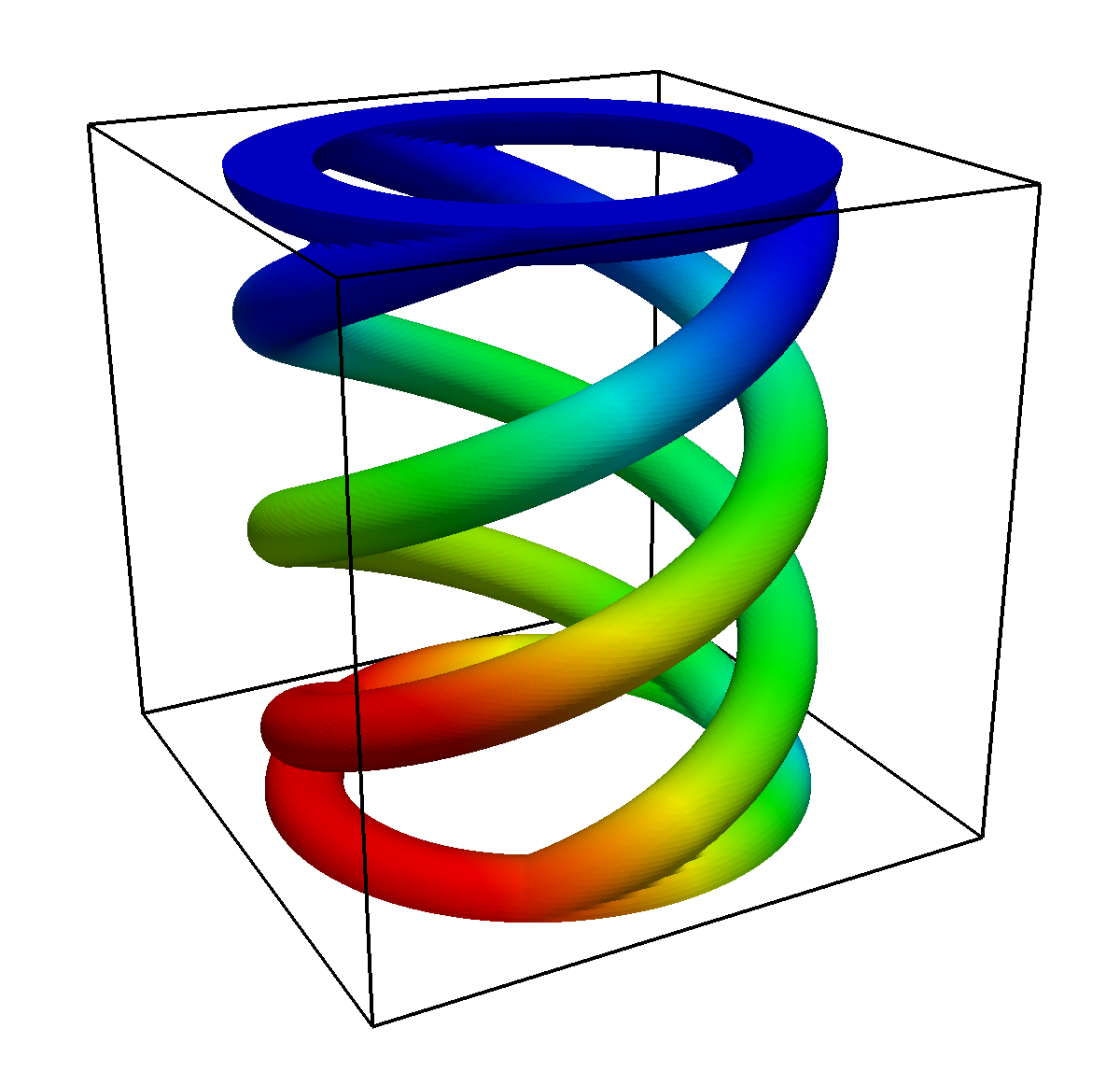} \\ 
 \vspace{-2mm}
  \includegraphics[width=.9\textwidth]{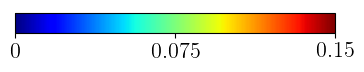}
 \vspace{-1mm}
  \caption{Displacement magnitude\label{fig:SpiralSmallDisplacement}}
 \vspace{-0mm}
 \end{subfigure}
 \hfill
 \begin{subfigure}{.49\textwidth}
  \centering
  \includegraphics[width=.9\textwidth]{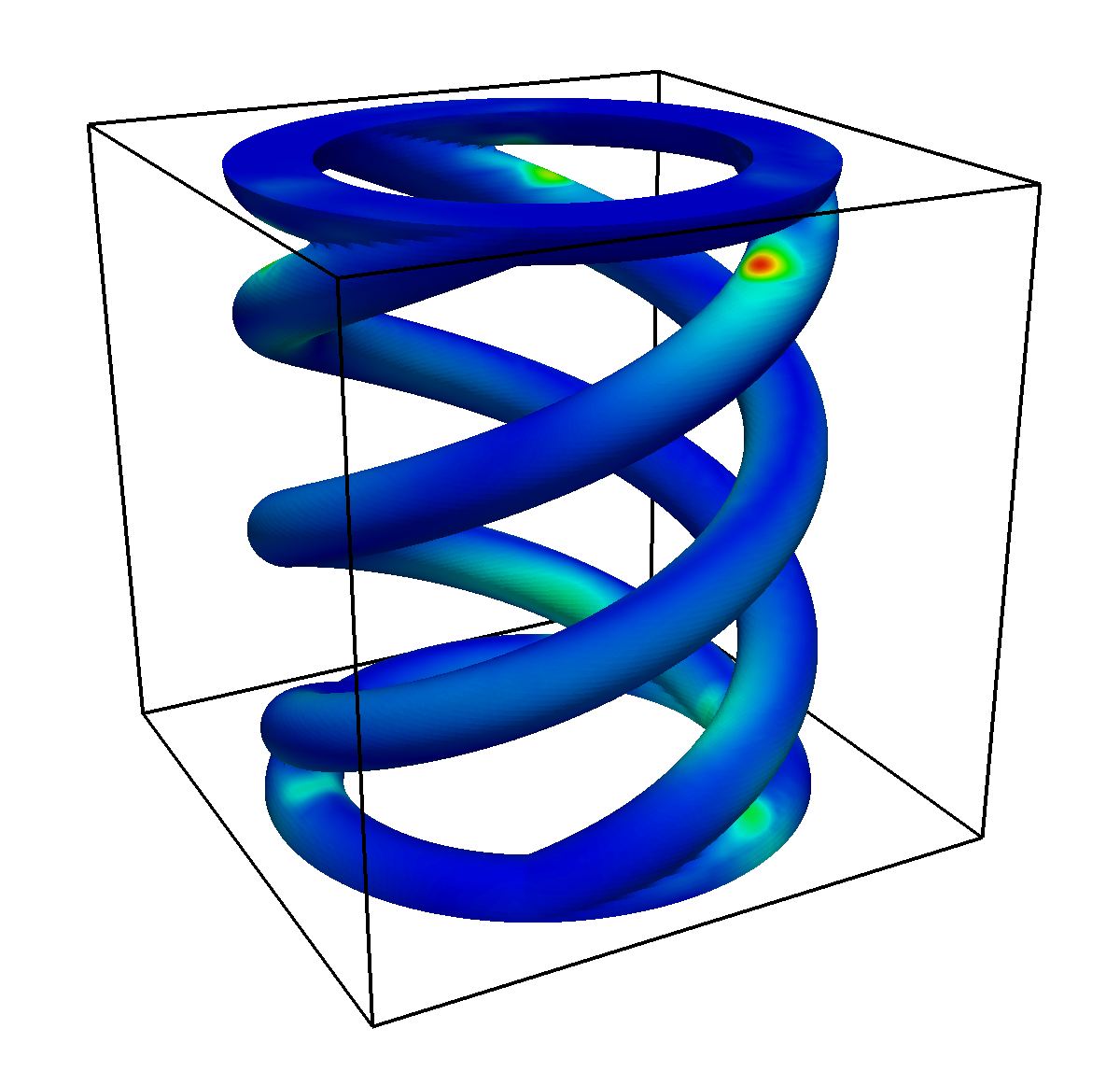} \\ 
 \vspace{-2mm}
  \includegraphics[width=.9\textwidth]{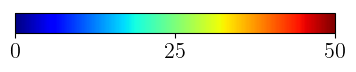}
 \vspace{-1mm}
  \caption{Frobenius norm of the stress tensor\label{fig:SpiralSmallStress}}
 \vspace{-0mm}
 \end{subfigure}
 \caption{Smallest eigenmode of a system with quadratic B-splines on a grid of $32^3$ elements that is preconditioned by a V-cycle with $3$ levels. The stress clearly indicates the artificial coupling that was also observed in the test case on the trabecular bone geometry. Furthermore, it is visible that the mismatch between the stiffness properties of the coarsest and the finer grids yields a smooth eigenmode that is not captured by the coarse grid correction.\label{fig:SpiralSmall}}
 \vspace{-2mm}
\end{figure}

The embedding domain is discretized with $128^3$ elements, yielding $6.96\cdot10^6$ and $1.11\cdot10^6$ DOFs in the quadratic Lagrange and B-spline bases, respectively. The integration depth is set to $0$, and the preconditioner applies $2$, $3$, $4$, and $5$ levels in the V-cycle. The convergence plots are shown in Figure~\ref{fig:SpiralResults}. It is clearly observable that the convergence is severely retarded with $5$ levels in the preconditioner. This could be anticipated, because the separation of the helixes is approximately equal to the mesh size on the coarsest level of $8^3$ elements. Hence, artificial coupling between the helixes will occur on the coarsest mesh, reducing the effectiveness of the coarse grid correction. To further illustrate this effect, Figure~\ref{fig:SpiralSmall} displays the smallest eigenmode in a system with B-splines on $32^2$ elements that is preconditioned with $3$ levels in the V-cycle, which behaves similarly as it also contains $8^3$ elements at the coarsest level. This mode is the equivalent of the smallest mode with the trabecular bone geometry in Figure~\ref{fig:TrabecularSmall}, as both these smallest eigenmodes clearly demonstrate the aforementioned artificial coupling effect. The preconditioner with $4$ levels -- for which approximately $2$ elements of the coarsest level fit between separate helixes -- also converges significantly slower. For the quadratic B-splines that span $3$ elements this can obviously be attributed to the artificial coupling effect. For the Lagrange basis functions that only span $2$ elements, it should be noted that basis functions covering the gap between a helix and a half ring can still artificially increase the stiffness of the connection between these. In the systems preconditioned by a V-cycle with $2$ and $3$ levels, the coarsest level contains $64^3$ and $32^3$ elements, respectively, such that the nonphysical coupling effect is not observed and the convergence is similar to that of the previous examples. We expect that the observed deterioration of the convergence behavior can be mitigated by a dedicated coarsening algorithm, which prevents disjoint supports in coarse basis functions by applying local refinements or XFEM-type enrichments on the coarse grid. In this regard the work presented in \cite{Hiriyur2012} is noteworthy, since it considers an algebraic multigrid approach in which a similar form of artificial coupling through a crack on the coarse level is precluded.


\clearpage

\renewcommand{\floatpagefraction}{.6}

\subsection{A level set based topology optimization problem with truncated hierarchical B-splines}\label{sec:topology}

\begin{figure}[pt]
  \centering
  \includegraphics[height=3.5cm]{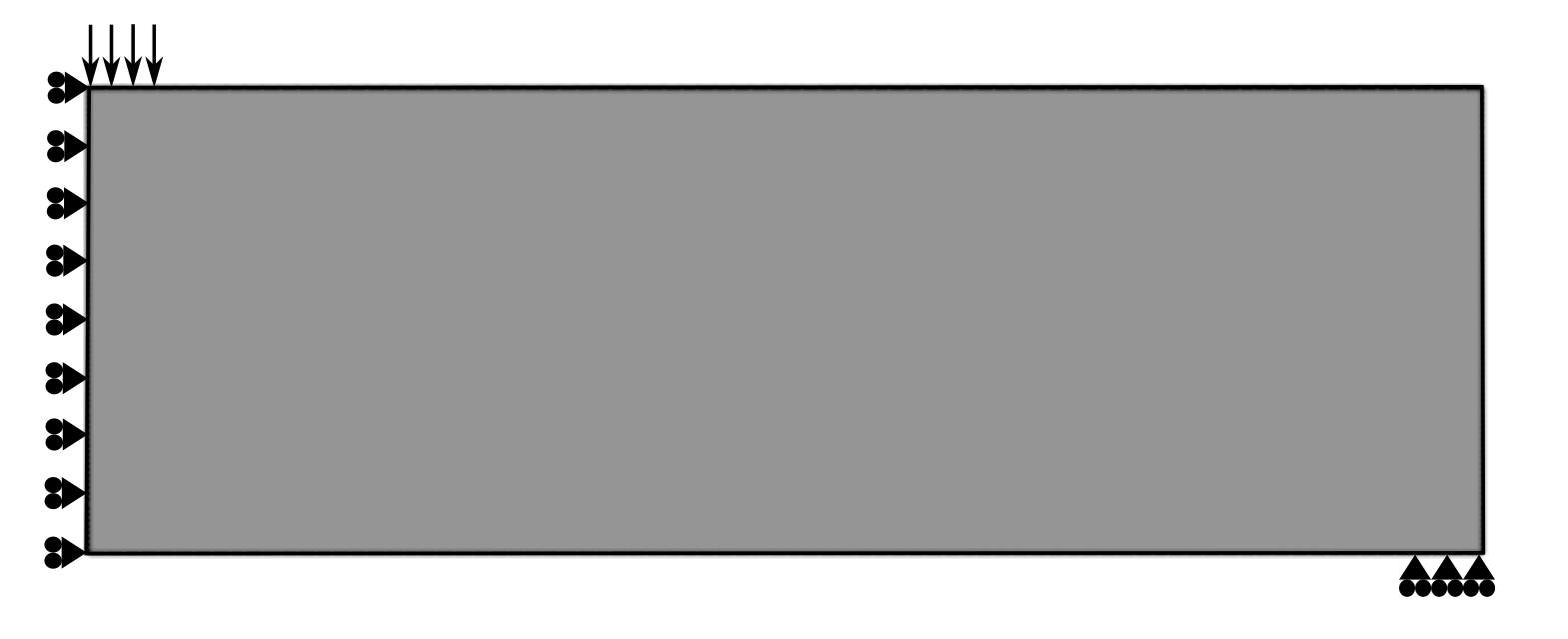}
  \caption{Design space of $3\times1$ and boundary conditions of the topology optimization problem. The vertical displacement is constrained at the right bottom, and a vertical load is applied at the left top. The horizontal displacement is constrained at the left boundary, which imposes a symmetry condition, such that the problem setup resembles a simply supported beam.\label{fig:topologyProblem}}
\end{figure}

\begin{figure}[pt]
  \vspace{-2cm}
  \centering
  \begin{subfigure}{.99\textwidth}
   \centering
   \includegraphics[height=2.5cm]{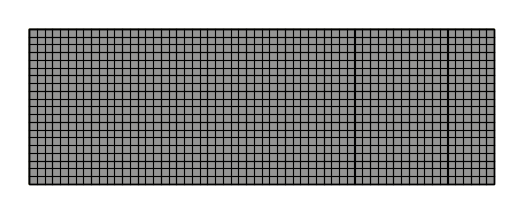}
   \hspace{.5cm}
   \includegraphics[height=2.5cm]{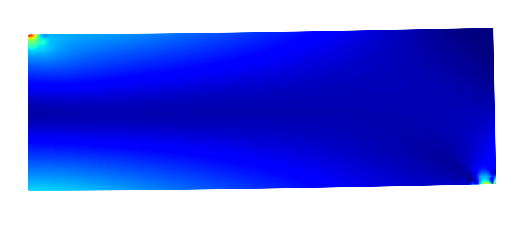}
   \caption{Step 0\label{fig:step001}}
  \end{subfigure}
  \begin{subfigure}{.99\textwidth}
   \centering
   \includegraphics[height=2.5cm]{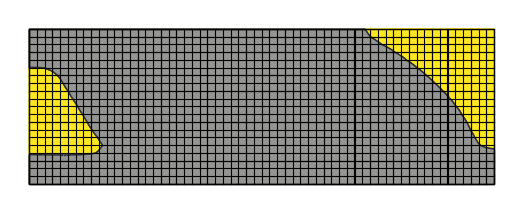}
   \hspace{.5cm}
   \includegraphics[height=2.5cm]{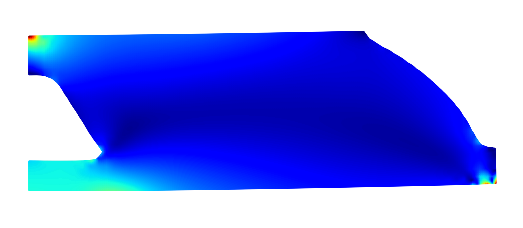}
   \caption{Step 30\label{fig:step030}}
  \end{subfigure}
  \begin{subfigure}{.99\textwidth}
   \centering
   \includegraphics[height=2.5cm]{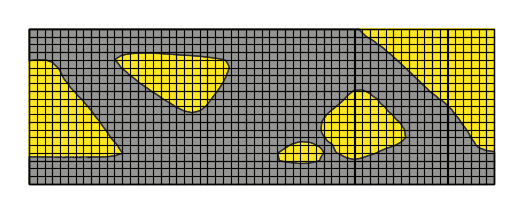}
   \hspace{.5cm}
   \includegraphics[height=2.5cm]{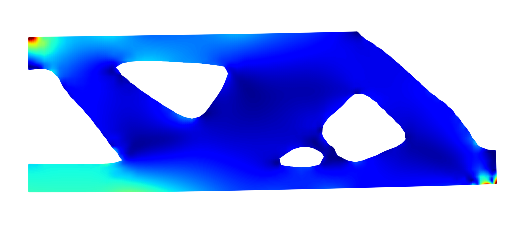}
   \caption{Step 45\label{fig:step045}}
  \end{subfigure}
  \begin{subfigure}{.99\textwidth}
   \centering
   \includegraphics[height=2.5cm]{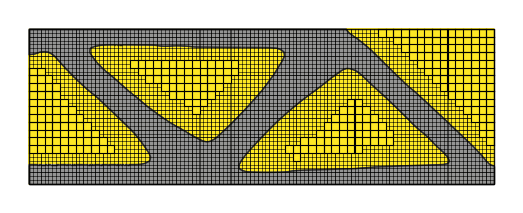}
   \hspace{.5cm}
   \includegraphics[height=2.5cm]{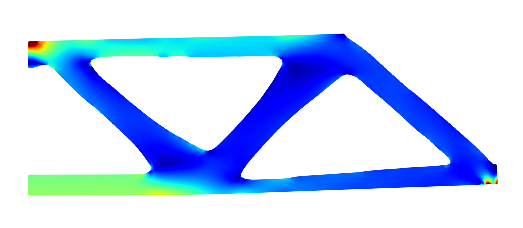}
   \caption{Step 70\label{fig:step070}}
  \end{subfigure}
  \begin{subfigure}{.99\textwidth}
   \centering
   \includegraphics[height=2.5cm]{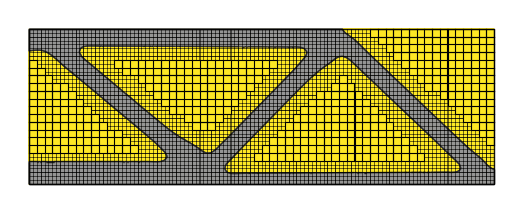}
   \hspace{.5cm}
   \includegraphics[height=2.5cm]{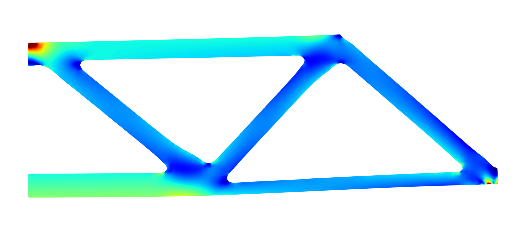}
   \caption{Step 120\label{fig:step120}}
  \end{subfigure}
  \begin{subfigure}{.99\textwidth}
   \centering
   \includegraphics[height=2.5cm]{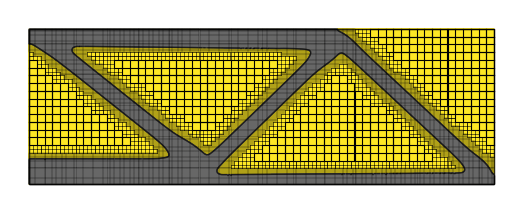}
   \hspace{.5cm}
   \includegraphics[height=2.5cm]{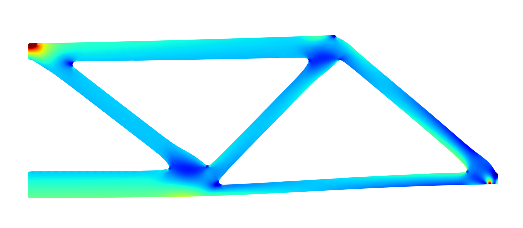}
   \caption{Step 150\label{fig:step150}}
  \end{subfigure}
  \begin{subfigure}{.99\textwidth}
   \hspace{8.3cm}\includegraphics[width=5.3cm]{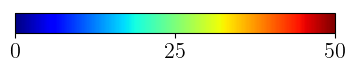}
   \caption*{\hspace{7cm}Frobenius norm of the stress tensor}
  \end{subfigure}
 \caption{Samples of grids, geometries and solutions during the optimization procedure.\label{fig:topologySamples}}
\end{figure}

In this example we apply the developed multigrid preconditioning technique to a dimensionless level set based topology optimization problem, which is inspired by the classical MBB beam \cite{Olhoff1991}. This test case is of particular interest, because it demonstrates the robustness to evolving geometries, and it establishes the suitability to locally refined grids with truncated hierarchical B-splines. Because this involves a large number of computations on an evolving geometry, this numerical experiment is performed in two dimensions to reduce the computational cost. The design space and boundary conditions of the optimization problem are shown in Figure~\ref{fig:topologyProblem}, and samples of the grids and geometries during the procedure are presented in Figure~\ref{fig:topologySamples}. The objective of the design problem is to minimize the strain energy of the structure, subject to a volume constraint that restricts the volume to 30\% of the design domain. The level set function is discretized by linear basis functions which are smoothed by a linear filter \cite{Kreissl2012}. To mitigate the dependence of the optimization results on the initial level set function, a hole seeding method is used that considers the co-evolution of a density field \cite{De2019}. The parameters of the discretized level set and density fields are treated as optimization variables and are updated in the optimization process by the globally convergent method of moving asymptotes (GCMMA) \cite{Svanberg2002}. This test case is intended to establish that elasticity problems on these geometries and grids can be robustly solved in an iterative manner by means of the developed multigrid preconditioner. This opens the doors to level set based topology optimization problems beyond the reach of direct solvers.

On the locally refined grids supplied by the topology optimization procedure, quadratic truncated hierarchical B-splines are constructed, which are trimmed with the provided level set functions.
%For the first $50$ steps, a uniform grid of $60\times20$ elements is applied. The grids in steps $50$ to $124$ contain a single-level hierarchical refinement, which is updated based on the evolution of the physical domain after $75$ and $100$ steps. The grid in steps $125$ until $149$ employs two levels of local refinements. The two-level locally refined grid is updated before the final level set function in step $150$, after which the optimization procedure had reached convergence and was terminated.
In steps 0 to 49 of the routine, a uniform grid of $60\times20$ elements is applied. In steps 50 to 124, grids with a single level of hierarchical refinements are used, see Figure~\ref{fig:topologySamples}. The refinement region is determined based on the design at iterations 49 (i.e., the final design on the uniform grid), 74, and 99. In steps 125 to 149, a grid with two levels of hierarchical refinements is employed. This two-level refinement is based on the final design with a single level of refinements in step 124. The two-level locally refined mesh is updated in step 150 based on the design in the previous step. After this mesh update, the topology optimization algorithm reaches convergence and terminates. The final design is displayed in Figure~\ref{fig:step150}.

For each level set and each mesh occurring in the topology-optimization procedure, we consider a linear elasticity problem with Lam\'{e} parameters $\lambda = \mu = 10^3$. The Dirichlet conditions are weakly imposed by the penalty method, using the penalty parameters $\beta_h^\lambda=\beta_h^\mu=4\frac2h$, with $h=\frac1{20}$ denoting the unrefined element size and the factor $4$ to cover for the local refinements. The vertical support is applied over the width of one unrefined element, and the load at the left top is applied over the same width and has a magnitude of $20$. We apply the conjugate gradient solver, preconditioned by the multigrid cycle with the tailored multiplicative Schwarz smoother. Three levels are applied in the preconditioner, such that the direct solver is applied to a coarsest grid of $15\times5$ elements, with a similar pattern of $0$, $1$, or $2$ levels of local refinements as supplied by the optimization procedure for that step.

\begin{figure}[pt]
  \centering
  \begin{subfigure}{.49\textwidth}
   \centering
   \includegraphics[height=3.5cm]{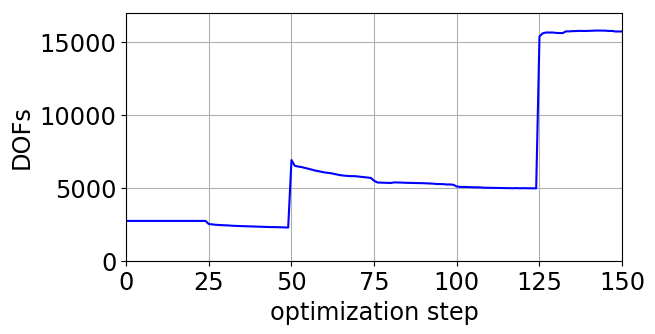}
   \caption{Number of DOFs\label{fig:topologyDOFs}}
  \end{subfigure}
  \hfill
  \begin{subfigure}{.49\textwidth}
   \centering
   \includegraphics[height=3.5cm]{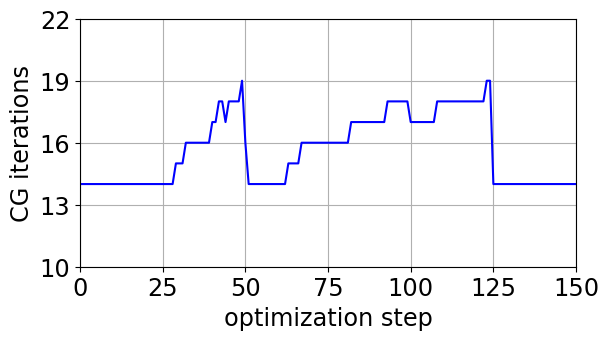}
   \caption{Number of CG iterations\label{fig:topologyIterations}}
  \end{subfigure}
  \caption{Number of DOFs and required number of CG iterations to reduce the residual by $10^{-10}$ during the different steps of the topology optimization procedure.\label{fig:topologyResults}}
\end{figure}

Figure~\ref{fig:topologyResults} plots the number of DOFs and the number of CG iterations during the procedure. The number of DOFs clearly shows a sharp increase after $50$ and $125$ iterations, when the number of local refinement levels is increased. The number of iterations in Figure~\ref{fig:topologyIterations} demonstrates that the preconditioning technique is robust to cut elements, and is not sensitive to changes in the geometry and topology. It can be observed that the number of iterations moderately increases from steps $25$ and $50$ and between steps $50$ and $125$, as the complexity of the evolving physical domain increases. Furthermore, the number of iterations is reduced when more levels of local refinements are applied. These effects are similar to those observed in the three-dimensional test cases; since the extra levels of local refinements are also applied to the coarser grids, the coarse grid correction terms resolve the smooth eigenmodes of the finest grid more accurately, which clarifies the reduction in the number of iterations.

\section{Conclusion}\label{sec:conclusionMulti}

This contribution develops a geometric multigrid preconditioner that enables iterative solutions for higher-order immersed finite element methods at a computational cost that is linear with the number of degrees of freedom. This preconditioning technique is robust to the cut elements in immersed methods, and is applicable to traditional, isogeometric, and locally refined discretizations. This is an improvement with respect to state-of-the-art preconditioning techniques for immersed finite element methods and immersed isogeometric analysis, as these are either restricted to linear discretizations \cite{Lang2014,Lehrenfeld2017,Badia2017} or provide convergence rates that are dependent on the grid size \cite{SIPIC,Prenter2019,John}.

A spectral analysis of immersed finite element methods reveals a spectrum that consists of a combination of \emph{i)} generic modes that can also be observed in mesh-fitting methods and \emph{ii)} modes that are characteristic to immersed finite elements and are only supported on small cut elements. Furthermore, immersed systems that are treated by a Schwarz-type preconditioner -- with blocks selected such that the linear dependencies on small cut elements are resolved -- possess essentially identical spectral properties as mesh-fitting methods. This confirms the observations with additive Schwarz preconditioners in \cite{Prenter2019}, and opens the door to the application of the established and highly efficient framework of multigrid techniques in immersed formulations.

The presented examples convey that the developed multigrid preconditioner results in a number of iterations that is \emph{i)} independent of the mesh size, \emph{ii)} only very slightly affected by the number of levels in the multigrid cycle, and \emph{iii)} only marginally affected by the geometry of the problem. These observations indicate that, while the examples in this contribution already contain multi-million degrees of freedom, the preconditioned iterative solution method enables large-scale computations with immersed finite element methods. While the presented results were obtained with a sequential implementation, further upscaling of the number of degrees of freedom requires an efficient parallel implementation, for which the procedure is suitable.% The procedure is suitable for this, however, and can even be simplified by applying Gauss-Seidel at the coarser levels.

The numerical results are obtained with a baseline multigrid algorithm and relatively straightforward selection of the Schwarz blocks. While the scaling of the computational cost is already optimal with respect to the system size, the framework of multigrid preconditioners with Schwarz-type smoothers allows for adjustments that can enhance the efficiency even further. First of all, as observed in \cite{BeiraoDaVeiga2012,BeiraoDaVeiga2013Schwarz,delaRiva2018}, the efficiency of Schwarz-type smoothers in isogeometric methods is sensitive to the size and the overlap in the selection of the blocks, and different block selections have not been investigated in this contribution. Second, smoothing with multiplicative Schwarz is computationally more expensive than smoothing with additive Schwarz. While multiplicative Schwarz generally results in superior spectral properties and fewer iterations, a detailed comparison between these in the context of the overall computational cost has not been performed. Finally, the V-cycle with a single pre-smoothing and post-smoothing operation is the simplest symmetric multigrid cycle, and it is anticipated that improvements of the cycle design are possible.

Another recommendation pertains to the robustness with respect to the geometrical complexity. Different grid sizes can yield very different stiffness properties when coarse grid basis functions cover the gap between disconnected parts of the geometry. This reduces the effectivity of the multigrid method. Therefore, the robustness can be enhanced by a dedicated coarsening algorithm that identifies coarse grid basis functions with disjoint supports and resolves this with a local refinement or an XFEM-type enrichment.

This contribution only considers symmetric positive definite problems. Based on the investigation of the conditioning problems for different partial differential equations in \cite{Prenter2019}, this is representative for the specific cut-element-related conditioning problems in immersed finite element methods and immersed isogeometric analysis. While symmetric positive definite problems cover a large variety of problems in computational mechanics, the developed preconditioning technique is not immediately applicable to nonsymmetric and mixed formulations in, in particular, flow problems. However, multigrid methods are commonly applied in mesh-fitting flow problems, see e.g., \cite{Brandt2011}, and have been observed to be effective with very similar Schwarz (or Vanka \cite{Vanka1986}) blocks in \cite{Coley2018}. Furthermore, it is demonstrated in \cite{Prenter2019} that Schwarz-type methods can effectively resolve the cut-element-specific conditioning problems in immersed flow problems. Therefore, it is anticipated that the developed multigrid preconditioner extends mutatis mutandis to problems that are not symmetric positive definite.

\section*{Acknowledgement}
\noindent The research of F.\ de Prenter was funded by NWO under the Graduate Program Fluid \& Solid Mechanics. The research of J.A.\ Evans, C.\ Messe, and K.\ Maute was funded by the Defense Advanced Research Projects Agency under Grant Number HR0011-17-2-0022. The simulations in this work were performed using the open source software package Nutils (\cite{nutils2018}, www.nutils.org).

\section*{References}
\bibliographystyle{elsarticle-num}
\bibliography{RefFile}

\end{document}